\documentclass[12pt]{article}
\usepackage[colorlinks,hyperfootnotes=false,pageanchor=false]{hyperref}
\hypersetup{colorlinks=true, linkcolor=black, citecolor=black}
\usepackage{amsmath, amsthm}
\usepackage{algorithm}
\usepackage[noend]{algpseudocode}
\usepackage[nameinlink,noabbrev]{cleveref}
\usepackage[doublespacing]{setspace}
\usepackage{tabularx}
\usepackage{tabulary}
\usepackage{placeins}
\usepackage{graphicx}
\usepackage{pgfplots}
\usepackage{tikz}
\pgfplotsset{compat=newest}
\pgfplotsset{plot coordinates/math parser=false}
\newlength\figureheight
\newlength\figurewidth
\usepackage{relsize}
\usepackage{pdflscape}
\usepackage{amssymb}
\usepackage[round]{natbib}
\usepackage{amsfonts}
\usepackage{amsmath}
\usepackage{bm}
\usepackage{latexsym}
\usepackage{epsfig}
\usepackage{rotating}
\usepackage{ctable}
\usepackage{array}
\usepackage{tabularx}
\usepackage{tabulary}
\usepackage{footnote}
\usepackage{threeparttable}
\usepackage{float}

\newcommand{\blind}{0}
\def\1{1\!{\rm l}}
\newtheorem{theorem}{Theorem}

\newtheorem{proposition}{Proposition}

\newtheorem{corollary}{Corollary}

\newtheorem{remark}{Remark}

\newcommand{\ignore}[1]{}

\newcommand{\y}{\mathbf{y}}
\newcommand{\z}{\mathbf{z}}
\newcommand{\GG}[1]{}

\begin{document}
\def\spacingset#1{\renewcommand{\baselinestretch}%
{#1}\small\normalsize} \spacingset{1}

\if0\blind
{
  \title{\bf Model Misspecification in ABC: Consequences and Diagnostics.\thanks{We would like to thank the Editor, David Dunson, an Associate Editor and two anonymous referees for their constructive comments that greatly improved the paper.}}
  \author{David T. Frazier\footnote{Monash University, Melbourne Australia. Email: david.frazier$@$monash.edu.}, Christian P. Robert\footnote{Universit\'e Paris Dauphine PSL, CEREMADE CNRS, Paris, France. Email: xian@ceremade.dauphine.fr}, and Judith Rousseau\footnote{University of Oxford, Universit\'e Paris Dauphine PSL, CEREMADE CNRS, Paris, France. Email: rousseau@ceremade.dauphine.fr} }
  \maketitle
} \fi

\if1\blind
{\title{\bf}
  
  \bigskip
  \begin{center}
    {\LARGE\bf{}}
\end{center}
  \medskip
} \fi

\begin{abstract}
We analyze the behavior of approximate Bayesian computation (ABC) when the model generating the simulated data differs from the actual data generating process; i.e., when the data simulator in ABC is misspecified. We demonstrate both theoretically and in simple, but practically relevant, examples that when the model is misspecified different versions of ABC can yield substantially different results. Our theoretical results demonstrate that even though the model is misspecified, under regularity conditions, the accept/reject ABC approach concentrates posterior mass on an appropriately defined pseudo-true parameter value. However, under model misspecification the ABC posterior does not yield credible sets with valid frequentist coverage and has non-standard asymptotic behavior. In addition, we examine the theoretical behavior of the popular local regression adjustment to ABC under model misspecification and demonstrate that this approach concentrates posterior mass on a completely different pseudo-true value than accept/reject ABC. Using our theoretical results, we suggest two approaches to diagnose model misspecification in ABC. All theoretical results and diagnostics are illustrated in a simple running example.
\end{abstract}

\vspace{2cm}


\newpage

\section{Introduction}
It is now routine in the astronomic, ecological and genetic sciences, as well as in economics and finance, that the models used to describe observed data are so complex that the likelihoods associated with these models can be computationally intractable. In a Bayesian inference paradigm, these situations have led to the rise of approximate Bayesian computation (ABC) methods that eschew calculation of the likelihood in favor of simulation; for reviews on ABC methods see, e.g., \citet{marinea2012}, \citet{Robert2016} and \cite{sisson2018handbook}.

ABC is predicated on the belief that the observed data $\mathbf{y}:=(y_{1},y_{2},...,y_{n})^{\intercal}$
 is drawn from the class of models $\{{\theta \in
\Theta }:P^n_{{\theta }}\}$, where $\theta\in\Theta\subset\mathbb{R}^{k_{\theta}}$ is an unknown vector of parameters and where $\pi(\theta)$ describes our prior beliefs about $\theta$. The goal of ABC is to conduct inference on the unknown $\theta$ by simulating pseudo-data $\mathbf{z}$, $\mathbf{z}:=(z_{1},...,z_{n})^{\intercal}$, from $P^n_\theta$ and then ``comparing'' $\mathbf{y}$ and $\mathbf{z}$. In most cases, this comparison is carried out using a vector of summary statistics $\eta(\cdot)$, and a metric $d(\cdot,\cdot)$. Generally speaking, in ABC values of $\theta$ are accepted, and used to build an approximation to the exact posterior, if they satisfy an acceptance rule that depends on the tolerance parameter $\epsilon$. 
\begin{algorithm}
\caption{ABC Algorithm}\label{ABC}
\begin{algorithmic}[1]
\State Simulate ${\theta }^{i}$, $i=1,2,...,N$, from $\pi({\theta }),$
\State Simulate $\mathbf{z}^{i}=(z_{1}^{i},z_{2}^{i},...,z_{n}^{i})^{\intercal }$, $i=1,2,...,N$, from $P^n_{\theta^i}$;
\State For each $i=1,...,N$, accept ${\theta }^{i}$  if $d(\eta(\mathbf{z}^{i}),\eta(\mathbf{y}))\leq \epsilon$, where $\epsilon$ denotes an user chosen tolerance parameter $\epsilon$. 
\end{algorithmic}
\end{algorithm}

Algorithm \ref{ABC} details the common accept/reject implementation of ABC, which can be augmented with additional steps to increase sampling efficiency; see, e.g., the MCMC-ABC approach of \citet{Marjoram2003},
 or the SMC-ABC approach of \citet{Sisson2007}. Post-processing of the simulated pairs $\{\theta^{i},\eta(\mathbf{z}^{i})\}$ has also been proposed as a means of obtaining more accurate posterior approximations (see, e.g., the local linear regression adjustment approach of \citealp{Beaumont2025}, the marginal adjustment approach of \citealp{nott2014approximate}, or the recalibration approach of \citealp{rodrigues2018recalibration}). 
	
{While several post processing strategies exist, the most common approach is the so-called local linear regression adjustment (\citealp{Beaumont2025}), which involves post-processing the output from Algorithm \ref{ABC} using a linear regression model to improve the resulting posterior approximation; we refer the interested reader to \citet{blum2018regression} for an overview of regression adjustment methods in ABC. For $\{\theta^{i},\eta(\mathbf{z}^{i})\}_{i\geq1}$ denoting a sample from the ABC posterior based on Algorithm \ref{ABC}, the local linear regression adjustment uses the sample $\{\theta^{i},\eta(\mathbf{z}^{i})\}_{i\geq1}$ to produce the adjusted posterior sample $\{\theta^i-\hat{\beta}^{\intercal}[\eta(\mathbf{z})-\eta(\mathbf{y})]\}_{i\geq1}$, where  $\hat{\beta}$ is obtained from a regression of $\theta^i$ on $\{\eta(\mathbf{z}^i)-\eta(\mathbf{y})\}$.}

Regardless of the ABC algorithm chosen, the very nature of ABC is such that the researcher must believe there are values of $\theta$ in the prior support that can yield simulated summaries $\eta(\mathbf{z})$ that are `close to' the observed summaries $\eta(\mathbf{y})$. Therefore, in order for ABC to yield meaningful inference about $\theta$ there must exist values of $\theta\in\Theta$ such that $\eta(\mathbf{z})$ and $\eta(\mathbf{y})$ are similar.

While complex models allow us to explain many features of the observed data, {it is unlikely that any researcher will be able to construct a model $P^n_{\theta}$ that perfectly reproduces all features of $\mathbf{y}$.} In other words, by the very nature of the complex models to which ABC is applied, the class of models $\{\theta\in\Theta:P^n_{\theta}\}$ used to simulate pseudo-data $\mathbf{z}$ is likely misspecified. Even when accounting for the use of summary statistics that are not sufficient, and which might be compatible with several models, the value these summaries take for the observed data may well be {incompatible, i.e., highly unlikely,} with the realised values of these statistics in the assumed model.

Given the likelihood of model misspecification in empirical applications, understanding the behavior of popular ABC approaches under model misspecification, and the consequences of this behavior, is of paramount importance for practitioners. As the following example illustrates, a particular consequence of model misspecification is that {different ABC approaches can yield significantly different results.}

\medskip

\noindent \textbf{Example 1:} Consider an artificially simple example where the assumed data generating process (DGP) is $z_1,\dots,z_n$ iid as $\mathcal{N}(\theta,1)$ but the actual DGP is $y_1,\dots,y_n$ iid as $\mathcal{N}(\theta,{\sigma}^2)$. That is, for ${\sigma}^2\neq1$, the assumed DGP maintains an incorrect assumption about the variance of the observed data. We consider as the basis of our ABC analysis the following summary statistics:
\begin{itemize}
\item the sample mean $\eta_{1}(\mathbf{y})=\frac{1}{n}\sum_{i=1}^{n}{y}_{i}$,
\item the sample variance  $\eta_{2}(\mathbf{y})=\frac{1}{n-1}\sum_{i=1}^{n}({y}_{i}-\eta_{1}(\mathbf{y}))^{2}$.
\end{itemize}
Consider conducting inference on $\theta$ based on two versions of ABC: the accept/reject approach (hereafter, ABC-AR), where we take $d(x,y)=\|x-y\|$ to be the Euclidean norm, {and a local linear regression adjustment approach to ABC (hereafter, ABC-Reg). ABC-Reg adjusts the accepted draws from ABC-AR using a weighted linear regression of $\theta^i$ on $\{\eta(\mathbf{z}^i)-\eta(\y)\}$, with weights $ K_\epsilon(\|\eta(\mathbf{z}^i)-\eta(\y)\|)$, where $K_\epsilon(\cdot)$ is a kernel function and where the tolerance $\epsilon$ operates as a bandwidth for the kernel function. Following \cite{Beaumont2025}, we take as the kernel function, $K_\epsilon(t)$, the Epanechnikov kernel:
$K_\epsilon(t)= \1_{t \leq \epsilon}\cdot c \epsilon ^ { - 1 } \left( 1 - ( t / \epsilon ) ^ { 2 } \right) $, where $\1_{t\leq \epsilon}$ denotes the indicator function on the event $t\leq \epsilon$,
 and where $c$ is a normalizing constant. }

{To determine how these two ABC approaches behave under varying levels of model misspecification}, we fix $\theta=1$ and simulate ``observed data'' $\mathbf{y}$ according to different values of ${\sigma}^2$. The sample size across the experiments is taken to be $n=100$. We consider a sequence of simulated data
sets for $\y$ such that each corresponds to a different value of
${\sigma}^{2}$, with ${\sigma}^{2}$ taking values from
${\sigma}^{2}=.5$ to ${\sigma}^{2}=5$ with evenly spaced
increments of $0.05$. Across all the datasets we fix the random numbers used to generate
the simulated observed data and only change the value of ${\sigma}^{2}$ to isolate
the impact of model misspecification; i.e., we generate one common set of
random numbers $\nu_{i}\sim \mathcal{N}(0,1)$, $i=1,...,100$,
then for  a value of ${\sigma}^{2}$ we generate observed data according to
$y_{i}=1+\nu_{i}\cdot{\sigma}$.

{Our prior beliefs are given by $\theta\sim \mathcal{N}(0,25)$. We implement ABC-AR using $N=25, 000$ simulated pseudo-datasets generated iid according to $z^{j}_{i}\sim \mathcal{N}(\theta^j,1)$.  For both ABC-AR and ABC-Reg, we set $\epsilon$ to be the 1\% quantile of the simulated distances $\|\eta(\mathbf{y})-\eta(\mathbf{z}^{j})\|$. To further isolate the impact of randomness on this procedure, we use the same simulated data across the different observed datasets; i.e., both ABC procedures will use the same simulated data across the different values of $\sigma^2$. By fixing the simulated data across the experiments, and by fixing the random numbers in the observed data, differences in the ABC output across the experiments can be attributed to the changing value of $\sigma^2$.}

Figure \ref{fig1} compares the posterior means of ABC-AR, and ABC-Reg across the different values for ${\sigma}^{2}$.\footnote{The posterior densities for ABC-AR and ABC-Reg obtained from these experiments also display a similar pattern of behavior and are presented in the appendix.}  The results demonstrate that model misspecification {induces dramatic differences between the two ABC approaches}, even at a relatively small sample size.\footnote{Even though the DGP for $\z$ is misspecified, because of the nature of the model misspecification and the limiting behavior of $\eta(\y)$, if one were to only use the first summary statistic (the sample mean) model misspecification would have little impact in this example. However, in general both the nature of the model misspecification and the precise limiting form of $\eta(\y)$ are unknown. Therefore, choosing a set of summaries that can mitigate the impact of model misspecification will be difficult, if not impossible, in practical applications of ABC. } We draw two specific conclusions from Figure \ref{fig1}: one, {the posterior mean of ABC-AR remains relatively stable across the different levels of misspecification, but does shift away from the true value ($\theta=1$) as the level of model misspecification increases}; two, the posterior mean of ABC-Reg becomes unstable even at relatively small levels of misspecification. The performance of the local linear regression adjustment is particularly interesting given that this method can have theoretical advantages over ABC-AR, i.e., Algorithm \ref{ABC}, when the model is correctly specified (\citealp{LF2016b}). We formally explore these issues in Sections two and three but note here that when ${\sigma}^2\approx1$  (i.e., correct model specification) both ABC approaches give similar results.
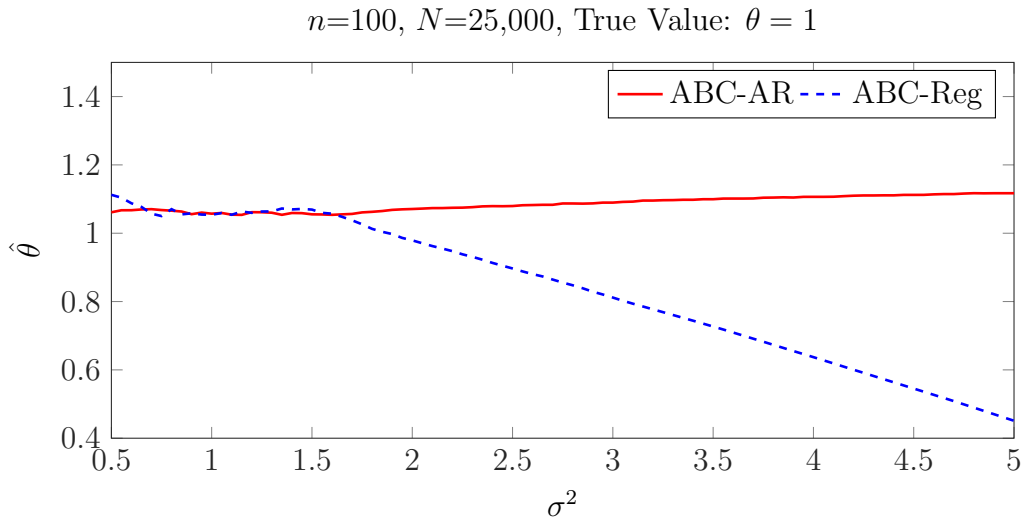
\begin{figure}[h!]
\centering 
\setlength\figureheight{5cm} 
\setlength\figurewidth{12cm} 
%
%
%
%
\begin{tikzpicture}

\begin{axis}[%
view={0}{90},
width=\figurewidth,
height=\figureheight,
scale only axis,
every outer x axis line/.append style={darkgray!60!black},
every x tick label/.append style={font=\color{darkgray!60!black}},
xmin=0.5, xmax=5,
xlabel={$\sigma{}^2$},
every outer y axis line/.append style={darkgray!60!black},
every y tick label/.append style={font=\color{darkgray!60!black}},
ymin=0.4, ymax=1.5,
ylabel={$\hat{\theta}$},
title={$n$=100, $N$=25,000, True Value: $\theta=1$ },
legend style={legend columns=3,align=left}]
\addplot [
color=red,
solid,
line width=1.0pt,
]
coordinates{
 (0.5,1.06126396)(0.55,1.06721446215139)(0.6,1.06734432)(0.65,1.06956744)(0.7,1.0704132)(0.75,1.06783988)(0.8,1.066057)(0.85,1.0635108)(0.9,1.05575216)(0.95,1.060392)(1,1.05722232)(1.05,1.05926524)(1.1,1.05428308)(1.15,1.05371492)(1.2,1.06161904)(1.25,1.06136492)(1.3,1.05999764)(1.35,1.05390792)(1.4,1.05918704)(1.45,1.05874816)(1.5,1.05529804)(1.55,1.05493172)(1.6,1.05383616)(1.65,1.05522176)(1.7,1.05650737051793)(1.75,1.06095964)(1.8,1.06265436)(1.85,1.06495368)(1.9,1.06834824)(1.95,1.06994336)(2,1.07089124)(2.05,1.07211292)(2.1,1.07352532)(2.15,1.07361612)(2.2,1.07448636)(2.25,1.07485404)(2.3,1.0762828)(2.35,1.07857204)(2.4,1.07930364)(2.45,1.07904244)(2.5,1.07957404)(2.55,1.0819014)(2.6,1.0824954)(2.65,1.08320256)(2.7,1.08320256)(2.75,1.08711087649402)(2.8,1.08706132)(2.85,1.08661944)(2.9,1.08734584)(2.95,1.08960372)(3,1.08960372)(3.05,1.0916736)(3.1,1.0924996)(3.15,1.09558424)(3.2,1.09558424)(3.25,1.09665688)(3.3,1.09694124)(3.35,1.09792628)(3.4,1.09792628)(3.45,1.0995232)(3.5,1.0995232)(3.55,1.10139008)(3.6,1.10147008)(3.65,1.10147008)(3.7,1.10176528)(3.75,1.10396908)(3.8,1.10435668)(3.85,1.10509468)(3.9,1.1045002)(3.95,1.10645856)(4,1.10652416)(4.05,1.10652416)(4.1,1.1068288)(4.15,1.1084646)(4.2,1.1098198)(4.25,1.110549)(4.3,1.110549)(4.35,1.1109122)(4.4,1.1108586)(4.45,1.1122046)(4.5,1.1122046)(4.55,1.1122874501992)(4.6,1.11375328)(4.65,1.11457492)(4.7,1.11457492)(4.75,1.1161036)(4.8,1.1171588)(4.85,1.11680724)(4.9,1.11708328)(4.95,1.11708328)(5,1.11708328) 
};
\addlegendentry{ABC-AR};

\addplot [
color=blue,
dashed,
line width=1.0pt,
]
coordinates{
 (0.5,1.1124192994594)(0.55,1.10403871113261)(0.6,1.0876694680244)(0.65,1.07702251515685)(0.7,1.05719524397423)(0.75,1.05026431951203)(0.8,1.07045013489805)(0.85,1.05513798551897)(0.9,1.05826067305938)(0.95,1.05455071974121)(1,1.05407409144824)(1.05,1.05908057819918)(1.1,1.05404021650992)(1.15,1.06350677603022)(1.2,1.05967526813763)(1.25,1.06349030723327)(1.3,1.06390284696284)(1.35,1.07235265997687)(1.4,1.06964969352786)(1.45,1.07110081165746)(1.5,1.06917894445536)(1.55,1.05960850536597)(1.6,1.0572814689328)(1.65,1.04800589005071)(1.7,1.03778141467116)(1.75,1.02616911090976)(1.8,1.01271338075073)(1.85,1.00423058104648)(1.9,0.99755416213836)(1.95,0.986241363983929)(2,0.979167460274823)(2.05,0.970640497572598)(2.1,0.962426475078956)(2.15,0.955156563775748)(2.2,0.94755072345351)(2.25,0.938307467910352)(2.3,0.930827131697382)(2.35,0.922527506840915)(2.4,0.913126621123698)(2.45,0.90491851176993)(2.5,0.896831256913158)(2.55,0.888232735522674)(2.6,0.880349502992181)(2.65,0.873472842862642)(2.7,0.86438715783062)(2.75,0.855591536302843)(2.8,0.847616904906305)(2.85,0.838851746495387)(2.9,0.828866059739685)(2.95,0.820700074208587)(3,0.811697781064918)(3.05,0.802390286182524)(3.1,0.793782632458642)(3.15,0.785556810032348)(3.2,0.777218677889794)(3.25,0.768973070116992)(3.3,0.760614446355674)(3.35,0.752246084054683)(3.4,0.743981151581867)(3.45,0.735427946862728)(3.5,0.726845247210867)(3.55,0.718266738221852)(3.6,0.709288551150594)(3.65,0.700303187750083)(3.7,0.691310984613284)(3.75,0.682458056509052)(3.8,0.673309164022934)(3.85,0.664256081664859)(3.9,0.655196931447296)(3.95,0.6462208510719)(4,0.637281644225029)(4.05,0.628062841462607)(4.1,0.618924930305093)(4.15,0.609683971649881)(4.2,0.600442388549977)(4.25,0.591302623090779)(4.3,0.582060305816915)(4.35,0.572873077606975)(4.4,0.563686406418195)(4.45,0.554500260163553)(4.5,0.545335321625686)(4.55,0.53615013789857)(4.6,0.526965397115943)(4.65,0.517781075884073)(4.7,0.508617801036804)(4.75,0.499007039052181)(4.8,0.489436212184759)(4.85,0.479865907383193)(4.9,0.470296097630052)(4.95,0.460623717519281)(5,0.451076349600749) 
};\addlegendentry{ABC-Reg};
\end{axis}
\end{tikzpicture}%
\caption{Comparison of posterior means, denoted by $\hat{\theta}$, for ABC-AR, and ABC-Reg across varying levels of model misspecification. } 
\label{fig1} 
\end{figure}
\hfill$\square$

In the remainder of the paper, we elaborate on the above issues and rigorously characterizes the asymptotic behavior {of accept/reject ABC and  local regression adjustment ABC} when the model generating the pseudo-data is misspecified. In Section two, we discuss model misspecification in the ABC context and demonstrate that under model misspecification, for a certain choice of the tolerance, the posterior associated with Algorithm \ref{ABC} asymptotically concentrates on an appropriately defined pseudo-true value. In addition, we demonstrate that the asymptotic shape of the ABC posterior is non-standard under model misspecification, and can yield credible sets with arbitrary levels of coverage. In Section three, we provide a rigorous justification for the behavior observed in Figure \ref{fig1}: we demonstrate that under model misspecification the posteriors for local linear (and nonlinear) regression adjustment ABC can asymptotically concentrate onto a completely different region of the parameter space than accept/reject ABC (Algorithm \ref{ABC}). We then use these theoretical results to devise an alternative regression adjustment approach that performs well regardless of model specification. Motivated by our asymptotic results, in Section four we develop two model misspecification detection procedures: a graphical detection approach based on comparing acceptance probabilities from Algorithm \ref{ABC} and an approach based on comparing the output from Algorithm \ref{ABC} and its linear regression adjustment counterpart. {Section five concludes with a brief overview and discussion of our results}. Proofs of all theoretical results are contained in the supplementary appendix.

\section{Model Misspecification in ABC\label{Aux}}
Before rigorously characterising the behavior of ABC under model misspecification, we first set the notation used throughout the remainder of the paper. For $\mathbf{y}$ denoting the observed data, we let $P^n_{0}$ denote the true distribution generating $\mathbf{y}$. The class of implied distributions used in ABC to simulate pseudo-data is denoted by $\mathcal{P}:=\{\theta\in\Theta\subseteq\mathbb{R}^{k_{\theta}}:P^n_{\theta}\}$, while $\mathbf{z}$ denotes pseudo-data with support $\mathcal{Z}$ generated according to $P_\theta^n$. The simulated summary statistics $\eta(\mathbf{z})=(\eta_{1}(\mathbf{z}),...,\eta_{k_{\eta}}(\mathbf{z}))^\intercal$ are a $k_{\eta}$-dimensional random vector with support $\mathcal{B}:=\{
{\eta }(\mathbf{z}):\mathbf{z}\in \mathcal{Z}\}\subseteq\mathbb{R}^{k_{\eta}}$. We let $d_{1}(\cdot,\cdot)$ denote a metric on ${\Theta }$, and $d_2(\cdot ,\cdot)$ a metric on $\mathcal{B}$. However, when no confusion will result we simply denote a generic metric by $d(\cdot,\cdot)$. $\Pi ({\theta })$ denotes the prior measure and $\pi({\theta })$ its corresponding density. {For economy of notation, in what follows, we disregard the dependence of $P^n_0$ and $P_\theta^n$ on $n$, and simply represent these quanties as $P_0$ and $P_\theta$, respectively. } 

\subsection{On the Notion of  Model Misspecification in ABC}\label{sec:def:mis}
Recall that, in likelihood-based inference, model misspecification is take to mean that  ${P}_{0}\notin \mathcal{P}$. The result of this model misspecification is that the Kullback-Leibler divergence, $$\mathcal{D}(P_{0}||P_{\theta})=-\int_{}\log\left\{\frac{dP_{\theta}(\mathbf{y})}{dP_{0}(\mathbf{y})}\right\}dP_{0}(\mathbf{y}),$$ satisfies $$\inf_{\theta\in\Theta}\mathcal{D}(P_{0}||P_{\theta})>0.$$ In this case, the parameter value $$\theta^{*}=\arg\inf_{\theta\in\Theta}\mathcal{D}(P_{0}||P_{\theta})$$ is referred to as the pseudo-true value. Even though the model is misspecified, under reasonable regularity conditions, Bayesian procedures predicated on the likelihood of $P_{\theta}$ yield posteriors that concentrate on $\theta^*$; see, e.g., \citet{kleijn2012} and \citet{Muller2013}.

In this paper, we assume the researcher conducts posterior inference on $\theta$ via ABC when the observed sample $\y$ is generated according to $P_{0}$, and in the case where $P_{0}\notin \mathcal{P}$. However, in contrast to likelihood-based procedures, ABC is not based on the full data $\mathbf y$ but on two separate approximations, the summary statistics $\eta(\mathbf y) $ and the threshold $\epsilon$. Therefore, even if $P_{0}\notin \mathcal{P}$ the model class $\mathcal{P}$ may still be capable of generating a simulated summary $\eta(\mathbf z) $ that is compatible with the observed summary $\eta(\mathbf y)$, or is within an $\epsilon$-neighbourhood of $\eta(\y)$. The approximate nature of ABC means that the notion of model misspecification in a likelihood-based setting, namely $\inf_{\theta}\mathcal{D}(P_{0}||P_{\theta})>0$, is not necessarily a meaningful notion of model misspecification associated with the output of an ABC algorithm, or ABC posterior distributions.

Recalling that the ABC posterior measure is given by, for $A\subset {\Theta }$,
\begin{flalign*}
\Pi_{\epsilon}[ A|\eta(\mathbf{y})] &=
{\int_{A}P_{\theta}\left[d(\eta(\mathbf{y}),\eta(\mathbf{z}))\leq \epsilon_{}\right]d\Pi(\theta)}\bigg{/}{\int_{\Theta}P_{\theta}\left[d(\eta(\mathbf{y}),\eta(\mathbf{z}))\leq \epsilon\right]d\Pi(\theta)},
\end{flalign*}
we see that misspecification in ABC will be driven by the behavior of $\eta(\mathbf{y}),\eta(\mathbf{z})$ and the set $\{\theta\in\Theta:d(\eta(\mathbf{y}),\eta(\mathbf{z}))\leq\epsilon\}$. 
To rigorously formulate the notion of model misspecification in ABC,  we must therefore study the limiting behaviour of the ABC likelihood $P_{\theta}\left[d(\eta(\mathbf{y}),\eta(\mathbf{z}))\leq \epsilon\right]$ as the amount of information in the data accumulates. 

To this end, we follow the framework of  \citet{Marinetal2014}, \citet{FMRR2016} and \citet{LF2016a,LF2016b}, where it is assumed that the summary statistics concentrate around some fixed value, namely, $\eta(\y)$ to $b_0 $ (under $P_0$) and $\eta(\z)$ to $b(\theta)$ (under $P_\theta$).  In \citet{Marinetal2014}, the authors study the case where $\epsilon = 0$, while \citet{ FMRR2016} and \citet{LF2016a,LF2016b} study $\epsilon >0$ but allow $\epsilon$ to vary with $n$ and set $\epsilon = \epsilon_n$. In the latter papers the authors demonstrate that  the
amount of information ABC obtains about a given $\theta$ depends
on: (1) the rate at which the observed (resp. simulated)
summaries converge to a well-defined limit counterpart ${b}_{0}$ (resp., ${b}({\theta })$); 
 (2) the rate at which the tolerance $\epsilon_{n}$ goes to zero; (3) the link between $b_{0}$ and $b(\theta)$. When $P_0 \in \mathcal P$, there exists some $\theta_0 $ such that $b(\theta_0) = b_0$ and the results of \citet{FMRR2016} completely characterize the asymptotic behaviour of the ABC posterior distribution. Furthermore, this analysis remains correct even if $P_0 \notin \mathcal P$, so long as there exists some $\theta_0 \in \Theta$ such that $b_0 = b(\theta_0)$. 
 
 Therefore, the meaningful concept of model misspecification in ABC is that there does not exist any $\theta_0\in\Theta$ satisfying $b_0=b(\theta_0)$, which is precisely the notion of model incompatibility defined in \citet{Marinetal2014}. Throughout the remainder, we say that the model is (ABC) misspecified if
 \begin{equation}
\epsilon^{*}=\inf_{{\theta\in\Theta}}d(b_{0},b(\theta))>0\label{bad1}
\end{equation} and note here that this condition is more likely to occur when $k_\theta < k_\eta$.

Heuristically, the implication of misspecification in ABC is that, {under concentration of $\eta(\mathbf{z})$ to $b(\theta)$ and $\eta(\y)$ to $b_0$, by the triangle inequality and the definition of $\epsilon^*$}
$$d(\eta(\mathbf{y}),\eta(\mathbf{z}))  \geq d(b_{0},b(\theta)) - o_{P_{\theta}}(1) - o_{P_{0}}(1) \geq \epsilon^* - o_p(1) \text{, for all }\epsilon_{n} = o(1),$$ and the event $\{\theta\in\Theta:d(\eta(\mathbf{y}),\eta(\mathbf{z}))\leq \epsilon_n\}$ becomes extremely rare, and corresponds to the event $$\left\{\theta\in\Theta:d(\eta(\mathbf z), b(\theta))> \epsilon^* - o(1) \right\}.$$ Consequently, for a sequence of tolerances $\epsilon_{n}=o(1)$, once $\epsilon_{n}<\epsilon^{*}+o(1)$ hardly any draws of $\theta$ will be selected regardless of how many simulated samples from $\pi(\theta)$ we generate, and the ABC posterior $\Pi_{\epsilon}[A|{\eta }(\mathbf{y})]$ will become ill-behaved as $n$ increases.

While tolerance sequences $\epsilon_{n}=o(1)$ will eventually cause  the posterior $\Pi_{\epsilon}[A|{\eta }(\mathbf{y})]$ to be ill-behaved, it is possible that other choices for $\epsilon_{n}$ will produce a well-behaved posterior. In the following section we show that (certain) tolerance sequences satisfying $\epsilon_{n}\rightarrow\epsilon^*$, as $n\rightarrow+\infty$, yield well-behaved ABC posteriors that concentrate posterior mass on an appropriately defined pseudo-true value.

\subsection{ABC Posterior Concentration Under Misspecification}\label{const_s}
Building on the intuition in the previous section, in this and the following section we rigorously characterize the asymptotic behaviour of
\begin{flalign*}
\Pi_{\epsilon}[ A|\eta(\mathbf{y})] &=
{\int_{A}P_{\theta}\left[d(\eta(\mathbf{y}),\eta(\mathbf{z}))\leq \epsilon_{n}\right]d\Pi(\theta)}\bigg{/}{\int_{\Theta}P_{\theta}\left[d(\eta(\mathbf{y}),\eta(\mathbf{z}))\leq \epsilon_{n}\right]d\Pi(\theta)}
\end{flalign*}when $P_{0}\notin\mathcal{P}$ and $\epsilon^*>0$. To do so, we first define the following additional notations: for sequences $\{a_{n}\}$ and $\{b_{n}\}$, real valued, $a_{n}\lesssim b_{n}$
denotes $a_{n}\leq Cb_{n}$ for some $C>0$, $a_{n}\asymp b_{n}$ denotes equivalent order of
magnitude, $a_{n}{\gg }b_{n}$ indicates a larger order of magnitude and the
symbols $o_{P}(a_{n}),O_{P}(b_{n})$ have their usual meaning. Unless otherwise noted, all limits are taken as $n\rightarrow+\infty$.

We maintain the following assumptions. 

\medskip 

\noindent \textbf{[A0]} {There exist a unique $b_0$ such that} $d(\eta(\mathbf{y} ), b_{0}) = o_{P_0}(1)$ and a positive sequence $v_{0,n} \rightarrow +\infty $ such that  
$$  \liminf_{n\rightarrow+\infty} P_0\left[ d(\eta(\mathbf{y}), b_{0}) \geq v_{0,n}^{-1} \right]=1.$$

\noindent \textbf{[A1]} There exist a continuous, injective map $b:\Theta\rightarrow\mathcal{B}\subset \mathbb R^{k_{\eta}}$ and a function $\rho_{n}(\cdot)$ satisfying: $\rho _{n}(u)\rightarrow 0$ as $%
n\rightarrow +\infty $ for all $u>0$, and $\rho _{n}(u)$ monotone
non-increasing in $u$ (for any given $n$), such that, for all $\theta \in \Theta$,  
\begin{equation*}
P_{{\theta }}\left[ d(\eta(\z),b(\theta )) >u\right] \leq c({\theta }%
)\rho _{n}(u),\quad \int_{\Theta }c({\theta })d\Pi ({%
\theta })<+\infty,
\end{equation*}%
where $\mathbf{z}\sim P_{{\theta }}$, and we assume
either of the following:

\begin{enumerate}
\item[\textbf{(i)}] \textit{Polynomial deviations:} There exist a positive sequence $v_{n}\rightarrow+\infty$ and $u_{0},\kappa>0$ such that $\rho_{n}(u)=v_{n}^{-\kappa}u^{-\kappa}$, for $u\leq u_{0}$.  
\item[\textbf{(ii)}] \textit{Exponential deviations:} There exists $h_{\theta}(\cdot)>0$ such that $P_{\theta}[d(\eta(\mathbf{z}),b(\theta))>u]\leq c(\theta)e^{-h_{\theta}(u v_{n})}$ and there exists $m,C>0$ such that $$\int_{\Theta} c(\theta)e^{-h_{\theta}(u v_{n})}d\Pi(\theta)\leq Ce^{-m\cdot(u v_{n})^{\tau}},\;\;\text{for }u\leq u_{0}.$$ 
\end{enumerate}


\noindent \textbf{[A2]} There exist some $D>0$ and $M_0, \delta_0>0$ such that, for all $\delta_0\geq \delta >0$ and  $M\geq M_0$,  there exists $S_\delta \subset \left\{\theta \in \Theta : d( b(\theta),b_{0})-\epsilon^* \leq \delta \right\}$ for which 
\begin{itemize}
\item[\textbf{(i)}]  In case (i) of \textbf{[A1]}, $D< \kappa $ and $\int_{S_\delta} \left( 1  - \frac{ c(\theta) }{ M }\right)d\Pi(\theta)  \gtrsim\delta^{D}.
$
\item[\textbf{(ii)}]  In case (ii) of \textbf{[A1]}, $\int_{S_\delta} \left( 1  - c(\theta)e^{-h_\theta( M )} \right) d\Pi(\theta)  \gtrsim\delta^{D}.$
\end{itemize}

{The above assumptions are similar to those given in \cite{FMRR2016}, and we refer the interested reader to Remarks 1 and 2, and Example 1 in that paper for a detailed discussion of these assumptions. Under the above assumptions, we have the following result.  }
\begin{theorem}\label{th:asymp0}
Assume that the data generating process for $\mathbf{y}$ satisfies \textbf{[A0]} and assume that equation \eqref{bad1} holds. 
Assume also that conditions \textbf{[A1]} and \textbf{[A2]} are satisfied and $\epsilon_n \downarrow \epsilon^* $ with 
$$\epsilon_n \geq \epsilon^*  +  M v_n^{-1}  + v_{0,n}^{-1}, $$ for $M$  large  enough. Let $M_n$ be any positive sequence going to infinity and 
$\delta_n \geq  M_n(\epsilon_n - \epsilon^*)$, then  
\begin{equation}\label{post:conc:good}
\Pi_{\epsilon}\left[ d(b(\theta), b_{0}) \geq \epsilon^*+ \delta_{n} |\eta( \mathbf{y} ) \right] = o_{P_0}(1),
\end{equation}
as soon as
 \begin{equation*}
 \begin{split}
\delta_n &\geq M_n v_n^{-1}u_n^{-D/\kappa}  = o(1) \quad \mbox{in case (i) of assumption [A1]}\\
\delta_n&\geq  M_n v_n^{-1}|\log(u_n)|^{1/\tau} = o(1)  \quad \mbox{in case (ii) of assumption [A1] }.
 \end{split}
 \end{equation*}
 with $u_n =\epsilon_n - (\epsilon^*  +  M v_n^{-1}  + v_{0,n}^{-1}  )\geq 0$.
\end{theorem}We remind the reader that the proofs of all theoretical results are contained in the appendix.
\begin{remark}
\normalfont{Theorem \ref{th:asymp0} states that even though the model is misspecified, the ABC posterior concentrates onto $$\arg\inf_{\theta \in \Theta}d(b(\theta),b_0),$$ under the assumption that $\epsilon_n$ is slightly larger than $\epsilon^*$. Under the more precise framework of Theorem \ref{normal_thm}, where the asymptotic shape of the posterior distribution is studied, this condition can be refined to allow $\epsilon_n$ to be slightly smaller than $\epsilon^*$. However, we demonstrate that if $\epsilon^* - \epsilon_n$ is bounded below by a positive constant, then the posterior distribution does not necessarily concentrate. }
\end{remark}
Using the posterior concentration in Theorem \ref{th:asymp0}, we have the following result.
\begin{corollary}\label{cor1}Assume the hypotheses of Theorem \ref{th:asymp0} are satisfied and define $\theta^*\in\Theta$ as $$\theta^*=\arg\inf_{\theta\in\Theta}d(b_{0},b(\theta)), $$ then, for any $\delta>0$,
$$\Pi_{\epsilon}[d(\theta,\theta^*)>\delta|\eta(\mathbf{y})]=o_{P_{0}}(1).$$
\end{corollary}
\begin{remark}\normalfont{Theorem \ref{th:asymp0} and Corollary \ref{cor1} demonstrate that, under an identification condition, the ABC posterior $\Pi_{\epsilon}[\cdot|\eta(\y)]$ concentrates on $\theta^*$ if the model is misspecified. Therefore, Theorem \ref{th:asymp0} is an extension of Theorem 1 in \citet{FMRR2016} to the case of misspecified models. In addition, we note that Theorem \ref{th:asymp0} above is similar to Theorem 4.3 in \cite{Bernton2017} for ABC inference based on the Wasserstein distance. {The validity of each of these results requires that the map $\theta\mapsto b(\theta)$ be injective. If  this condition is not satisfied, there can exist a continuum of values under which $d(b(\theta),b_0)=\epsilon^*$. In this case, the ABC posterior will no longer converge to a point mass, but will concentrate onto the set $\{\theta\in\Theta:\epsilon^*=d(b(\theta),b_0)\}$.}}
\end{remark}

\begin{remark} \normalfont{It is crucial to note that the pseudo-true value $\theta^*$ depends on the choice of $d(\cdot,\cdot)$. This implies that ABC based on two different metrics $d(\cdot,\cdot)$ and $\tilde{d}(\cdot,\cdot)$ will produce two different pseudo-true values, unless if by happenstance $\inf\{\theta\in\Theta:d(b(\theta),b_0)\}$ and $\inf\{\theta\in\Theta:\tilde{d}(b(\theta),b_0)\}$ coincide. This lies in stark contrast to the posterior concentration result in \cite{FMRR2016}, which demonstrated that under correct model specification the posterior $\Pi_{\epsilon}[\cdot|\eta(\y)]$ concentrates on the same true value regardless of the choice of $d(\cdot,\cdot)$. }
\end{remark}

\subsection{Shape of the Asymptotic Posterior Distribution\label{norm}}

In this section, we analyse the asymptotic shape of the ABC posterior under model misspecification. For simplicity, we take the rate at which the simulated and observed summaries converge to their limit counterparts to be the same, i.e., we take $v_{0,n}=v_{n}$,\footnote{{Allowing $v_{0,n}$ and $v_{n}$ to differ will not greatly alter the following result. The result presented will still be valid, but only for the slower of the two rates.}} and we consider as the distance $d(\eta(\z) , \eta(\y) ) = \| \eta(\z)- \eta(\y) \|$ where $\| \cdot \| $ denotes is the norm associated to a given scalar product $\langle \cdot , \cdot \rangle$.  Denote by $I_{k_\eta}$ the $(k_\eta\times k_\eta)$ dimensional identity matrix and let 
 $$ \Phi(B) =\text{ Pr}\left[ \mathcal N(0, I_{k_\eta}) \in B\right], $$
 for any measurable subset $B$ of $\mathbb R^{k_\eta}$. 

The following conditions are needed to establish the results of this section.
 
\medskip

\noindent\textbf{[A0$^{\prime}$]} Assumption \textbf{[A0]} is satisfied, and $\epsilon^*=d(b(\theta^*),b_0)>0$, where $\theta^*=\arg\inf_{\theta\in\Theta}d(b(\theta),b_0)$.

\medskip

\noindent \textbf{[A1$^{\prime }$]} Assumption \textbf{[A1] }holds and for some positive-definite matrix $\Sigma _{n}({\theta }^{*})$, $c_0>0$, $\kappa >1$ and $\delta >0$, for all $\Vert {\theta }-{%
\theta }^{*}\Vert \leq \delta $, $P_{{\theta }}\left[ \Vert \Sigma _{n}({%
\theta }^{*})\{{\eta }(\mathbf{z})-{b}({\theta })\}\Vert >u\right] \leq {c_{0}%
}{u^{-\kappa }}$ for all $0<u\leq \delta v_{n}$.\medskip

\noindent \textbf{[A3]} The
map $\theta\mapsto{b}(\theta)$ is twice continuously differentiable at ${\theta
^{*}}$ and the Jacobian $\nabla _{\theta }b({\theta }^{*})$ has full column
rank $k_{\theta }$. The Hessian of $\| b(\theta) - b_0 \|^2$ evaluated at $\theta^* $, and denoted by $H^*$, is a positive-definite matrix.

\medskip

\noindent \textbf{[A4]} There exists a sequence of $(k_{\eta }\times k_{\eta
}) $ positive-definite matrices ${\Sigma }_{n}({\theta})$ such that for all $M>0$ there exists $u_0 >0$ for which
\begin{equation*}
\sup_{|x| \leq M} \sup_{\|\theta - \theta^*\|\leq u_0 } \left|P_\theta \left(\langle Z_n, e  \rangle \leq x \right)  - \Phi(x) \right| = o(1) ,  
\end{equation*}%
where $Z_n = \Sigma_n(\theta) ( \eta(\mathbf{z} ) - b(\theta) )$ and $e = (b(\theta^*) - b_0 ) /\| b(\theta^*) - b_0\|$. 

\medskip

\noindent \textbf{[A5]} There exist $v_n $ going to infinity  and $u_0>0$ such that for all $\Vert \theta -\theta ^{*}\Vert \leq u_0$%
, the sequence of functions ${\theta }\mapsto {\Sigma }_{n}({\theta }%
)v_{n}^{-1}$ converges to some positive-definite matrix $A(\theta )$ and is
equicontinuous at ${\theta }^{*}$.

\medskip

\noindent \textbf{[A6]} $\pi(\theta)$, the density of the prior measure $\Pi(\theta)$, is continuous and positive  at $\theta^* $. 
\medskip 

\noindent \textbf{[A7]} For $Z_{n}^{0}={\Sigma }_{n}({\theta }^{*})\{{\eta }(\mathbf{y})-{b_{0}}\}$  and all $M_n $ going to infinity 
$$P_0\left( \| Z_{n}^{0}\| > M_n \right) = o(1) .$$ 

\medskip 

The above assumptions are similar to those used in \cite{FMRR2016} to deduce the limiting shape of the ABC posterior under correct model specification, and we refer the interested reader to Remarks 3 and 4 in that paper for a detailed discussion of these assumptions. Under the above assumptions, we have the following result.

\begin{theorem}
\label{normal_thm} Assume \textbf{[A0\,{{$^{\prime }$}}]}, \textbf{[A1}{{$^{\prime }$}}\textbf{]} (with $\kappa
\geq k_{\theta }$),\textbf{\
[A2] }and \textbf{[A3}\textbf{]}-\textbf{[A7]} are satisfied. We then have the following results. 

\begin{itemize}
\item[\textbf{(i)}] If $\lim_{n}v_{n}(\epsilon _{n}-\epsilon^*)= 2c$, with $c\in \mathbb R$, then for $\|\cdot\|_{{TV}}$ the total-variation norm 
\begin{equation*}
\| \Pi_{v_n^{1/2}, \epsilon} - Q_c \|_{{TV}} =  o_{P_0}(1) 
\end{equation*}
where $\Pi_{z_n, \epsilon}$ is the ABC posterior distribution of $z_n(\theta - \theta^*) $ for any sequence $z_n>0$ and $Q_c$ has density $q_c$ with respect to Lebesgue measure on $\mathbb R^{k_\theta}$ proportional to 
$$ q_c(x) \propto \Phi\left(\frac{ c - \langle Z_n^0 , A(\theta^*)e \rangle \epsilon^* }{ \|A(\theta^*) e\| \epsilon^* } - \frac{ x^{\intercal} H^* x }{ 4 \|A(\theta^*) e \| \epsilon^* }\right) $$

\item[\textbf{(ii)}] If $\lim_{n}v_{n}(\epsilon _{n}-\epsilon^*)= + \infty $ with $u_n = \epsilon _{n}-\epsilon^* = o(1) $, for $\mathcal U_{\{\|x\|\leq M\}}$ the uniform measure over the set $\{\|x\|\leq M\}$,
\begin{equation*}
\|  \Pi_{u_n^{-1}, \epsilon} - \mathcal U_{\{x^{\intercal} H^* x \leq 2 \}} \|_{TV} =  o_{P_0}(1) , 
\end{equation*}

\end{itemize}
\end{theorem}

\begin{remark} \label{rem:geneshape}\normalfont{
As is true in the case where the model is correctly specified, if $\epsilon_n$ is too large, which here means that $(\epsilon_n - \epsilon^*) \gg 1/v_n$, then the asymptotic distribution of the ABC posterior is uniform with a radius that is of the order $\epsilon_n - \epsilon^*$. In contrast to the case of correct model specification, if $\epsilon^*>0$ and if $v_n\{\epsilon_n - \epsilon^* \}\rightarrow 2c \in \mathbb R$, then the limiting distribution is no longer Gaussian. Moreover, this result maintains even if $c= 0$.}
\end{remark} 

\begin{remark}\normalfont{
In likelihood-based Bayesian inference, credible sets are not generally valid confidence sets if the model is misspecified, however, the resulting posterior is still asymptotically normal (see, e.g., \citealp{kleijn2012} and \citealp{Muller2013}). In the case of ABC, not only will credible sets not be valid confidence sets, but the asymptotic shape of the ABC posterior is not even Gaussian. }
\end{remark}

\begin{remark}\normalfont{
In practice $\epsilon^* $ is unknown, and it is therefore not possible to choose $\epsilon_n$ directly. However, we note that the application of ABC is most often implemented by accepting draws of $\theta$ within some pre-specified (and asymptotically shrinking) quantile threshold; i.e., one accepts a simulated draw $\theta^i$ if $d( \eta (\mathbf z^i ) , \eta(\mathbf y) )$ is smaller than the $\alpha$-th empirical quantile of the simulated values $d( \eta (\mathbf z^j ),  \eta(\mathbf y) )$, $j\leq N$. However, as discussed in \citealp{ FMRR2016}, the two representations of the ABC approach are dual in the sense that  choosing a value of $\alpha $ on the order of $\delta v_n^{-k_\theta}$, with $\delta$ small, corresponds to choosing $|\epsilon_n-\epsilon^*| \lesssim \delta^{1/k_\eta} v_n$ and choosing $\alpha_n \gtrsim M  v_n^{-k_\theta}$ corresponds to choosing $\epsilon_n- \epsilon^* \gtrsim M v_n$. We further elaborate on the equivalence between both approaches in Section \ref{sec:detection}. }
\end{remark}

Interestingly, the proof of Theorem \ref{normal_thm} (provided in the appendix) demonstrates that if $v_n(\epsilon_n - \epsilon^*) \rightarrow - \infty$, in particular when $\epsilon_n= o(1)$ and $\epsilon^*>0$, posterior concentration of $\Pi_{\epsilon}[\cdot|\eta(\y)]$ need not occur. We present an illustration of this phenomena in the following simple example. 
\medskip 

\noindent\textbf{Example 2:}
Consider the case where $k_\theta = 1$ and $k_\eta = 2$. Let $\tilde Z_y = \sqrt{n}( \eta(\y) - b_0)$ and $\tilde Z_n = \sqrt{n} (\eta (\mathbf z) - b(\theta)) $, where $\tilde Z_n \sim \mathcal N( 0, v_\theta^2 I_2)$, for $v_{\theta}$ some known function of $\theta$,   and $b(\theta) = ( \theta, \theta)^{\intercal}$. In addition, assume that $b_0 = (\bar{b}_0, -\bar{b}_0)$, with $\bar{b}_0\neq0$. Under this setting, and when $\|\cdot\|$ is the Euclidean norm, it follows that the unique pseudo-true value is $\theta^* = 0.$ However, depending on $v_\theta$, the approximate posterior need not concentrate on $\theta^* = 0$. This is summarized in the following Proposition.
\begin{proposition}\label{prop:counter}
In the setup described above, if $v_\theta/v_{\theta^*} = \sigma(\theta ) $, for $v_{\theta^*}$ some known function, such that $\sigma$ is continuous and  $\sigma(\bar b_0/2)^2\geq 3$, and if the prior has positive and continuous  density on $[-\bar b_0, \bar b_0]$, then 
\begin{equation*}
\Pi_\epsilon\left\{ | \theta - \theta^* | \leq  \delta | \eta(\mathbf y) \right\}   = o\left( \Pi_\epsilon\left\{ | \theta - \bar b_0/2| \leq  \delta | \eta(\mathbf y) \right\}  \right) = o(1) .
\end{equation*}
\end{proposition}

\section{Local Regression Adjustment under Misspecification}
\subsection{Posterior Concentration}
{Local regression adjustments to ABC have found broad applicability with practitioners. However,  we caution against the blind application of local regression adjustment when one is willing to entertain the idea of model misspecification. As demonstrated by the introductory example, the use of this particular adjustment can lead to point estimators that behave very differently to those obtained from Algorithm \ref{ABC}, even in small samples. }

In this section, we first rigorously characterize posterior concentration of local linear regression adjustment ABC (ABC-Reg) under model misspecification. Using this initial result, we then extend the conclusion to local nonlinear regression adjustment approaches. For simplicity, we only consider the case of scalar $\theta$, however, we allow $\eta(\mathbf{y})$ to be multi-dimensional.\footnote{This result can be extended at the cost of more complicated arguments but we refrain from this setting to simplify the interpretation of our results.} 

ABC-Reg first runs Algorithm \ref{ABC}, with tolerance $\epsilon_{n}$, to obtain a set of accepted draws and summaries $\{\theta^i,\eta(\mathbf{z}^i)\}$, and then uses a linear regression model to adjust the accepted values of $\theta$. {In this way, the original accepted value $\theta^i$ is artificially related to $\eta(\mathbf{y})$ and $\eta(\mathbf{z})$ through the} linear regression model $$\theta^i=\mu+\beta^\intercal \{\eta(\mathbf{y})-\eta(\mathbf{z}^i)\}+\nu_i,$$ where $\nu_i$ denotes the model residual. Define $\bar{\theta}=\sum_{i=1}^{N}\theta^i/N$ and $\bar{\eta}=\sum_{i=1}^{N}\eta(\mathbf{z}^i)/N$. Given $\theta^i$, ABC-Reg then produces an adjusted parameter draw according to 
\begin{flalign*}
\tilde{\theta}^i&=\theta^i-\hat{\beta}^{\intercal}\{\eta(\mathbf{z}^i)-\eta(\mathbf{y}^{})\},\\
\hat{\beta}&= \left[\frac{1}{N_{}}\sum_{i=1}^{N}\left(\eta(\mathbf{z}^{i})-\bar{\eta}\right)\left(\eta(\mathbf{z}^{i})-\bar{\eta}\right)^{\intercal}\right]^{-1}\left[\frac{1}{N_{}}\sum_{i=1}^{N}\left(\eta(\mathbf{z}^{i})-\bar{\eta}\right)\left(\theta^i-\bar{\theta}\right)\right]=\widehat{\text{Var}}^{-1}(\eta(\mathbf{z}^{i}))\widehat{\text{Cov}}(\eta(\mathbf{z}^{i}),\theta^i)
\end{flalign*}
Therefore, for ${\theta}^i\sim \Pi_{\epsilon}[\cdot|\eta(\mathbf{y})]$, the posterior measure for $\tilde{\theta}^i$ is nothing but a scaled and shifted version of $\Pi_{\epsilon}[\cdot|\eta(\mathbf{y})]$. Consequently, the asymptotic behavior of the ABC-Reg posterior, denoted by $\widetilde{\Pi}_{\epsilon}[\cdot|\eta(\mathbf{y})]$, is determined by the behavior of ${\Pi}_{\epsilon}[\cdot|\eta(\mathbf{y})]$, $\hat{\beta}$, and $\{\eta(\mathbf{y})-\eta(\mathbf{z}^i)\}$. 

The following result describes the asymptotic behavior of the ABC-Reg posterior $\widetilde{\Pi}_{\epsilon}[\cdot|\eta(\mathbf{y})]$.

\begin{corollary}\label{th:asympAdj}
Assume that \textbf{[A0\,{\em{$^{\prime}$}}]}, \textbf{[A1]} and \textbf{[A2]} are satisfied and $\epsilon_n \downarrow \epsilon^* $ with $$\epsilon_n \geq \epsilon^* + M v_n^{-1} + v_{0,n}^{-1},$$ for $M$ large enough. Furthermore, assume that for some $\beta_0$ with $\|\beta_0\|>0$, $\|\hat{\beta}-\beta_0\|=o_{P_{\theta}}(1)$. Define $\tilde{\theta}^*=\theta^*-\beta_0^{\intercal}(b(\theta^*)-b_0)$.  Let $M_n$ be any positive sequence going to infinity and 
$\delta_n \geq  M_n(\epsilon_n - \epsilon^*)$, then
\begin{equation*}
\widetilde{\Pi}_{\epsilon}[ |{\theta}-\tilde{\theta}^*|> \delta |\eta( \mathbf{y} ) ] = o_{P_0}(1),
\end{equation*}as soon as
 \begin{equation*}
\begin{split}
\delta_n &\geq M v_n^{-1}u_n^{-D/\kappa}  = o(1) \quad \mbox{in case (i) of assumption [A1]}\\
\delta_n&\geq  M v_n^{-1}|\log(u_n)|^{1/\tau} = o(1)  \quad \mbox{in case (ii) of assumption [A1] }.
\end{split}
\end{equation*}
with $u_n =\epsilon_n - (\epsilon^*  +  M v_n^{-1}  + v_{0,n}^{-1}  )\geq 0$.
\end{corollary}

\begin{remark}\label{remark:cons}\normalfont{
An immediate consequence of Corollary \ref{cor1} and  \ref{th:asympAdj} is that the ABC posterior $\Pi_{\epsilon}[\cdot|\eta(\mathbf{y})]$ concentrates mass on $$\theta^*=\arg\inf_{\theta\in\Theta}d(b(\theta),b_{0}),$$while the ABC-Reg posterior $\widetilde{\Pi}_{\epsilon}[\cdot|\eta(\mathbf{y})]$ concentrates mass on 
$$\tilde{\theta}^*=\theta^*-\beta^{\intercal}_{0}(b(\theta^*)-b_{0}).$$
{As a consequence, ABC-Reg takes draws of $\theta$ that are (asymptotically) optimal in terms of minimizing the chosen distance, $d(\cdot,\cdot)$, between observed and simulated summaries, and perturbs them in a (linear) manner that need not preserve the optimality of the original draws. Furthermore, for $\|\beta_0\|$ large, the pseudo-true value $\tilde{\theta}^*$, onto which ABC-Reg concentrates, can easily lie outside $\Theta$. Therefore, if the model is misspecified, there is no guarantee that ABC-Reg returns draws that are optimal in terms of minimizing $d(\cdot,\cdot)$ and there is even no guarantee that ABC-Reg returns values in $\Theta$.\footnote{In the appendix, we give a concrete example of this later behavior in the confines of the $g$-and-$k$ distribution.}}}
\end{remark}

\begin{remark}\normalfont{
	{Crucially, the result of Corollary \ref{th:asympAdj} and the phenomena discussed in Remark \ref{remark:cons}, are not restricted to local linear regression adjustment, and extend to the nonlinear variety of regression adjustment. For brevity, we only sketch the general idea here, and note that a rigorous proof follows along the same lines as Corollary \ref{th:asympAdj}, and is therefore omitted.} Consider the nonlinear regression model $$\theta = m(\eta(\mathbf z) ) + v, $$ for some unknown function $m(\cdot)$. Denote by $\hat m(\cdot )  $ a nonparametric estimator of the unknown regression function that is constructed using the accepted draws $\{\theta^\ell, \eta(\mathbf z^\ell)\}_{\ell \leq L} $ obtained from the ABC posterior $\Pi_{\epsilon}[\cdot|\eta(\mathbf{y})]$. A nonlinear regression post-processing approach transforms the accepted $\theta^\ell $ into 
	\begin{equation}\label{nonlinear}
	\tilde{\theta}^\ell = \theta^\ell+\left\{\hat m(\eta(\mathbf y) ) - \hat m(\eta(\mathbf z^\ell) ) \right\}, \text{ for }\ell=1,\dots, L. 
	\end{equation}
	Under the regularity conditions \textbf{[A0{$^{\prime}$}]}, \textbf{[A1]} and \textbf{[A2]}, if the (nonparametric) estimator $\hat m(\cdot)$ converges to a function $m^*(\cdot)$, $\hat m ( \eta ) = m^*(\eta) + o_p(1)$, uniformly in a neighbourhood of $b_0$ and $b(\theta^*)$, then by the concentration of $\eta(\mathbf y)$ and $\eta(\mathbf z^\ell)$, equation \eqref{nonlinear} becomes
	\begin{flalign*}
	\tilde{\theta}^\ell &= \theta^\ell+\left\{m^*(b_0)  - m^*(b(\theta^*) )\right\} + o_p(1) = \theta^* + \left\{m^*(b_0) - m^*(b(\theta^*) )\right\} + o_p(1) .
	\end{flalign*}
	The last equality follows from the posterior concentration of $\Pi_{\epsilon}[\cdot|\eta(\mathbf{y})]$ toward $\theta^* $, i.e., $\theta^\ell = \theta^*+o_p(1)$. Hence, as soon as $m^*(b_0) - m^*(b(\theta^*) )\neq 0$ the nonlinear regression post-processed ABC posterior concentrates onto a different value than $\Pi_{\epsilon}[\cdot|\eta(\mathbf{y})]$, which offers no particular justification. Moreover, given that $\hat m(\cdot) $ is constructed using the simulated data, there is no reason to suspect that $m^*(b_0) = m^*(b(\theta^*) )$. Lastly, we note that the above computations are not significantly altered if instead we had considered a nonlinear conditional heteroskedastic regression model, as in \cite{blumF2010non}, from the outset. As a consequence, the nonlinear conditional heteroskedastic regression adjustment will have similar asymptotic behavior to the nonlinear regression adjustment.}	
\end{remark}

\begin{remark}\normalfont{
	An additional consequence of Corollary \ref{cor1} and \ref{th:asympAdj} is that the ABC posterior $\Pi_{\epsilon}[\cdot|\eta(\mathbf{y})]$ and the ABC-Reg posterior  $\widetilde{\Pi}_{\epsilon}[\cdot|\eta(\mathbf{y})]$ will yield different posterior expectations. This difference between expectations calculated under $\Pi_{\epsilon}[\cdot|\eta(\mathbf{y})]$  and $\widetilde{\Pi}_{\epsilon}[\cdot|\eta(\mathbf{y})]$ explains the divergence between the ABC-AR and ABC-Reg posterior means observed in Figure \ref{fig1}. In the following section, we use this behavior to derive a procedure for detecting model misspecification. }
\end{remark}

\subsection{Adjusting Local Regression Adjustment}
The difference between accept/reject ABC and ABC-Reg under model misspecification is related to the regression adjustments re-centering of the accepted draws $\theta^i$ by $\hat{\beta}^{\intercal}\{\eta(\y)-\eta(\z)\}$.  Whilst useful under correct model specification, when the model is misspecified the adjustment can force $\theta^i$ away from $\theta^*$ and towards $\tilde{\theta}^*$, which need not lie in $\Theta$ or be optimal in terms of minimizing $d(\cdot,\cdot)$. 

The cause of this behavior is the inability of $\eta(\z)$ to replicate the asymptotic behavior of $\eta(\y)$, which in the terminology of \cite{Marinetal2014} means that the model is incompatible with the observed summaries. This incompatibility of the summary statistics ensures that the influence of the centering term $\hat{\beta}^{\intercal}\{\eta(\y)-\eta(\z)\}$ can easily dominate that of the accepted draws $\theta^i$, with the introductory example being just one example of this behavior. 

In an attempt to maintain the broad applicability of local linear regression adjustment in ABC, and still ensure it gives sensible results under model misspecification, we propose a useful modification of the regression adjustment approach. To motivate this modification recall that, under correct model specification and regularity conditions, at first-order the linear regression adjustment approach ensures (see Theorem 4 in \citealp{FMRR2016}):
\begin{flalign}
\tilde{\theta}^{i}&=\theta^i+\hat{\beta}^{\intercal}\{\eta(\y)-\eta(\z)\}\nonumber\\&=\theta^i+\hat{\beta}^{\intercal}\{b_0-b(\theta^i)\}+O_{p}(1/v_{n})\nonumber \\&=\theta^i-\left[\nabla_{\theta}b(\theta^*)^{\intercal}V_0^{-1}\nabla_{\theta}b(\theta^*)\right]^{-1}\nabla_{\theta}b({\theta}^*)^{\intercal}V_{0}^{-1}\nabla_{\theta}b(\bar{\theta})(\theta^i-\theta^*)+O_{p}(1/v_{n})\label{wow},
\end{flalign} where $b_0=b(\theta^*)$ by correct model specification, $\bar{\theta}$ is an intermediate value satisfying $|\bar{\theta}-\theta^*|\leq|\theta^i-\theta^*|$, $V_{0}=\lim_n\text{Var} [\sqrt{n}\{\eta(\y)-b_0\}]$, {and the third line follows from a mean-value expansion and the definition of the local linear regression adjustment}. Therefore, it follows from \eqref{wow} that, even if $k_{\eta}>k_{\theta}$, the dimension of $\eta(\y)$ will not impact the asymptotic variance of the ABC-Reg posterior mean. This result, at least in part, helps explain (from a technical standpoint) the popularity of the ABC-Reg approach as a dimension reduction method. 

However, under model misspecification,  $b_0\neq b(\theta)$ for any $\theta\in\Theta$, and hence there does not exist an intermediate value $\bar{\theta}$ such that $$b_0-b(\theta^i)\neq \nabla_{\theta}b(\bar{\theta})(\theta^*-\theta^i).$$As a consequence, equation \eqref{wow} cannot be valid (in general) if the model is misspecified. 

The behavior of ABC-Reg under correct and incorrect model specifications suggests that the methods poor behavior  under the latter can be mitigated by replacing $\eta(\y)$ with an alternative term. To this end, define $\hat{\theta}=\int \theta d\Pi_\epsilon[\theta|\eta(\y)]$ to be the posterior mean of accept/reject ABC. Let $\hat{\z}^{m}$, $m=1,...,M$, be a pseudo-data set of length $n$ simulated under the assumed DGP and at the value $\hat{\theta}$, and define$$\hat{\eta}=\sum_{m=1}^{M}\eta(\hat{\z}^m)/M.$$ Using $\hat{\eta}$, we can then implement the modified local linear regression adjustment
$$\breve{\theta}^{i}=\theta^{i}+\hat{\beta}^{\intercal}\{\hat{\eta}-\eta(\z^{i})\}.$$

The key to this modified approach is that under correct specification $\hat{\eta}$ behaves like $\eta(\y)$, while under incorrect specification $\hat{\eta}$ behaves like $\eta(\z)$.\footnote{{While the choice of $M$ will not matter asymptotically, we argue that $M$ should be chosen so that the variability of $\hat{\eta}$ is small relative to $\eta(\mathbf{z})$.}} A direct consequence of this construction is that this approach avoids the incompatibility issue that arises under model misspecification. {In addition, since this new regression adjustment approach uses a centering sequence calculated from the accept/reject ABC posterior mean, $\hat{\theta}$, the asymptotic behavior of this new approach is similar to accept/reject ABC.}

\medskip 

\noindent \textbf{Example 1 (Continued):} 
 Recall that the assumed DGP is $z_1,\dots,z_n$ iid as $\mathcal{N}(\theta,1)$ but the actual DGP is $y_1,\dots,y_n$ iid as $\mathcal{N}(\theta,{\sigma}^2)$. ABC is conducted using the following summary statistics:
\begin{itemize}
\item the sample mean $\eta_{1}(\mathbf{y})=\frac{1}{n}\sum_{i=1}^{n}{y}_{i}$,
\item the sample variance  $\eta_{2}(\mathbf{y})=\frac{1}{n-1}\sum_{i=1}^{n}({y}_{i}-\eta_{1}(\mathbf{y}))^{2}$.
\end{itemize}
We consider three different DGPs corresponding to $\sigma^{2}\in\{1,2,3\}$. For each of these cases we generate 1000 artificial samples for $\y$ of length $n=100$ and apply four different ABC approaches: the accept/reject ABC approach (ABC-AR), the local linear regression adjustment of \cite{Beaumont2025} (ABC-Reg), our new local linear regression adjustment (ABC-RegN), and the nonlinear regression adjustment of \cite{blumF2010non} that is fit using neural nets (ABC-NN). Each procedure relies on $N=25, 000$ pseudo-datasets generated according to ${z}^{j}_{i}\sim \mathcal{N}(\theta^{j},1)$, where the prior is again given by $\theta\sim \mathcal{N}(0,25)$. {For each of the procedures we set the tolerance $\epsilon$ to be the 1\% quantile of the simulated distances $\|\eta(\y)-\eta(\z^i)\|$.}

\begin{figure}[h]
	\resizebox{\textwidth}{!}{\input{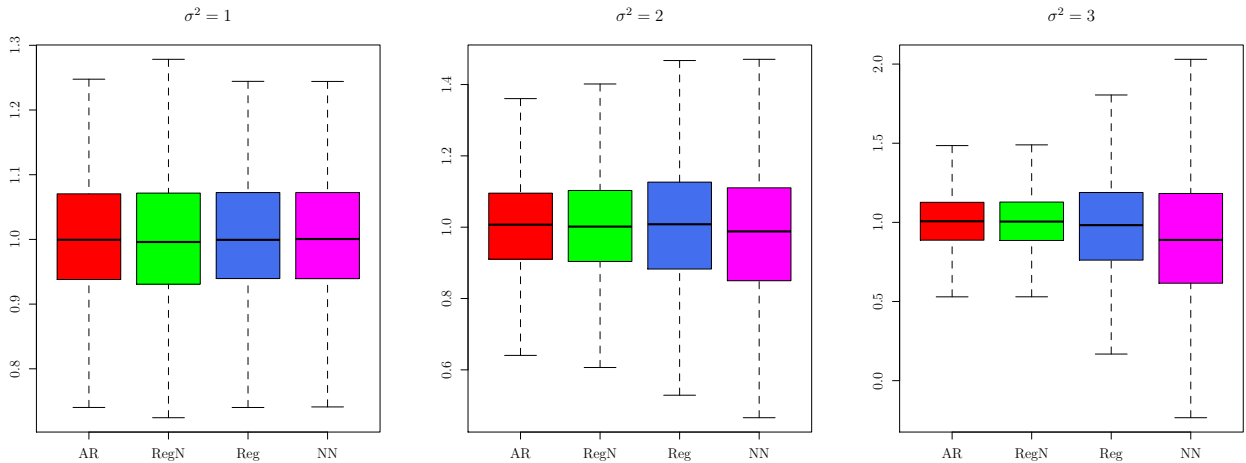}  }  
\caption{{Posterior mean comparison of ABC-AR (AR), standard local linear regression adjustment (Reg), the new regression adjustment approach (RegN), and the local nonlinear regression adjustment (NN) across $\sigma^2\in\{1,2,3\}$. Recall that $\sigma^2=1$ corresponds to correct model specification. }} 
\label{AbcNew}
\end{figure}

Figure \ref{AbcNew} plots the posterior mean of each approach across the Monte Carlo replications and
across all designs. The results demonstrate that the new regression adjustment maintains
stable performance across both correct and incorrect model specification, while the point estimators obtained from ABC-Reg and ABC-NN are more varied. More specifically, for $\sigma^2\in\{1,2\}$ we see that all regression adjustments tend to give similar results. However, for $\sigma^2=3$, it is clear that the traditional linear and nonlinear adjustment approaches produce point estimators with more variability, across the repeated samples, than the other ABC approaches. This additional variability is a direct consequence of the fact that ABC-Reg and ABC-NN enforce a regression relationship between the accepted draws $\theta^i$ and $\eta(\y)-\eta(\mathbf{z})$ when one does not necessarily exist. As we have seen in Corollary \ref{th:asympAdj}, enforcing this additional (wrong) model to produce new values of $\theta^i$, which need not respect the actual relationship between $\theta^i$ and $\eta(\y)-\eta(\mathbf{z})$, will (randomly) shift the adjusted draws away from their initial center of posterior mass, $\theta^*$, and can yield point estimators with more variability in a repeated sampling context.

{In addition to the results for the posterior means, in Table \ref{tab:newTable} we record for each method the posterior standard deviation, the length of the corresponding 95\% credible set and the Monte Carlo coverage across the different designs. The values given in Table \ref{tab:newTable} represent the average values of these quantities across the Monte Carlo replications.}

The results demonstrate that all of the local regression adjustments, ABC-Reg, ABC-RegN, and ABC-NN, have much smaller posterior variability and much shorter credible sets than ABC-AR (on average). As a consequence, when the model is misspecified, this behavior gives researchers a false sense of precision, and leads to poor coverage rates (for the pseudo-true value) across all the local regression adjustment procedures (linear and nonlinear). Therefore, even though our new regression adjustment procedure gives stable performance under correct and incorrect model specification, it still suffers from the coverage issues alluded to in the remarks given after Theorem \ref{normal_thm}. 

Furthermore, additional numerical experiments conducted in the supplemental appendix demonstrate that the use of the so-called heteroskedasticity correction for the local regression adjustment, as in \cite{blumF2010non}, does not significantly alter these results. More specifically, the resulting heteroskedasticity corrected local regression adjustment approaches yield results that are very similar to those obtained without the heteroskedasticity correction. We refer the interested reader to the appendix for these details.

From this simple example we can conclude that under model misspecification local regression adjustment approaches, both linear and nonlinear, with or without a heteroskedasticity correction, can lead to significant overconfidence in the resulting point estimates obtained from these methods, and can result in poor coverage rates for the pseudo-true value.

\vspace{.5cm}

\begin{table}[h!]
\caption{{Monte Carlo coverage (Cov), credible set length (Len), and posterior standard deviation (Std) for the simple normal example under various levels of model misspecification. Cov is the percentage of times that the 95\% credible set contained $\theta=1$. Len is the average length of the credible set across the Monte Carlo trials. Std is the average posterior standard deviation across the Monte Carlo trials.}}
\centering%
	\begin{tabular}{lrrrrrrr}
	&       & \multicolumn{1}{c}{ABC-AR} &       &       &       & \multicolumn{1}{c}{ABC-RegN} &  \\
	& \multicolumn{1}{c}{Cov} & \multicolumn{1}{c}{Len} & \multicolumn{1}{c}{Std} &       & \multicolumn{1}{c}{Cov} & \multicolumn{1}{c}{Len} & \multicolumn{1}{c}{Std} \\
	$\sigma^2=1$    & 0.9820 &\multicolumn{1}{c} {0.4666} & 0.1221 &       & 0.9380 & \multicolumn{1}{c}{0.3851} & 0.1001 \\
	$\sigma^2=2$    & 0.9610 & \multicolumn{1}{c}{0.6147} & 0.1576 &       & 0.8020 & \multicolumn{1}{c}{0.3837} & 0.0998 \\
	$\sigma^2=3$   & 0.9130 & \multicolumn{1}{c}{0.6164} & 0.1581 &       & 0.7070 & \multicolumn{1}{c}{0.3839} & 0.0997 \\
	&       &       &       &       &       &       &  \\
	&       & \multicolumn{1}{c}{ABC-Reg} &       &       &       & \multicolumn{1}{c}{ABC-NN} &  \\
	& \multicolumn{1}{c}{Cov} & \multicolumn{1}{c}{Len} & \multicolumn{1}{c}{Std} &       & \multicolumn{1}{c}{Cov} & \multicolumn{1}{c}{Len} & \multicolumn{1}{c}{Std} \\
	$\sigma^2=1$    & 0.9410 & \multicolumn{1}{c}{0.3820} & 0.0997 &       & 0.9500 & \multicolumn{1}{c}{0.3853} & 0.1006 \\
	$\sigma^2=2$   & 0.7170 & \multicolumn{1}{c}{0.3826} & 0.0998 &       & 0.7290 & \multicolumn{1}{c}{0.4440} & 0.1228 \\
	$\sigma^2=3$    & 0.4600 & \multicolumn{1}{c}{0.3821} & 0.0997 &       & 0.4190 & \multicolumn{1}{c}{0.5043} & 0.1490 \\
\end{tabular}%
	\label{tab:newTable}%
\end{table}\hfill$\square$

\section{Detecting Misspecification}In this section we propose two methods to detect model misspecification in ABC. The first approach is based on the behavior of the acceptance probability under correct and incorrect model specification. The second approach is based on comparing posterior expectations calculated under $\Pi_{\epsilon}[\cdot|\eta(\mathbf{y})]$ (obtained from Algorithm \ref{ABC}) and $\widetilde{\Pi}_{\epsilon}[\cdot|\eta(\mathbf{y})]$ (obtained using local linear regression adjustment). 
\subsection{A Simple Graphical Approach to Detecting Misspecification }\label{sec:detection}
From the results of \citet{FMRR2016}, under regularity and correct model specification, the acceptance probability ${\alpha}_{n}=\text{Pr}\left[d(\eta(\mathbf{y}),\eta(\mathbf{z}))\leq \epsilon_{n}\right]$ satisfies, for $n$ large and $\epsilon_{n}\gg v_{n}^{-1},$ $$\alpha_{n}=\text{Pr}\left[d(\eta(\mathbf{y}),\eta(\mathbf{z}))\leq \epsilon_{n}\right]\asymp \epsilon_{n}^{k_{\theta}}.$$ In this way, as $\epsilon_{n}\rightarrow0$ the acceptance probability $\alpha_{n}\rightarrow0$ in a manner that is approximately linear in $\epsilon_{n}^{k_{\theta}}$.  

However, this relationship between $\alpha_{n}$ and $\epsilon_{n}$ does not extend to the case where the model is misspecified. In particular, if $\epsilon^{*}>0$, once $\epsilon_{n}<\epsilon^*$ the acceptance probability $\alpha_{n}$ will be small or zero, even for a large number of simulations $N$.

The behavior of $\alpha_{n}$ under correct and incorrect model specifications means that one can potentially diagnose misspecification by comparing the behavior of $\alpha_{n}$ over a decreasing sequence of tolerance values. In particular, if we take a decreasing sequence of equally spaced tolerances $\epsilon_{1,n}< \epsilon_{2,n}<\cdots<\epsilon_{J,n}$ we can construct and plot the resulting sequence $\{\alpha_{j,n}\}_{j}$ to determine if $\{\alpha_{j,n}\}_{j}$ decays in an (approximately) linear fashion as $\epsilon_{j,n}$ decreases. 

While $\alpha_{n}$ is infeasible to obtain in practice, the same procedure can be applied with $\alpha_{n}$ replaced by the estimator $\hat{\alpha}_{n}=\sum_{i=1}^{N}\1_{[d(\eta(\mathbf{y}),\eta(\mathbf{z}))\leq\epsilon_{n}]}/N$. In this way, such a graphical check can easily be performed using the ABC reference table. The only difference is that, instead of considering a single tolerance $\epsilon$, one would consider a sequence of tolerances $\{\epsilon_{j,n}\}_{j}$ and record, for each $j$, $$\hat{\alpha}_{j,n}=\sum_{i=1}^{N}\1_{[d(\eta(\mathbf{y}),\eta(\mathbf{z}))\leq\epsilon_{j,n}]}/N.$$ Once $\hat{\alpha}_{j,n}$ has been obtained, it can be plotted against ${\epsilon}_{j,n}^{k_\theta}$ (in some fashion) and the relationship can be analyzed to determine if deviations from linearity are in evidence. 

To understand exactly how such a procedure can be implemented, we return to the simple normal example. 

\medskip 

\noindent \textbf{Example 1 (Continued):} The assumed DGP is $z_1,\dots,z_n$ iid as $\mathcal{N}(\theta,1)$ but the actual DGP is $y_1,\dots,y_n$ iid as $\mathcal{N}(\theta,{\sigma}^2)$. 
We again consider ABC using the following summary statistics:
\begin{itemize}
\item the sample mean $\eta_{1}(\mathbf{y})=\frac{1}{n}\sum_{i=1}^{n}{y}_{i}$,
\item the sample variance  $\eta_{2}(\mathbf{y})=\frac{1}{n-1}\sum_{i=1}^{n}({y}_{i}-\eta_{1}(\mathbf{y}))^{2}$.
\end{itemize}
Taking $\sigma^{2}\in\{1,1+1/9,...,1+8/9\}$, we generate observed samples of size $n=100$ according to ${y}_i\sim\mathcal{N}(1,\sigma^{2})$, iid, where, for each of the nine different simulated datasets, we keep the random numbers fixed and only change $\sigma^{2}$. We consider $N=25, 000$ simulated datasets generated iid according to ${z}^{j}_{i}\sim \mathcal{N}(\theta^{j},1)$, with $\theta^{j}\sim \mathcal{N}(0,25)$, and for $d(\cdot,\cdot)$ we take the Euclidean norm. { For the sequence of $\epsilon_{j,n}$ values, we consider $J=100$ evenly spaced increments, with $\epsilon_{J,n}$ chosen to correspond to the $10\%$ quantile of the simulated distances, and where $\epsilon_{1,n}$ is taken to be the $0.1\%$ quantile of the simulated distances.}
	
\begin{figure}[h]
\centering
\resizebox{\textwidth}{!}{\input{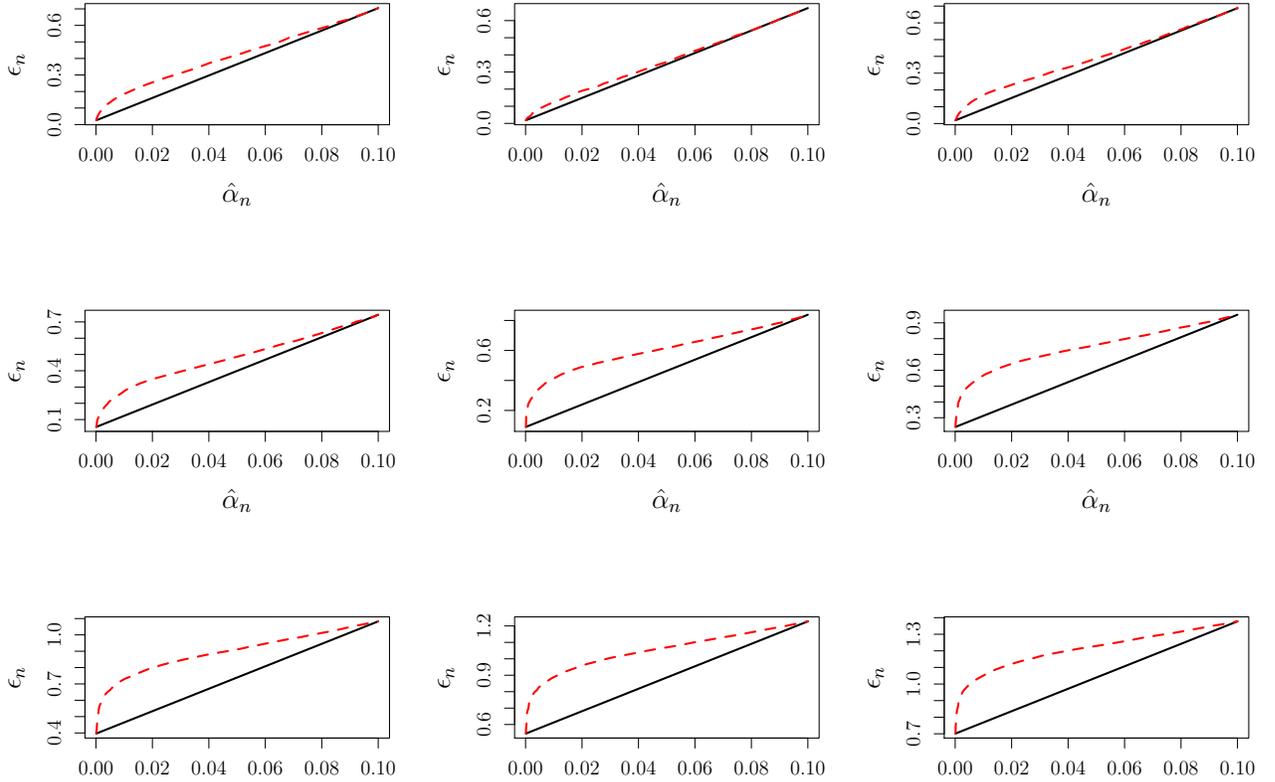} } 
\caption{Graphical comparison of empirical acceptance probabilities $\hat{\alpha}_{j,n}$ (dashed line) and theoretical acceptance probabilities (solid line) for a decreasing sequence of tolerance values $\epsilon_{j,n}$. } 
\label{fig_diag}
\end{figure}

{In Figure \ref{fig_diag}, we plot the results over the nine different levels of misspecification. Each figure contains two different curves: the dashed curve represents the observed relationship between $\hat{\alpha}_{j,n}$ and $\epsilon_{j,n}$, while the solid curve plots a linear relationship between  $\hat{\alpha}_{j,n}$ and $\epsilon_{j,n}$ that can be used to help visually diagnose departures from linearity.\footnote{Across the different values of $\sigma^2$, the solid line is constructed using the endpoint pairs $(\hat{\alpha}_{J,n},\epsilon_{J,n})=(0.10,\epsilon_{J,n})$ and $(\hat{\alpha}_{J,n},\epsilon_{J,n})=(0.001,\epsilon_{1,n})$.} }

We recall that, in this example correct specification warrants a linear relationship between $\hat\alpha_{n}$ and $\epsilon_{n}$, since we are only conducting inference on a single parameter. More generally, under correct model specification, we would expect a linear relationship between $\hat\alpha_{n}$ and $\epsilon_{n}^{k_{\theta}}$.

{ Analyzing Figure \ref{fig_diag}, we see that the relationship is fairly linear for $\sigma^2\leq1+1/3$. However, for $\sigma^2>1+1/3$ the relationship between the acceptance probabilities $\hat{\alpha}_{j,n}$ and $\epsilon_{j,n}$ exhibits significant nonlinear behavior. Therefore, in this example the diagnostic would suggest that the model is misspecified once $\sigma^2>1+1/3$, which is evidenced by the nonlinear relationship between $\hat{\alpha}_{n,j}$ and $\epsilon_{j,n}$.}

Clearly, obtaining broad conclusions about model misspecification from this graphical approach depends on many features of the underlying model, the dimension of $\theta$, and the exact nature of misspecification. However, it is always possible to benchmark the results for the observed data against those obtained under the ABC reference table. That is, using the ABC reference table, we can easily apply the above diagnostic to one, or even several, ``observed'' series from the reference table. Then, if the relationship between $(\hat{\alpha}_{j,n},\epsilon_{j,n})$ in the observed data deviates from that observed in the reference table, this is evidence that the model is misspecified.

\hfill$\square$

\subsection{Detecting Model Misspecification Using Regression Adjustment}
Corollary \ref{cor1} and \ref{th:asympAdj} demonstrate that accept/reject ABC (ABC-AR) and local linear regression adjustment ABC ( ABC-Reg), place posterior mass in different regions of the parameter space. Therefore, for $\theta\mapsto h(\theta)$ a known, smooth function, under model misspecification the posterior expectation of $h(\cdot)$ under ABC-AR and ABC-Reg,
\begin{flalign*}
\hat{h}&=\int h(\theta)d\Pi_{\epsilon}[\theta|\eta(\mathbf{y})], \quad 
\tilde{h}=\int  h(\theta)d\widetilde{\Pi}_{\epsilon_{}}[\theta|\eta(\mathbf{y})],
\end{flalign*}will converge in probability, as $n\rightarrow+\infty$ and $\epsilon_{n}\downarrow\epsilon^*$, to distinct values. However, if the model is correctly specified, it can be shown that $\hat{h}_{}$ and $\tilde{h}_{}$ will not differ, up to first-order, so long as $\epsilon_{n}=o(1/\sqrt{n})$ {(this result follows from Theorem 4 in \citealp{FMRR2016}, or Theorem 3.1 in \citealp{LF2016a}).} Therefore, a useful approach for detecting model misspecification is to compare various posterior expectations, such as moments or quantiles, calculated under the two posteriors. 

More specifically, if the model is correctly specified and if we use Algorithm 1 based on quantile thresholding with $\alpha_{n}= \delta n^{-k_\theta/2}$, for $\delta >0$, then $$\sqrt{n}\|\hat{h}-\tilde{h}\|=o_{P_{0}}(1).$$ However, if $\epsilon^*=\inf_{\theta\in\Theta}d(b_{0},b(\theta))>0$, under regularity conditions, we can deduce that $$\|\hat{h}-\tilde{h}\|=O_{P_{0}}(1).$$Therefore, if $\sqrt{n}\|\hat{h}-\tilde{h}\|$ is large, this is meaningful evidence that the model may be misspecified.

Detecting misspecification by analyzing the magnitude of $\sqrt{n}\|\hat{h}-\tilde{h}\|$ requires the specification of a cutoff value, denoted by $t_n$, such that if $\sqrt{n}\|\hat{h}-\tilde{h}\|$ is larger than $t_n$, we conclude that the model is likely misspecified. While there are several ways to choose the cutoff value $t_n$, we propose a simulation-based approach that uses the ABC reference table. Namely, we use the fact that for any fixed value of $\theta$, which for instance could be drawn from the prior or from the ABC-AR posterior, we can always simulate ``observed data'' from the assumed model  to determine what the magnitude of $\sqrt{n}\|\hat{h}-\tilde{h}\|$ should be under correct specification. More specifically, for a given value of $\theta$, we can always simulate an ``observed data series'', and then run both ABC-AR and ABC-Reg on this simulated observed data to calculate $\sqrt{n}\|\hat{h}-\tilde{h}\|$ under correct model specification. 

To operationalize this approach, we generate $b=1,\dots,B$ such simulated ``observed'' datasets, all at the same value of $\theta$, and for each of these datasets we calculate $\sqrt{n}\|\hat{h}_{b}-\tilde{h}_{b}\|$. The cutoff value $t_n$ can then be defined as an empirical quantile of the simulated distances $\{\sqrt{n}\|\hat{h}_{b}-\tilde{h}_{b}\|\}_{b=1}^{B}$, which have all been calculated under correct specification. For instance, $t_n$ could be defined as the 95\% quantile of $\{\sqrt{n}\|\hat{h}_{b}-\tilde{h}_{b}\|\}_{b=1}^{B}$. The value of $t_n$ obtained from this procedure can then be compared with the corresponding value of $\sqrt{n}\|\hat{h}-\tilde{h}\|$ obtained from the actual observed data $\y$, and if $\sqrt{n}\|\hat{h}-\tilde{h}\|>t_n$ we conclude that the model is likely misspecified.

We demonstrate this approach to diagnosing model misspecification in our simple running example.

\medskip 
\noindent \textbf{Example 1 (Continued):} The assumed DGP is $z_1,\dots,z_n$ iid as $\mathcal{N}(\theta,1)$ but the actual DGP is $y_1,\dots,y_n$ iid as $\mathcal{N}(\theta,{\sigma}^2)$. We again consider the following summary statistics:
\begin{itemize}
\item the sample mean $\eta_{1}(\mathbf{y})=\frac{1}{n}\sum_{i=1}^{n}{y}_{i}$,
\item the sample variance $\eta_{2}(\mathbf{y})=\frac{1}{n-1}\sum_{i=1}^{n}({y}_{i}-\eta_{1}(\mathbf{y}))^{2}$.
\end{itemize}
For the observed data series, $\y$, we simulate $n=100$ observed data points from a normal random variable with mean $\theta=1$, variance $\sigma^2$, and take $\sigma^2=\{2,3\}$. For both ABC-AR and ABC-Reg, we again take $N=25,000$ simulated datasets generated according to ${z}^{j}_{i}\sim \mathcal{N}(\theta^{j},1)$, with $\theta^{j}\sim\mathcal{N}(0,25)$. For $d(\cdot,\cdot)$ we take the Euclidean norm and we again take the tolerance to be the 1\% quantile of the simulated distances. 

For the choice of the function $h(\theta)$, we consider $h(\theta)=(\theta^2,\theta^3)^\intercal$, so that $\hat{h}$ and $\tilde{h}$ represent the second and third posterior moments calculated under ABC-AR and ABC-Reg:
\begin{flalign*}
\hat{h}=\left(\int \theta^2 d\Pi_{\epsilon}[\theta|\eta(\mathbf{y})],\;\int \theta^3 d\Pi_{\epsilon}[\theta|\eta(\mathbf{y})]\right)^\intercal,\;\;
\tilde{h}=\left(\int \theta^2 d\tilde{\Pi}_{\epsilon}[\theta|\eta(\mathbf{y})],\;\int \theta^3 d\tilde{\Pi}_{\epsilon}[\theta|\eta(\mathbf{y})]\right)^\intercal.
\end{flalign*}Recall that under correct specification $\sqrt{n}\|\hat{h}-\tilde{h}\|=o_{P_{0}}(1)$. Thus, if the model was correctly specified, we would expect to realize a value of $\sqrt{n}\|\hat{h}-\tilde{h}\|$ that is relatively small.

To gauge the magnitude of  $\sqrt{n}\|\hat{h}-\tilde{h}\|$ under correct specification, we first simulate $b=1,\dots,B$ ``observed data'' series, where each series contains $n=100$ observations generated iid from a normal random variable with mean $\theta=1$ and $\sigma^2=1$. Then, using these series as our observed data, we run ABC-AR and ABC-Reg and calculate $\{\sqrt{n}\|\hat{h}_b-\tilde{h}_b\|\}_{b=1}^{B}$. For this experiment, we take $B=100$ replications. The cutoff value $t_n$ is then defined to be the $95\%$ empirical quantile.\footnote{The choice of $\theta$ used to simulate the pseudo ``observed data'' needed to find $t_n$ has little impact on the magnitude of the cutoff value in large samples. Therefore, it is possible to randomly choose a value of $\theta$ from the ABC-AR posterior without significantly altering the results presented here.}

Having obtained $t_n$, we now analyze the ability of this approach to detect misspecification. For the two different misspecified DGPs, which corresponds to $\sigma^2\in\{2,3\}$, we simulate 100 Monte Carlo replications and calculate $\sqrt{n}\|\hat{h}-\tilde{h}\|$ in each replication. The resulting sampling distributions of $\sqrt{n}\|\hat{h}-\tilde{h}\|$ are given in Figure \ref{fig_diag_reg1}. For comparison purposes, in Panel A of Figure \ref{fig_diag_reg1}, we also give the sampling distribution of $\sqrt{n}\|\hat{h}-\tilde{h}\|$ calculated under correct specification ($\sigma^2=1$). The results demonstrate that there are dramatic differences between the sampling distribution of $\sqrt{n}\|\hat{h}-\tilde{h}\|$ under correct and incorrect specification. 

Using the chosen cutoff value $t_n$, we detect misspecification 91\% of the time when $\sigma^2=2$ and 97\% when $\sigma^2=3$. Given these results, and the results in Figure \ref{fig_diag_reg1}, it is clear that this simple diagnostic allows us to detect model misspecification in this example.

\begin{figure}[h!]
	\centering
\scalebox{.70}{
\begin{tikzpicture}[x=1pt,y=1pt]
\definecolor{fillColor}{RGB}{255,255,255}
\path[use as bounding box,fill=fillColor,fill opacity=0.00] (0,0) rectangle (722.70,289.08);
\begin{scope}
\path[clip] ( 49.20, 61.20) rectangle (336.15,239.88);
\definecolor{fillColor}{RGB}{255,0,0}

\path[fill=fillColor] ( 73.11, 69.51) --
	(179.39, 69.51) --
	(179.39, 72.69) --
	( 73.11, 72.69) --
	cycle;
\definecolor{drawColor}{RGB}{0,0,0}

\path[draw=drawColor,line width= 1.2pt,line join=round] ( 73.11, 70.88) -- (179.39, 70.88);

\path[draw=drawColor,line width= 0.4pt,dash pattern=on 4pt off 4pt ,line join=round,line cap=round] (126.25, 67.82) -- (126.25, 69.51);

\path[draw=drawColor,line width= 0.4pt,dash pattern=on 4pt off 4pt ,line join=round,line cap=round] (126.25, 76.58) -- (126.25, 72.69);

\path[draw=drawColor,line width= 0.4pt,line join=round,line cap=round] ( 99.68, 67.82) -- (152.82, 67.82);

\path[draw=drawColor,line width= 0.4pt,line join=round,line cap=round] ( 99.68, 76.58) -- (152.82, 76.58);

\path[draw=drawColor,line width= 0.4pt,line join=round,line cap=round] ( 73.11, 69.51) --
	(179.39, 69.51) --
	(179.39, 72.69) --
	( 73.11, 72.69) --
	( 73.11, 69.51);
\definecolor{fillColor}{RGB}{67,110,238}

\path[fill=fillColor] (205.96, 88.91) --
	(312.24, 88.91) --
	(312.24,157.89) --
	(205.96,157.89) --
	cycle;

\path[draw=drawColor,line width= 1.2pt,line join=round] (205.96,111.58) -- (312.24,111.58);

\path[draw=drawColor,line width= 0.4pt,dash pattern=on 4pt off 4pt ,line join=round,line cap=round] (259.10, 68.76) -- (259.10, 88.91);

\path[draw=drawColor,line width= 0.4pt,dash pattern=on 4pt off 4pt ,line join=round,line cap=round] (259.10,233.26) -- (259.10,157.89);

\path[draw=drawColor,line width= 0.4pt,line join=round,line cap=round] (232.53, 68.76) -- (285.67, 68.76);

\path[draw=drawColor,line width= 0.4pt,line join=round,line cap=round] (232.53,233.26) -- (285.67,233.26);

\path[draw=drawColor,line width= 0.4pt,line join=round,line cap=round] (205.96, 88.91) --
	(312.24, 88.91) --
	(312.24,157.89) --
	(205.96,157.89) --
	(205.96, 88.91);
\end{scope}
\begin{scope}
\path[clip] (  0.00,  0.00) rectangle (722.70,289.08);
\definecolor{drawColor}{RGB}{0,0,0}

\path[draw=drawColor,line width= 0.4pt,line join=round,line cap=round] (126.25, 61.20) -- (259.10, 61.20);

\path[draw=drawColor,line width= 0.4pt,line join=round,line cap=round] (126.25, 61.20) -- (126.25, 55.20);

\path[draw=drawColor,line width= 0.4pt,line join=round,line cap=round] (259.10, 61.20) -- (259.10, 55.20);

\node[text=drawColor,anchor=base,inner sep=0pt, outer sep=0pt, scale=  1.00] at (126.25, 39.60) {$\sigma^2$=1};

\node[text=drawColor,anchor=base,inner sep=0pt, outer sep=0pt, scale=  1.00] at (259.10, 39.60) {$\sigma^2$=2};

\path[draw=drawColor,line width= 0.4pt,line join=round,line cap=round] ( 49.20, 67.52) -- ( 49.20,220.36);

\path[draw=drawColor,line width= 0.4pt,line join=round,line cap=round] ( 49.20, 67.52) -- ( 43.20, 67.52);

\path[draw=drawColor,line width= 0.4pt,line join=round,line cap=round] ( 49.20, 98.09) -- ( 43.20, 98.09);

\path[draw=drawColor,line width= 0.4pt,line join=round,line cap=round] ( 49.20,128.66) -- ( 43.20,128.66);

\path[draw=drawColor,line width= 0.4pt,line join=round,line cap=round] ( 49.20,159.22) -- ( 43.20,159.22);

\path[draw=drawColor,line width= 0.4pt,line join=round,line cap=round] ( 49.20,189.79) -- ( 43.20,189.79);

\path[draw=drawColor,line width= 0.4pt,line join=round,line cap=round] ( 49.20,220.36) -- ( 43.20,220.36);

\node[text=drawColor,rotate= 90.00,anchor=base,inner sep=0pt, outer sep=0pt, scale=  1.00] at ( 34.80, 67.52) {0};

\node[text=drawColor,rotate= 90.00,anchor=base,inner sep=0pt, outer sep=0pt, scale=  1.00] at ( 34.80, 98.09) {5};

\node[text=drawColor,rotate= 90.00,anchor=base,inner sep=0pt, outer sep=0pt, scale=  1.00] at ( 34.80,128.66) {10};

\node[text=drawColor,rotate= 90.00,anchor=base,inner sep=0pt, outer sep=0pt, scale=  1.00] at ( 34.80,159.22) {15};

\node[text=drawColor,rotate= 90.00,anchor=base,inner sep=0pt, outer sep=0pt, scale=  1.00] at ( 34.80,189.79) {20};

\node[text=drawColor,rotate= 90.00,anchor=base,inner sep=0pt, outer sep=0pt, scale=  1.00] at ( 34.80,220.36) {25};
\end{scope}
\begin{scope}
\path[clip] (  0.00,  0.00) rectangle (361.35,289.08);
\definecolor{drawColor}{RGB}{0,0,0}

\node[text=drawColor,anchor=base,inner sep=0pt, outer sep=0pt, scale=  1.20] at (192.68,260.34) {\bfseries Panel A};
\end{scope}
\begin{scope}
\path[clip] (  0.00,  0.00) rectangle (722.70,289.08);
\definecolor{drawColor}{RGB}{0,0,0}

\path[draw=drawColor,line width= 0.4pt,line join=round,line cap=round] ( 49.20, 61.20) --
	(336.15, 61.20) --
	(336.15,239.88) --
	( 49.20,239.88) --
	( 49.20, 61.20);
\end{scope}
\begin{scope}
\path[clip] (410.55, 61.20) rectangle (697.50,239.88);
\definecolor{fillColor}{RGB}{67,110,238}

\path[fill=fillColor] (434.46, 72.39) --
	(540.74, 72.39) --
	(540.74, 88.03) --
	(434.46, 88.03) --
	cycle;
\definecolor{drawColor}{RGB}{0,0,0}

\path[draw=drawColor,line width= 1.2pt,line join=round] (434.46, 77.53) -- (540.74, 77.53);

\path[draw=drawColor,line width= 0.4pt,dash pattern=on 4pt off 4pt ,line join=round,line cap=round] (487.60, 67.82) -- (487.60, 72.39);

\path[draw=drawColor,line width= 0.4pt,dash pattern=on 4pt off 4pt ,line join=round,line cap=round] (487.60,105.13) -- (487.60, 88.03);

\path[draw=drawColor,line width= 0.4pt,line join=round,line cap=round] (461.03, 67.82) -- (514.17, 67.82);

\path[draw=drawColor,line width= 0.4pt,line join=round,line cap=round] (461.03,105.13) -- (514.17,105.13);

\path[draw=drawColor,line width= 0.4pt,line join=round,line cap=round] (434.46, 72.39) --
	(540.74, 72.39) --
	(540.74, 88.03) --
	(434.46, 88.03) --
	(434.46, 72.39);
\definecolor{fillColor}{RGB}{0,255,0}

\path[fill=fillColor] (567.31, 78.81) --
	(673.59, 78.81) --
	(673.59,142.56) --
	(567.31,142.56) --
	cycle;

\path[draw=drawColor,line width= 1.2pt,line join=round] (567.31, 90.18) -- (673.59, 90.18);

\path[draw=drawColor,line width= 0.4pt,dash pattern=on 4pt off 4pt ,line join=round,line cap=round] (620.45, 68.19) -- (620.45, 78.81);

\path[draw=drawColor,line width= 0.4pt,dash pattern=on 4pt off 4pt ,line join=round,line cap=round] (620.45,233.26) -- (620.45,142.56);

\path[draw=drawColor,line width= 0.4pt,line join=round,line cap=round] (593.88, 68.19) -- (647.02, 68.19);

\path[draw=drawColor,line width= 0.4pt,line join=round,line cap=round] (593.88,233.26) -- (647.02,233.26);

\path[draw=drawColor,line width= 0.4pt,line join=round,line cap=round] (567.31, 78.81) --
	(673.59, 78.81) --
	(673.59,142.56) --
	(567.31,142.56) --
	(567.31, 78.81);
\end{scope}
\begin{scope}
\path[clip] (  0.00,  0.00) rectangle (722.70,289.08);
\definecolor{drawColor}{RGB}{0,0,0}

\path[draw=drawColor,line width= 0.4pt,line join=round,line cap=round] (487.60, 61.20) -- (620.45, 61.20);

\path[draw=drawColor,line width= 0.4pt,line join=round,line cap=round] (487.60, 61.20) -- (487.60, 55.20);

\path[draw=drawColor,line width= 0.4pt,line join=round,line cap=round] (620.45, 61.20) -- (620.45, 55.20);

\node[text=drawColor,anchor=base,inner sep=0pt, outer sep=0pt, scale=  1.00] at (487.60, 39.60) {$\sigma^2$=2};

\node[text=drawColor,anchor=base,inner sep=0pt, outer sep=0pt, scale=  1.00] at (620.45, 39.60) {$\sigma^2$=3};

\path[draw=drawColor,line width= 0.4pt,line join=round,line cap=round] (410.55, 67.54) -- (410.55,233.92);

\path[draw=drawColor,line width= 0.4pt,line join=round,line cap=round] (410.55, 67.54) -- (404.55, 67.54);

\path[draw=drawColor,line width= 0.4pt,line join=round,line cap=round] (410.55, 95.27) -- (404.55, 95.27);

\path[draw=drawColor,line width= 0.4pt,line join=round,line cap=round] (410.55,123.00) -- (404.55,123.00);

\path[draw=drawColor,line width= 0.4pt,line join=round,line cap=round] (410.55,150.73) -- (404.55,150.73);

\path[draw=drawColor,line width= 0.4pt,line join=round,line cap=round] (410.55,178.46) -- (404.55,178.46);

\path[draw=drawColor,line width= 0.4pt,line join=round,line cap=round] (410.55,206.19) -- (404.55,206.19);

\path[draw=drawColor,line width= 0.4pt,line join=round,line cap=round] (410.55,233.92) -- (404.55,233.92);

\node[text=drawColor,rotate= 90.00,anchor=base,inner sep=0pt, outer sep=0pt, scale=  1.00] at (396.15, 67.54) {0};

\node[text=drawColor,rotate= 90.00,anchor=base,inner sep=0pt, outer sep=0pt, scale=  1.00] at (396.15, 95.27) {20};

\node[text=drawColor,rotate= 90.00,anchor=base,inner sep=0pt, outer sep=0pt, scale=  1.00] at (396.15,123.00) {40};

\node[text=drawColor,rotate= 90.00,anchor=base,inner sep=0pt, outer sep=0pt, scale=  1.00] at (396.15,150.73) {60};

\node[text=drawColor,rotate= 90.00,anchor=base,inner sep=0pt, outer sep=0pt, scale=  1.00] at (396.15,178.46) {80};

\node[text=drawColor,rotate= 90.00,anchor=base,inner sep=0pt, outer sep=0pt, scale=  1.00] at (396.15,206.19) {100};

\node[text=drawColor,rotate= 90.00,anchor=base,inner sep=0pt, outer sep=0pt, scale=  1.00] at (396.15,233.92) {120};
\end{scope}
\begin{scope}
\path[clip] (361.35,  0.00) rectangle (722.70,289.08);
\definecolor{drawColor}{RGB}{0,0,0}

\node[text=drawColor,anchor=base,inner sep=0pt, outer sep=0pt, scale=  1.20] at (554.02,260.34) {\bfseries Panel B};
\end{scope}
\begin{scope}
\path[clip] (  0.00,  0.00) rectangle (722.70,289.08);
\definecolor{drawColor}{RGB}{0,0,0}

\path[draw=drawColor,line width= 0.4pt,line join=round,line cap=round] (410.55, 61.20) --
	(697.50, 61.20) --
	(697.50,239.88) --
	(410.55,239.88) --
	(410.55, 61.20);
\end{scope}
\end{tikzpicture}}
\caption{Monte Carlo sampling distribution of $\sqrt{n}\|\hat{h}-\tilde{h}\|$ in the normal example for $\sigma^{2}\in\{1,2,3\}$, where $h=(\theta^2,\theta^3)^{\intercal}$. Recall that $\sigma^2=1$ corresponds to correct specification.} 
	\label{fig_diag_reg1}
\end{figure}

{It is important to realize that the behavior of $\sqrt{n}\|\hat{h}-\tilde{h}\|$ observed in Figure \ref{fig_diag_reg1} is not a function of the specific $h(\theta)$ we have chosen but is driven by the differences in the ABC-AR and ABC-Reg posteriors under model misspecification. In particular, the results of Corollary \ref{cor1} and \ref{th:asympAdj} imply that similar behavior will be observed for virtually any well-behaved function $\theta\mapsto h(\theta)$.

To illustrate this fact, we repeat the example above using the alternative function $h(\theta)=\theta$. For this experiment, we use  the exact same datasets as in the above example, and implement the same procedure to find the cutoff value $t_n$, but now in the case where $h(\theta)=\theta$. Across the Monte Carlo replications, this approach detects misspecification 96\% of the time when $\sigma^2=2$ and 99\% when $\sigma^2=3$. Similarly, Figure \ref{fig_diag_reg2} plots the sampling distributions of $\sqrt{n}\|\hat{h}-\tilde{h}\|$ when $h(\theta)=\theta$. While the scales are different to those in Figure \ref{fig_diag_reg1}, the results are qualitatively the same: there are significant differences between the sampling distribution of $\sqrt{n}\|\hat{h}-\tilde{h}\|$ under correct and incorrect specification. We reiterate that this result is not surprising since, under our results, we know that ABC-AR and ARC-Reg concentrate posterior mass on different values and therefore any sufficiently smooth function $h(\theta)$ will concentrate onto different values under ABC-AR and ABC-Reg.\hfill$\square$
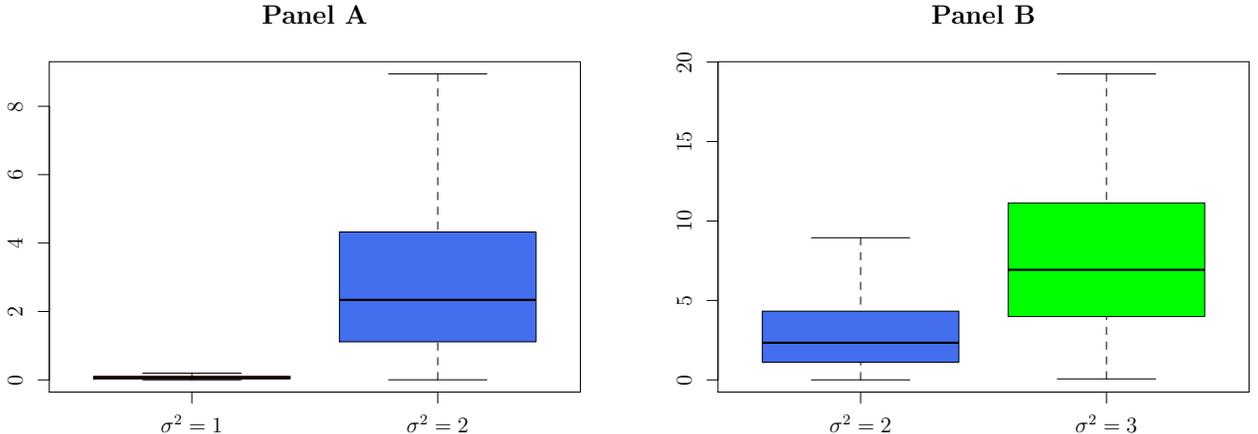
\begin{figure}[H]
	\centering
	\scalebox{.70}{
\begin{tikzpicture}[x=1pt,y=1pt]
\definecolor{fillColor}{RGB}{255,255,255}
\path[use as bounding box,fill=fillColor,fill opacity=0.00] (0,0) rectangle (722.70,289.08);
\begin{scope}
\path[clip] ( 49.20, 61.20) rectangle (336.15,239.88);
\definecolor{fillColor}{RGB}{255,0,0}

\path[fill=fillColor] ( 73.11, 68.40) --
	(179.39, 68.40) --
	(179.39, 69.67) --
	( 73.11, 69.67) --
	cycle;
\definecolor{drawColor}{RGB}{0,0,0}

\path[draw=drawColor,line width= 1.2pt,line join=round] ( 73.11, 68.87) -- (179.39, 68.87);

\path[draw=drawColor,line width= 0.4pt,dash pattern=on 4pt off 4pt ,line join=round,line cap=round] (126.25, 67.82) -- (126.25, 68.40);

\path[draw=drawColor,line width= 0.4pt,dash pattern=on 4pt off 4pt ,line join=round,line cap=round] (126.25, 71.46) -- (126.25, 69.67);

\path[draw=drawColor,line width= 0.4pt,line join=round,line cap=round] ( 99.68, 67.82) -- (152.82, 67.82);

\path[draw=drawColor,line width= 0.4pt,line join=round,line cap=round] ( 99.68, 71.46) -- (152.82, 71.46);

\path[draw=drawColor,line width= 0.4pt,line join=round,line cap=round] ( 73.11, 68.40) --
	(179.39, 68.40) --
	(179.39, 69.67) --
	( 73.11, 69.67) --
	( 73.11, 68.40);
\definecolor{fillColor}{RGB}{67,110,238}

\path[fill=fillColor] (205.96, 88.39) --
	(312.24, 88.39) --
	(312.24,147.82) --
	(205.96,147.82) --
	cycle;

\path[draw=drawColor,line width= 1.2pt,line join=round] (205.96,111.04) -- (312.24,111.04);

\path[draw=drawColor,line width= 0.4pt,dash pattern=on 4pt off 4pt ,line join=round,line cap=round] (259.10, 67.85) -- (259.10, 88.39);

\path[draw=drawColor,line width= 0.4pt,dash pattern=on 4pt off 4pt ,line join=round,line cap=round] (259.10,233.26) -- (259.10,147.82);

\path[draw=drawColor,line width= 0.4pt,line join=round,line cap=round] (232.53, 67.85) -- (285.67, 67.85);

\path[draw=drawColor,line width= 0.4pt,line join=round,line cap=round] (232.53,233.26) -- (285.67,233.26);

\path[draw=drawColor,line width= 0.4pt,line join=round,line cap=round] (205.96, 88.39) --
	(312.24, 88.39) --
	(312.24,147.82) --
	(205.96,147.82) --
	(205.96, 88.39);
\end{scope}
\begin{scope}
\path[clip] (  0.00,  0.00) rectangle (722.70,289.08);
\definecolor{drawColor}{RGB}{0,0,0}

\path[draw=drawColor,line width= 0.4pt,line join=round,line cap=round] (126.25, 61.20) -- (259.10, 61.20);

\path[draw=drawColor,line width= 0.4pt,line join=round,line cap=round] (126.25, 61.20) -- (126.25, 55.20);

\path[draw=drawColor,line width= 0.4pt,line join=round,line cap=round] (259.10, 61.20) -- (259.10, 55.20);

\node[text=drawColor,anchor=base,inner sep=0pt, outer sep=0pt, scale=  1.00] at (126.25, 39.60) {$\sigma^2=1$};

\node[text=drawColor,anchor=base,inner sep=0pt, outer sep=0pt, scale=  1.00] at (259.10, 39.60) {$\sigma^2=2$};

\path[draw=drawColor,line width= 0.4pt,line join=round,line cap=round] ( 49.20, 67.81) -- ( 49.20,215.78);

\path[draw=drawColor,line width= 0.4pt,line join=round,line cap=round] ( 49.20, 67.81) -- ( 43.20, 67.81);

\path[draw=drawColor,line width= 0.4pt,line join=round,line cap=round] ( 49.20,104.80) -- ( 43.20,104.80);

\path[draw=drawColor,line width= 0.4pt,line join=round,line cap=round] ( 49.20,141.79) -- ( 43.20,141.79);

\path[draw=drawColor,line width= 0.4pt,line join=round,line cap=round] ( 49.20,178.79) -- ( 43.20,178.79);

\path[draw=drawColor,line width= 0.4pt,line join=round,line cap=round] ( 49.20,215.78) -- ( 43.20,215.78);

\node[text=drawColor,rotate= 90.00,anchor=base,inner sep=0pt, outer sep=0pt, scale=  1.00] at ( 34.80, 67.81) {0};

\node[text=drawColor,rotate= 90.00,anchor=base,inner sep=0pt, outer sep=0pt, scale=  1.00] at ( 34.80,104.80) {2};

\node[text=drawColor,rotate= 90.00,anchor=base,inner sep=0pt, outer sep=0pt, scale=  1.00] at ( 34.80,141.79) {4};

\node[text=drawColor,rotate= 90.00,anchor=base,inner sep=0pt, outer sep=0pt, scale=  1.00] at ( 34.80,178.79) {6};

\node[text=drawColor,rotate= 90.00,anchor=base,inner sep=0pt, outer sep=0pt, scale=  1.00] at ( 34.80,215.78) {8};
\end{scope}
\begin{scope}
\path[clip] (  0.00,  0.00) rectangle (361.35,289.08);
\definecolor{drawColor}{RGB}{0,0,0}

\node[text=drawColor,anchor=base,inner sep=0pt, outer sep=0pt, scale=  1.20] at (192.68,260.34) {\bfseries Panel A};
\end{scope}
\begin{scope}
\path[clip] (  0.00,  0.00) rectangle (722.70,289.08);
\definecolor{drawColor}{RGB}{0,0,0}

\path[draw=drawColor,line width= 0.4pt,line join=round,line cap=round] ( 49.20, 61.20) --
	(336.15, 61.20) --
	(336.15,239.88) --
	( 49.20,239.88) --
	( 49.20, 61.20);
\end{scope}
\begin{scope}
\path[clip] (410.55, 61.20) rectangle (697.50,239.88);
\definecolor{fillColor}{RGB}{67,110,238}

\path[fill=fillColor] (434.46, 77.36) --
	(540.74, 77.36) --
	(540.74,104.97) --
	(434.46,104.97) --
	cycle;
\definecolor{drawColor}{RGB}{0,0,0}

\path[draw=drawColor,line width= 1.2pt,line join=round] (434.46, 87.89) -- (540.74, 87.89);

\path[draw=drawColor,line width= 0.4pt,dash pattern=on 4pt off 4pt ,line join=round,line cap=round] (487.60, 67.82) -- (487.60, 77.36);

\path[draw=drawColor,line width= 0.4pt,dash pattern=on 4pt off 4pt ,line join=round,line cap=round] (487.60,144.67) -- (487.60,104.97);

\path[draw=drawColor,line width= 0.4pt,line join=round,line cap=round] (461.03, 67.82) -- (514.17, 67.82);

\path[draw=drawColor,line width= 0.4pt,line join=round,line cap=round] (461.03,144.67) -- (514.17,144.67);

\path[draw=drawColor,line width= 0.4pt,line join=round,line cap=round] (434.46, 77.36) --
	(540.74, 77.36) --
	(540.74,104.97) --
	(434.46,104.97) --
	(434.46, 77.36);
\definecolor{fillColor}{RGB}{0,255,0}

\path[fill=fillColor] (567.31,102.06) --
	(673.59,102.06) --
	(673.59,163.51) --
	(567.31,163.51) --
	cycle;

\path[draw=drawColor,line width= 1.2pt,line join=round] (567.31,127.33) -- (673.59,127.33);

\path[draw=drawColor,line width= 0.4pt,dash pattern=on 4pt off 4pt ,line join=round,line cap=round] (620.45, 68.34) -- (620.45,102.06);

\path[draw=drawColor,line width= 0.4pt,dash pattern=on 4pt off 4pt ,line join=round,line cap=round] (620.45,233.26) -- (620.45,163.51);

\path[draw=drawColor,line width= 0.4pt,line join=round,line cap=round] (593.88, 68.34) -- (647.02, 68.34);

\path[draw=drawColor,line width= 0.4pt,line join=round,line cap=round] (593.88,233.26) -- (647.02,233.26);

\path[draw=drawColor,line width= 0.4pt,line join=round,line cap=round] (567.31,102.06) --
	(673.59,102.06) --
	(673.59,163.51) --
	(567.31,163.51) --
	(567.31,102.06);
\end{scope}
\begin{scope}
\path[clip] (  0.00,  0.00) rectangle (722.70,289.08);
\definecolor{drawColor}{RGB}{0,0,0}

\path[draw=drawColor,line width= 0.4pt,line join=round,line cap=round] (487.60, 61.20) -- (620.45, 61.20);

\path[draw=drawColor,line width= 0.4pt,line join=round,line cap=round] (487.60, 61.20) -- (487.60, 55.20);

\path[draw=drawColor,line width= 0.4pt,line join=round,line cap=round] (620.45, 61.20) -- (620.45, 55.20);

\node[text=drawColor,anchor=base,inner sep=0pt, outer sep=0pt, scale=  1.00] at (487.60, 39.60) {$\sigma^2=2$};

\node[text=drawColor,anchor=base,inner sep=0pt, outer sep=0pt, scale=  1.00] at (620.45, 39.60) {$\sigma^2=3$};

\path[draw=drawColor,line width= 0.4pt,line join=round,line cap=round] (410.55, 67.80) -- (410.55,239.67);

\path[draw=drawColor,line width= 0.4pt,line join=round,line cap=round] (410.55, 67.80) -- (404.55, 67.80);

\path[draw=drawColor,line width= 0.4pt,line join=round,line cap=round] (410.55,110.77) -- (404.55,110.77);

\path[draw=drawColor,line width= 0.4pt,line join=round,line cap=round] (410.55,153.74) -- (404.55,153.74);

\path[draw=drawColor,line width= 0.4pt,line join=round,line cap=round] (410.55,196.70) -- (404.55,196.70);

\path[draw=drawColor,line width= 0.4pt,line join=round,line cap=round] (410.55,239.67) -- (404.55,239.67);

\node[text=drawColor,rotate= 90.00,anchor=base,inner sep=0pt, outer sep=0pt, scale=  1.00] at (396.15, 67.80) {0};

\node[text=drawColor,rotate= 90.00,anchor=base,inner sep=0pt, outer sep=0pt, scale=  1.00] at (396.15,110.77) {5};

\node[text=drawColor,rotate= 90.00,anchor=base,inner sep=0pt, outer sep=0pt, scale=  1.00] at (396.15,153.74) {10};

\node[text=drawColor,rotate= 90.00,anchor=base,inner sep=0pt, outer sep=0pt, scale=  1.00] at (396.15,196.70) {15};

\node[text=drawColor,rotate= 90.00,anchor=base,inner sep=0pt, outer sep=0pt, scale=  1.00] at (396.15,239.67) {20};
\end{scope}
\begin{scope}
\path[clip] (361.35,  0.00) rectangle (722.70,289.08);
\definecolor{drawColor}{RGB}{0,0,0}

\node[text=drawColor,anchor=base,inner sep=0pt, outer sep=0pt, scale=  1.20] at (554.02,260.34) {\bfseries Panel B};
\end{scope}
\begin{scope}
\path[clip] (  0.00,  0.00) rectangle (722.70,289.08);
\definecolor{drawColor}{RGB}{0,0,0}

\path[draw=drawColor,line width= 0.4pt,line join=round,line cap=round] (410.55, 61.20) --
	(697.50, 61.20) --
	(697.50,239.88) --
	(410.55,239.88) --
	(410.55, 61.20);
\end{scope}
\end{tikzpicture}}
	\caption{Monte Carlo sampling distribution of $\sqrt{n}\|\hat{h}-\tilde{h}\|$ in the normal example for $\sigma^{2}\in\{1,2,3\}$, where $h(\theta)=\theta$. Recall that $\sigma^2=1$ corresponds to correct specification.} 
	\label{fig_diag_reg2}
\end{figure}

Lastly, we note that the choice of $\theta$ used to construct $t_n$ will not significantly alter the reported results. To demonstrate this, we rerun the diagnostic procedure at two additional values of $\theta$ used to construct $t_n$. If we use the value $\theta=0$ to obtain $t_n$, in the case where $h(\theta)=(\theta^2,\theta^3)^{\intercal}$, the procedure would lead us to conclude in favor of misspecification 100\% of the time for $\sigma^2=2$ and $\sigma^2=3$, while using the function $h(\theta)=\theta$ we would conclude in favor of misspecification 97\% of the time under $\sigma^2=2$ and 100\% under $\sigma^2=3$. Similarly, if we consider the value $\theta=2$ and construct $t_n$, in the case where $h(\theta)=(\theta^2,\theta^3)^{\intercal}$, we would conclude in favor of misspecification 83\% of the time for $\sigma^2=2$ and 92\% under $\sigma^2=3$, while taking  $h(\theta)=\theta$ we would conclude in favor of misspecification 95\% of the time under $\sigma^2=2$ and 98\% under $\sigma^2=3$.

\section{Discussion}

{At first glance, ABC techniques seem less than appropriate in misspecified settings as they {rely heavily on the assumed models capacity to reproduce the features of the observed dataset.} Furthermore, they scale the proximity (or tolerance) in terms of the actual data as a percentile, henceforth providing a relative as opposed to absolute measure of proximity. This paper demonstrates that ABC is indeed open to poor performances {in these settings}, with the performance seemingly worse the more involved the version of ABC. This pattern is not immensely surprising in that more complex methods try to extract more information from the simulated data. What we find more exciting, and of particular interest, are the convergence results for the more basic ABC versions, provided an identifiability constraint on the summary statistic {is satisfied}. The role of the tolerance {sequence} also appears to be more crucial than in the well-specified case, in connection with the fact that there is a minimal non-zero limit for this tolerance. {Furthermore, we demonstrate that} post-processing the ABC output by local regression may lead to {poor inference in misspecified models}, and we propose an alternative approach that is less sensitive to the correctness of the model specification. At this stage, it is unclear if other post-processing ABC approaches, such as, for instance, the marginal adjustment approach (\citealp{nott2014approximate}) or the recalibration approach (\citealp{rodrigues2018recalibration}), will perform similarly to local regression post-processing under model misspecification, and more research on this front is necessary to obtain conclusive results. 
	
In addition to the above, we demonstrate that, rather naturally, the discrepancies between the different ABC versions can be exploited to {detect model} misspecification, despite the inherent difficulty in separating Monte Carlo variability from a genuine difference. Both the examples in the paper and in the appendix illustrate that detection is achievable, since the reference tables produced by ABC {can be used to} automatically calibrate the expected differences {between the different versions of ABC}.

A potential direction of {future} research on misspecification not investigated in this paper would be to connect it with generative adversarial networks (GAN, \citealp{goodfellow:etal:2014}) since this machine-learning technique aims both at learning a data generating mechanism (the generative network) and at separating the actual data from the generated one (the discriminating network). While the proximity between ABC and GANs has already been noted, there is however little known from a statistical viewpoint in this area about misspecification properties. 

Concerning other versions of ABC based on a machine-learning approach to the intractable likelihood, for example, {such as in \cite{papamakarios2016fast}}, our intuition is that, {similar to ABC based on local nonlinear regression adjustment,} the massive and unstructured dependence of the calibration on the model will similarly yield poor performance under misspecification. Apart from a GAN strategy, {at this stage we} fail to see a robust solution to this difficulty.} 

\singlespacing
\bibliographystyle{apalike}
\bibliography{refs_mispec_jude}

\appendix

\section{Proofs of Main Results}  \label{sec:sup:proof}
This section contains the proofs of the theoretical results given in the main text. 
\subsection*{Proof of Theorem \ref{th:asymp0}}

This theorem is an adaptation of \citet{FMRR2016}. Let $\delta_{n}\geq M_{n}( \epsilon_n- \epsilon^*) \geq 3M_nv^{-1}_{0,n}$, then  $P_{0}(\Omega_{d})=1+o(1)$  for $\Omega_{d}:=\{\mathbf{y}:d_{}(\eta(\mathbf{y}),b_{0})\leq \delta_n/2\}$. Assume that $\mathbf{y} \in \Omega_d$.
Consider the event 
$$A_{d}(\delta_n):=\left\{(\mathbf{z},\theta):\{d(\eta(\mathbf{z}), \eta (\mathbf{y}) ) \leq \epsilon_n\}\cap\{d(b(\theta),b_{0} ) \geq \epsilon^*+ \delta_n \}\right\}.$$

Note that, by definition $d(b(\theta),b_{0})\geq\epsilon^{*}$, with $\epsilon^{*}>0$. For all $(\mathbf{z},\theta)\in A_{d}(\delta_n)$ and if $\mathbf{y} \in \Omega_d$, 
\begin{flalign*}
\delta_{n}<d(b(\theta),b_{0})-\epsilon^{*}&\leq d(b(\theta),\eta(\mathbf{z}))+d(\eta(\mathbf{z}),\eta(\mathbf{y}) ) +  d(\eta(\mathbf{y}), b_{0})-\epsilon^{*}\\
&\leq d(b(\theta),\eta(\mathbf{z})) + \epsilon_n - \epsilon^* + \delta_{n}/2
\end{flalign*}
so that $ \delta_{n} \leq 4 d(b(\theta),\eta(\mathbf{z})).$
This implies in particular that 
\begin{equation}\label{pr-eq1}
\begin{split}
\text{Pr}(A_{d}(\delta_{n}))  &=  \int_{\{d(b(\theta),b_{0} ) \geq \epsilon^*+ \delta_{n} \}} P_{\theta}\left[d(\eta(\mathbf{z}), \eta (\mathbf{y}) ) \leq \epsilon_n \right] d\Pi(\theta) \\
& \leq \int_\Theta P_\theta \left[d(b(\theta),\eta(\mathbf{z})) \geq \delta_{n}/4 \right] d\Pi(\theta) .
\end{split}
\end{equation} 
In case (i) of polynomial tails, 
\begin{equation}\label{pr-eqcasei}
\begin{split}
\text{Pr}(A_{d}(\delta_{n}))  
& \leq (v_{n}\delta_{n})^{-\kappa}  \int_\Theta c(\theta)  d\Pi(\theta) = o(1)
\end{split}
\end{equation} 
as soon as $v_{n} \delta_{n} \rightarrow + \infty$, or in case (ii) of exponential tails
\begin{equation}\label{pr-eqcaseii}
\begin{split}
\text{Pr}(A_{d}(\delta_{n}))  
& \leq C  e^{-c (\delta_{n} v_{n})^\tau} .
\end{split}
\end{equation} 
Moreover, we can bound from below 
\begin{equation*}
\begin{split}
\alpha_{n} &=   \int_\Theta P_{\theta}\left[d(\eta(\mathbf{z}), \eta (\mathbf{y}) ) \leq \epsilon_n \right] d\Pi(\theta)\\
\end{split}
\end{equation*}
Note that on $\{d(\eta(\mathbf{z}),b(\theta)) \leq M v_{n}^{-1}/2 \}\cap \Omega_d$
\begin{equation*}
d(\eta(\mathbf{y}),\eta(\mathbf{z})) \leq d(\eta(\mathbf{z}),b(\theta))+d(\eta(\mathbf{y}),b_{0})+ d(b(\theta),b_{0}) 
\leq v_{0,T}^{-1} + Mv_{n}^{-1}/2 + d(b(\theta),b_{0})  \leq \epsilon_n
\end{equation*}
as soon as $\epsilon^*\leq d(b(\theta),b_{0}) \leq \epsilon_n-v_{0,T}^{-1} + Mv_{n}^{-1}/2 $. Since $\epsilon_n - \epsilon^* \geq v_{0,T}^{-1} + Mv_{n}^{-1}$, on $\Omega_d$, 
\begin{equation*}
\begin{split}
\int_{\Theta }P_{{\theta }}\left( d({\eta }(%
\mathbf{z}),{\eta }(\mathbf{y}))\leq \epsilon_n \right) d\Pi ({\theta })& \geq \int_{ d(b(\theta),b_{0}) \leq (\epsilon_n - \epsilon^*)/4 \vee v_{n}^{-1}M/2 }( 1 - P_\theta \left( d(\eta(\mathbf{z}), b(\theta) ) \geq M v_{n}^{-1}/2 \right) d\Pi(\theta) \\
&\geq \int_{ d(b(\theta),b_{0}) \leq (\epsilon_n - \epsilon^*)/4 \vee v_{n}^{-1}M/2 }\left( 1 - \frac{ c(\theta)2^\kappa }{ M^\kappa} \right)d\Pi(\theta) \\
&\gtrsim (\epsilon_n - \epsilon^*)^D \vee v_{n}^{-D}\gtrsim (\epsilon_n - \epsilon^*)^D
\end{split}
\end{equation*}%
in case (i) of \textbf{[A1]}, under \textbf{[A2]}. If case (ii) of \textbf{[A1]} holds, under \textbf{[A2]}, we have 
\begin{equation*}
\begin{split}
\int_{\Theta }P_{{\theta }}\left( d({\eta }(%
\mathbf{z}),{\eta }(\mathbf{y}))\leq \epsilon_n \right) d\Pi ({\theta })
&\geq \int_{ d(b(\theta),b_{0}) \leq (\epsilon_n - \epsilon^*)/4 \vee v_{n}^{-1}M/2 }\left( 1 -  c(\theta)e^{-h_\theta(M/2)} \right)d\Pi(\theta) \\
&\gtrsim (\epsilon_n - \epsilon^*)^D 
\end{split}
\end{equation*}%
Combining these two inequality with the upper bounds \eqref{pr-eqcasei} or \eqref{pr-eqcaseii} leads to 
$$\Pi_{\epsilon}\left[d(b(\theta), b_{0})\geq  \epsilon^*+\delta_{n} |\eta(\mathbf{y})\right] \lesssim (\epsilon_n - \epsilon^*)^{-D}(v_{n}\delta_{n})^{-\kappa},$$
in case (i) and 
$$\Pi_{\epsilon}\left[d(b(\theta), b_{0})\geq  \epsilon^*+\delta_{n} |\eta(\mathbf{y})\right] \lesssim (\epsilon_n - \epsilon^*)^{-D}e^{-c (\delta_{n} v_{n})^\tau} , $$
in case (ii). These are of order $o(1)$ if $\delta_{n}\geq M_{n} v_{n}^{-1} (\epsilon_n-\epsilon^*)^{-D/\kappa}$ in case (i), or if $\delta_{n}\geq M_{n} v_{n}^{-1} |\log(\epsilon_n-\epsilon^*)  |^{1/\tau}$ in case (ii).   

\subsection*{Proof of Corollary \ref{cor1}}
\begin{proof}
	Define $Q(\theta)=|d(b(\theta),b_{0})-d(b(\theta^*),b_{0})|$. From the continuity of $\theta\mapsto b(\theta)$ and the definition of $\theta^*$, for any $\delta>0$ there exists a $\gamma(\delta)>0$ such that $\inf_{\theta:d\{\theta,\theta^*\}>\delta}Q(\theta)\geq\gamma(\delta)>0.$ Then,
	\begin{flalign*}
	\Pi_{\epsilon}[d_{}(\theta,\theta^{*})>\delta |\eta(\mathbf{y})]\leq\Pi_{\epsilon}[|Q(\theta)-Q(\theta^{*})|>\gamma(\delta)|\eta(\mathbf{y})]&=\Pi_{\epsilon}[|d(b(\theta),b_{0})-d(b(\theta^*),b_{0})|>\gamma(\delta)|\eta(\mathbf{y})]\\&=\Pi_{\epsilon_{}}[d(b(\theta),b_{0})>\epsilon^{*}+\gamma(\delta)|\eta(\mathbf{y})].
	\end{flalign*}The result follows if $\Pi_{\epsilon}[|d_{}(b(\theta),b_{0})>\epsilon^{*}+\gamma(\delta)|\eta(\mathbf{y})]=o_{P_{0}}(1)$. For $\delta_{n}>0$ and $\delta_{n}=o(1)$ as defined in Theorem \ref{th:asymp0}, by the conclusion of Theorem \ref{th:asymp0}, the result follows once $\gamma(\delta)\geq\delta_{n}$. 
\end{proof}

\subsection*{Proof of Theorem \ref{normal_thm}}
\begin{proof}
	For the sake of simplicity and without loss of generality we write $v_n = \sqrt{n}$, $\tilde Z_n = \sqrt{n} (\eta(\mathbf z) - b(\theta) )$ and $\tilde Z_y = \sqrt{n} ( \eta(\mathbf y) - b_0) $. Denote by $ B_n(K) = \{ \|\theta - \theta^* \| \leq K\} $. Throughout the proof $C$ denotes a generic constant which may vary from line to line.     We have for all $\theta$ 
	\begin{equation*}
	\begin{split}
	P_\theta&  \left( \|\eta(\z) - \eta(\y)\|^2\leq \epsilon_n^2  \right)  = P_\theta \left( \| \tilde Z_n - \tilde Z_y +\sqrt{n}(b(\theta ) - b_0)\|^2 \leq n \epsilon_n^2  \right) \\
	& = P_\theta\left( \| \tilde Z_n\|^2 +2\langle Z_n , \sqrt{n}(b(\theta) - b_0)-\tilde Z_y \rangle \leq n[\epsilon_n^2 - \|b(\theta)- b_0 -\tilde Z_y/\sqrt{n}\|^2] \right) \\
	&= 
	P_\theta\left(\langle\tilde Z_n ,b(\theta) - b_0 \rangle   \leq \frac{\sqrt{n}[\epsilon_n^2 - \|b(\theta)- b_0 -\tilde Z_y/\sqrt{n}\|^2]}{2}  -\frac{  \|\tilde Z_n\|^2 - 2\langle\tilde Z_n, \tilde Z_y\rangle}{ 2\sqrt{n}} \right)
	\end{split}
	\end{equation*}
	Now on $\Omega_n =  \{ \| \tilde Z_y \| \leq M_n/2\}$ with $M_n$ a sequence going to infinity arbitrarily slowly and such that $M_n = o(n^{1/4})$,
	\begin{equation*}
	\begin{split} \sqrt{n}\|b(\theta)- b_0 -\tilde Z_y/\sqrt{n}\|^2& = \sqrt{n} \|b(\theta^*)- b_0\|^2 +  \frac{\sqrt{n}(\theta - \theta^*)^{\intercal} H^*(\theta - \theta^*)}{2 }   - 2 \langle  b(\theta^*)- b_0, \tilde Z_y \rangle   \\
	&+ O(M_n^2 /\sqrt{n} ) + O(\sqrt{n} \|\theta - \theta^*\|^3+ \|\theta - \theta^*\| M_n) 
	\end{split}
	\end{equation*}
	where $H^*$ is the second derivative of $\theta \mapsto \|b(\theta)- b_0\|^2$ at $\theta^*$, noting that the first derivative is equal to 0 at $\theta^*$. 
	Let  $\epsilon^* = \|b(\theta^*)- b_0\|$, $e' = (b(\theta^*)- b_0)$, and $\epsilon >0$. If $\|\theta - \theta^*\|\leq \epsilon$ and on the event $\Omega_n = \{ \|Z_y \| \leq M_n \}$ where $M_n^2  = o(\sqrt{n})$,  
	\begin{flalign}\label{deve:likeli}
	P_\theta&  \left( \|\eta(\mathbf z ) - \eta(\y)\|^2\leq \epsilon_n^2  \right) \nonumber\\
	&  \leq  P_\theta\left(\langle\tilde Z_n ,e' \rangle \leq \frac{\sqrt{n}[\epsilon_n^2 - (\epsilon^*)^2 - (1+C\epsilon)(\theta - \theta^*)^{\intercal} H^*(\theta- \theta^*) /2 ]}{2} + \langle\tilde Z_y ,e' \rangle  + C\epsilon + \| \theta - \theta^*\| M_n \right) \nonumber\\
	& + P_\theta \left( \|  Z_n\|^2 > \epsilon^2 \sqrt{n}/4  \right) \nonumber\\
	& \geq 
	P_\theta\left(\langle\tilde Z_n ,e' \rangle \leq \frac{\sqrt{n}[\epsilon_n^2 - (\epsilon^*)^2 - (1-C\epsilon)(\theta - \theta^*)^{\intercal} H^*(\theta- \theta^*)/2  ]}{2} + \langle\tilde Z_y ,e' \rangle  - C\epsilon - \| \theta - \theta^*\| M_n \right)\nonumber\\
	&  - P_\theta \left( \|  \tilde Z_n\|^2 > \epsilon^2 \sqrt{n}/4  \right) 
	\end{flalign}
	
	Consider the case where $\sqrt{n}(\epsilon_n^2 - (\epsilon^*)^2 ) \rightarrow 2c \in \mathbb{ R}$. 
	We split $\Theta$ into $\{ \| \theta - \theta^* \|^2 \sqrt{n } \leq M\}$, $\{ \epsilon > \| \theta  - \theta^*\| > M /n^{1/4}\}$ and $\{\epsilon \leq \| \theta  - \theta^*\| \}$, where $\epsilon $ is arbitrarily small.

	First if $\| \theta - \theta^* \|^2 \sqrt{n } \leq M$
	\begin{equation*}
	\begin{split}
	&  \frac{\sqrt{n}[\epsilon_n^2 - (\epsilon^*)^2 - (1-C\epsilon)(\theta - \theta^*)^{\intercal} H^*(\theta- \theta^*)/2  ]}{2}   + \langle\tilde Z_y ,e' \rangle + \frac{ M_n^2}{\sqrt{n} } + \| \theta - \theta^*\| M_n \\
	& \quad \leq c + \langle\tilde Z_y ,e' \rangle 
	-\frac{ (1-C\epsilon)\sqrt{n} (\theta - \theta^*)^{\intercal} H^*(\theta- \theta^*) }{4} + \epsilon 
	\end{split}
	\end{equation*}
	and 
	\begin{equation*}
	\begin{split}
	& \frac{\sqrt{n}[\epsilon_n^2 - (\epsilon^*)^2 - (1+C\epsilon)(\theta - \theta^*)^{\intercal} H^*(\theta- \theta^*)/2  ]}{2}   + \frac{ M^2}{\sqrt{n} } + \| \theta - \theta^*\| M_n \\
	& \quad \geq c 
	+ \langle\tilde Z_y ,e' \rangle -\frac{ (1+C\epsilon)\sqrt{n} (\theta - \theta^*)^{\intercal} H^*(\theta- \theta^*)}{4}  - \epsilon 
	\end{split}
	\end{equation*}
	Moreover, using assumption \textbf{[A5]},  
	$$\langle\tilde Z_n ,e' \rangle  = \sqrt{n} < \Sigma_n(\theta)^{-1} Z_n, e' \rangle =  \langle Z_n, \sqrt{n} \Sigma_n(\theta)^{-1} e' \rangle = \langle Z_n, A(\theta^*)  e' \rangle+ o( \|Z_n\|).$$
	We then have with $c' = c + \langle\tilde Z_y ,e' \rangle $, $x  = n^{1/4} (1-\epsilon)^{1/2}  ( \theta- \theta^* ) $ and   $\| \theta - \theta^* \| \leq M/n^{1/4}\leq u_0$ if $n$ is large enough
	\begin{equation*}
	\begin{split}
	P_\theta&  \left( \|\eta(\mathbf z ) - \eta(\y)\|^2\leq \epsilon_n^2  \right)  \leq P_\theta\left(\langle Z_n , A(\theta^*) e' \rangle \leq c' - \frac{x^{\intercal} H^* x }{ 4 }   + \epsilon \right)  +\epsilon\\
	& \leq \Phi \left( \frac{ c'+M\epsilon}{ \|A(\theta^*) e'\|  }  -\frac{x^{\intercal} H^* x }{ 4 \|A(\theta^*) e'\| }  \right) + \epsilon \\ 
	&+  \sup_{\|\theta' - \theta^*\|\leq u } \left| P_{\theta'}\left(\langle Z_n , A(\theta^*) e' \rangle \leq c' - \frac{x^{\intercal} H^* x}{ 4 }   + \epsilon \right) - \Phi \left( \frac{ c'+\epsilon}{ \|A(\theta^*) e'\|  }  -\frac{x^{\intercal} H^* x }{ 4 \|A(\theta^*) e'\| }  \right)\right| \\
	& \leq \Phi \left( \frac{ c'+\epsilon}{ \|A(\theta^*) e'\|  }  -\frac{x^{\intercal} H^* x }{ 4 \|A(\theta^*) e'\| }  \right)  +C \epsilon . 
	\end{split}
	\end{equation*}
	Similarly with $y = n^{1/4} (1+\epsilon)^{1/2} (\theta- \theta^*)$
	\begin{equation*}
	\begin{split}
	P_\theta&  \left( \|\eta(\mathbf z ) - \eta(\y)\|^2\leq \epsilon_n^2  \right) \\
	& \geq 
	\Phi \left( \frac{ c'-\epsilon}{ \|A(\theta^*) e'\|  }  -\frac{y^{\intercal} H^* y }{ 4 \|A(\theta^*) e'\| }  \right) -\epsilon \\ 
	&-   \sup_{\|\theta - \theta^*\|\leq u } \left| P_\theta\left(\langle Z_n , A(\theta^*) e' \rangle \leq c' - \frac{y^{\intercal} H^* y}{ 4 }   - \epsilon \right) - \Phi \left( \frac{ c'-\epsilon}{ \|A(\theta^*) e'\|  }  -\frac{y^{\intercal} H^* y }{ 4 \|A(\theta^*) e'\| }  \right)\right| \\
	& \geq \Phi \left( \frac{ c'-M\epsilon}{ \|A(\theta^*) e'\|  }  -\frac{y^{\intercal} H^* y }{ 4 \|A(\theta^*) e'\| }  \right)  + C\epsilon 
	\end{split}
	\end{equation*}
	Moreover for all $t \in \mathbb{ R}$, writing $x(\theta)$ to emphasize its dependence in $\theta$, 
	\begin{equation*}
	\begin{split}
	\Delta_1 &= \int_{B_n( \frac{M}{n^{1/4}}) }  \sup_{\|\theta' - \theta^*\|\leq u_0 } \left| P_{\theta'} \left(\langle Z_n , A(\theta^*) e' \rangle \leq t - \frac{x(\theta) ^{\intercal} H^* x(\theta) }{ 4 }    \right) - \Phi \left( \frac{t }{ \|A(\theta^*) e'\|  } -\frac{x(\theta)^{\intercal} H^* x (\theta)}{ 4 \|A(\theta^*) e'\| }  \right)  \right| d\theta \\
	& =  n^{-k_\theta/4} \int_{\|x\| \leq M } \sup_{\| \theta' - \theta^*\| \leq u_0   } \left| P_{\theta'} \left(\langle Z_n , A(\theta^*) e' \rangle \leq t - \frac{x^{\intercal} H^* x }{ 4 }    \right) - \Phi \left( \frac{t }{ \|A(\theta^*) e'\|  } -\frac{x^{\intercal} H^* x}{ 4 \|A(\theta^*) e'\| }  \right)  \right| dx\\
	& =   n^{-k_\theta/4} M o(1) 
	\end{split}
	\end{equation*}
	where the last equality follows from \textbf{[A4]}
	We then have  
	\begin{equation}\label{ub1}
	\begin{split}
	\int_{B_n( M /n^{1/4})} &P_\theta   \left( \|\eta(\mathbf z ) - \eta(\y)\|^2\leq \epsilon_n^2  \right)\pi(\theta) d\theta\\
	& \leq \pi(\theta^*)( 1 +o(1))n^{-k_\theta /4} \int_{\|x\| \leq M} \Phi \left( \frac{ c'+M\epsilon}{ \|A(\theta^*) e'\|  }  -\frac{x^{\intercal} H^* x }{ 4 \|A(\theta^*) e'\| }  \right) dx  + o(n^{-k_\theta /4}M)
	\end{split}
	\end{equation}
	Similarly 
	\begin{equation}\label{lb}
	\begin{split}
	& \int_{\Theta} P_\theta   \left( \| \eta(\mathbf z ) - \eta(\y)\|^2\leq \epsilon_n^2  \right)\pi(\theta) d\theta \\
	& \geq 
	\int_{\| \theta - \theta^*\| \leq n^{-1/4} M} P_\theta   \left( \|\eta(\mathbf z ) - \eta(\y)\|^2\leq \epsilon_n^2  \right)\pi(\theta) d\theta\\
	&\geq  \pi(\theta^*)(1+o(1) )n^{-k_\theta /4} \int_{\|x\| \leq M} \Phi \left( \frac{ c'- M\epsilon}{ \|A(\theta^*) e'\|  }  -\frac{x^{\intercal} H^* x }{ 4 \|A(\theta^*) e'\| }  \right) dx(1 + o_{P_{0}}(1) ).
	\end{split}
	\end{equation}
	
	Also if $\epsilon > \| \theta  - \theta^*\| > M /n^{1/4}$, since there exists $a>0$ such that  $z^{\intercal} H^* z \geq  a \|z\|^2$ for all $z$, if $M$ is large enough,
	\begin{equation*}
	\begin{split}
	P_\theta  \left( \|\eta(\mathbf z ) - \eta(\y)\|^2\leq \epsilon_n^2  \right)  \leq P_\theta\left(\langle Z_n , A(\theta^*) e' \rangle\leq   - \frac{ a M^2  }{ 8 }    \right) 
	\leq P_\theta \left( \|Z_n\|  >  \frac{ a M^2  }{ 8 \|A(\theta^*)e'\| }\right) \lesssim M^{-2\kappa } 
	\end{split}
	\end{equation*}
	Now let $j \geq 0$ and set $M_j = 2^j M$. On $M_j n^{-1/4}\leq \|\theta -\theta^*\| \leq M_{j+1} n^{-1/4} $ 
	\begin{equation*}
	\begin{split}
	P_\theta  \left( \|\eta(\mathbf z ) - \eta(\y)\|^2\leq \epsilon_n^2  \right)  \leq P_\theta\left(\langle Z_n , A(\theta^*) e' \rangle\leq   - \frac{ a M_j^2  }{ 8 }    \right) 
	\leq P_\theta \left( \|Z_n\|  >  \frac{ a M_j^2  }{ 8 \|A(\theta^*)e\| }\right) \lesssim M_j^{-2\kappa }
	\end{split}
	\end{equation*}
	so that 
	\begin{equation}\label{ub2}
	\begin{split}
	\int_{Mn^{-1/4} \leq \|\theta - \theta^*\|\leq \epsilon } P_\theta  \left( \|\eta(\mathbf z ) - \eta(\y)\|^2\leq \epsilon_n^2  \right)\pi(\theta) d\theta & \lesssim n^{-k_\theta/4} \sum_{j=0}^{J_n}  M_j^{-2\kappa } (M_{j+1}- M_j) \\
	& \lesssim n^{-k_\theta/4} \sum_{j=0}^{J_n} M_j^{-2\kappa +1} \lesssim n^{-k_\theta/4} M^{-2\kappa+ 1}. 
	\end{split}
	\end{equation}
	Finally if $\| \theta - \theta^*\| > \epsilon $, $\|b(\theta) - b_0 - Z_y/\sqrt{n} \|^2 - (\epsilon^*)^2  \geq C\epsilon $ on $\Omega_n$ and when $n$ is large enough $\sqrt{n}( \epsilon_n^2 - (\epsilon^*)^2 ) \leq c + C \epsilon  $ so that 
	\begin{equation*}
	\begin{split}
	P_\theta  \left( \|\eta(\mathbf z ) - \eta(\y)\|^2\leq \epsilon_n^2  \right) \\
	&  \leq  P_\theta\left(\langle\tilde Z_n ,b(\theta) - b_0\rangle \leq   -  \sqrt{n} C\epsilon  + c'+ C\epsilon   \right)+ O( n^{-\kappa/4} )\\
	& \leq  P_\theta ( \|Z_n\| > C n^{1/2} \|b(\theta) -b_0 \|^{-1}/2) + O( n^{-\kappa/4} ) \\
	& \leq  c(\theta) C^{-\kappa}n^{-\kappa/2}\| b(\theta) - b_0\|^{\kappa} + O( n^{-\kappa/4} ) .
	\end{split}
	\end{equation*}
	Therefore
	\begin{equation}\label{ub3}
	\begin{split}
	\int_{\|\theta - \theta^*\|\geq \epsilon} P_\theta  \left( \|\eta(\mathbf z ) - \eta(\y)\|^2\leq \epsilon_n^2  \right)\pi(\theta) d\theta  
	& \leq  C^{-\kappa}n^{-\kappa/2} \int_{\|\theta - \theta^*\|\geq \epsilon} c(\theta) \| b(\theta) - b_0\|^{\kappa} \pi(\theta) d\theta 
	\\
	& \quad  + C n^{-\kappa/4}\int_{\|\theta - \theta^*\|\geq \epsilon} c(\theta)  \pi(\theta) d\theta 
	\end{split}
	\end{equation}
	Finally combining \eqref{ub1}, \eqref{ub2}, \eqref{ub3} and \eqref{lb} we obtain that if $\kappa > k_\theta$
	\begin{equation*}
	\int_\Theta P_\theta  \left( \|\eta(\mathbf z ) - \eta(\y)\|^2\leq \epsilon_n^2  \right)\pi(\theta) d\theta  =  n^{-k_\theta /4} \int_{\mathbb R} \Phi \left( \frac{ c'}{ \|A(\theta^*) e'\|  }  -\frac{x^{\intercal} H^* x }{ 4 \|A(\theta^*) e'\| }  \right) dx + o(n^{-k_\theta /4})
	\end{equation*}
	and for all  $x = n^{1/4} (\theta- \theta^*)  \in \mathbb R^{k_\theta}$ fixed, writing $\pi_{z_n, \epsilon}(\cdot)$ the density of $\Pi_{z_n, \epsilon}$, the ABC posterior distribution of $z_n (\theta -\theta^*)$,
	\begin{equation*}
	\begin{split}
	\pi_{n^{1/4},\epsilon}(x)  &= \frac{  \Phi \left( \frac{ c'}{ \|A(\theta^*) e\|  }  -\frac{x^{\intercal} H^* x }{ 4 \|A(\theta^*) e'\| }  \right) }{  \int_{\mathbb R} \Phi \left( \frac{ c'}{ \|A(\theta^*) e'\|  }  -\frac{x^{\intercal} H^* x }{ 4 \|A(\theta^*) e'\| }  \right) dx} + o(1)   := q_c(x) + o(1) 
	\end{split}
	\end{equation*}
	so that $\left\|  \pi_{n^{1/4},\epsilon} - q_c \right\|_1 = o(1) .$

	We now study the case where $\sqrt{n}(\epsilon_n^2 - (\epsilon^*)^2 ) := \sqrt{n} u_n^2 \rightarrow + \infty$ with $u_n = o(1)$ and we show that the limiting distribution is uniform. 
	Using \eqref{deve:likeli}, we have that if $B_{0,n} = \{(\theta -\theta^*)^{\intercal} H^*(\theta- \theta^*)  \leq 2 u_n^2 - 4M_n/\sqrt{n} \} $, with $M_n < u_n^2 \sqrt{n} \epsilon $ going to infinity 
	\begin{equation*}
	\begin{split}
	P_\theta ( \| \eta(\mathbf z) -\eta(\y) \|\leq \epsilon_n ) & \leq 1 \\
	& \geq P_\theta\left( \langle \tilde Z_n, e'\rangle\leq  2M_n + \langle\tilde Z_y, e'\rangle- \epsilon - \|\theta - \theta^*\|M_n\right) \\
	& \geq P_\theta\left( \langle \tilde Z_n, e'\rangle\leq  M_n/2\right) \geq 1 - \frac{ c_1 }{ M_n^{-\kappa} }   
	\end{split}
	\end{equation*}
	for some $c_1>0$ on the event $\{|\langle\tilde Z_y, e'\rangle| \leq M_n/2\}$, which has probability going to $1 $. 
	
	This implies in particular that
	\begin{equation} \label{unif:equiv}
	\begin{split}
	\int_{B_{0,n} } P_\theta ( \| \eta(\mathbf z) -\eta(\y) \|\leq \epsilon_n ) \pi(\theta) d\theta & \leq \pi(\theta^*)( 1 + o(1)) \mbox{Vol}( B_{0,n} )   \\
	& \geq \pi(\theta^*)( 1 + o(1)) \mbox{Vol}( B_{0,n} )  (1 - \frac{ c_1 }{ M_n^{-\kappa} }   )  
	\end{split}
	\end{equation}
	
	Also 
	\begin{equation}\label{volB0}
	\mbox{Vol}(B_{0,n} ) \asymp u_n^{k_\theta}
	\end{equation}

	Let $K_n \geq  4$, if  $K_n u_n^2  \geq (\theta -\theta^*)^{\intercal} H^*(\theta- \theta^*)  > 2 u_n^2( 1-C\epsilon)^{-1} + 4M_n/\sqrt{n}$, then there exists $C'>0$ such that 
	\begin{equation*}
	\begin{split}
	& \frac{\sqrt{n}[u_n^2 - (1-C\epsilon)(\theta - \theta^*)^{\intercal} H^*(\theta- \theta^*)/2  ]}{2} + \langle\tilde Z_y ,e' \rangle  - \epsilon - \| \theta - \theta^*\| M_n\\
	& \leq - 2 M_n+ \langle\tilde Z_y ,e' \rangle  - \epsilon - C' \epsilon^2 M_n \leq - M_n 
	\end{split}
	\end{equation*}
	on the event $\{|\langle\tilde Z_y, e'>| \leq M_n/2\}$. Therefore writing $B_{1,n} = \{ K_n u_n^2  \geq (\theta -\theta^*)^{\intercal} H^*(\theta- \theta^*)  > 2 u_n^2( 1-C\epsilon)^{-1} + 4M_n/\sqrt{n}\}$
	\begin{equation}\label{unifub1}
	\begin{split}
	\int_{B_{1,n} } P_\theta ( \| \eta(\mathbf z) -\eta(\y) \|\leq \epsilon_n ) \pi(\theta) d\theta & \lesssim M_n^{-\kappa}  \mbox{Vol}( B_{1,n} )   \lesssim M_n^{-\kappa} K_n^{k_\theta/2}  \mbox{Vol}( B_{0,n} )
	\end{split}
	\end{equation}
	
	Moreover 
	\begin{equation}\label{DeltaVol1}
	\mbox{Vol}\left(\{ (\theta -\theta^*)^{\intercal} H^*(\theta- \theta^*)  \leq 2 u_n^2( 1-C\epsilon)^{-1} + 4M_n/\sqrt{n}\}\right) - \mbox{Vol}(B_{0,n}) \lesssim \epsilon \mbox{Vol}(B_{0,n}).
	\end{equation}
	
	If $K_n u_n^2  \leq (\theta -\theta^*)^{\intercal} H^*(\theta- \theta^*) \leq \epsilon^2 $, then 
	\begin{equation*}
	\begin{split}
	& \frac{\sqrt{n}[u_n^2 - (1-C\epsilon)(\theta - \theta^*)^{\intercal} H^*(\theta- \theta^*)/2  ]}{2} + \langle\tilde Z_y ,e' \rangle  - \epsilon - \| \theta - \theta^*\| M_n\\
	& \leq - \frac{ \sqrt{n} (\theta - \theta^*)^{\intercal} H^*(\theta- \theta^*) }{ 8 } 
	\end{split}
	\end{equation*}
	when $n$ is large enough and there exists $b>0$ such that 
	\begin{equation} \label{unifub2}
	\begin{split}
	\int_{B_{1,n}^c }\1_{(\theta -\theta^*)^{\intercal} H^*(\theta- \theta^*) \leq \epsilon^2 } P_\theta ( \| \eta(\mathbf z) -\eta(\y) \|\leq \epsilon_n ) \pi(\theta) d\theta & \lesssim 
	n^{-\kappa /2} \int_{B_n(A \epsilon^2) }\1_{ \| \theta - \theta^*\| \geq b\sqrt{K_n} u_n } \|\theta - \theta^*\|^{ -2\kappa} d\theta \\
	& \lesssim n^{-\kappa /2}\int_{b \sqrt{K_n}u_n}^{\sqrt{A}\epsilon} r^{k_\theta - 2\kappa -1 } dr 
	\end{split}
	\end{equation} 
	Since $k_\theta < 2 \kappa $ then the above term is of order 
	$$K_n^{(k_\theta - 2\kappa)/2} (\sqrt{n}u_n^2)^{ -\kappa/2} u_n^{k_\theta} \asymp  K_n^{(k_\theta - 2\kappa)/2} (\sqrt{n}u_n^2)^{ -\kappa/2} \mbox{Vol} (B_{0,n}) = o( \mbox{Vol} (B_{0,n}))$$
	
	Finally if $\|\theta - \theta^*\|\geq \epsilon$, similarly to the case where $\sqrt{n} u_n^2 \rightarrow c \in \mathbb R$, we obtain \eqref{ub3}
	and this term is $o(\mbox{Vol} (B_{0,n}))$ as soon as $n^{-\kappa/4} = o(u_n^{k_\theta})$. Since $n^{-1/4} = o(u_n )$ the latter is true as soon as $\kappa \geq k_\theta$. 
	Combining \eqref{unif:equiv}, \eqref{unifub1}, \eqref{unifub2} , \eqref{DeltaVol1} and \eqref{ub3}, we obtain that 
	\begin{equation}\label{Denom:unif}
	\left \| \int_\Theta P_\theta ( \| \eta(Z) -\eta(\y) \|\leq \epsilon_n ) \pi(\theta) d\theta - \pi(\theta^*) \mbox{Vol}(\tilde B_{0,n} ) \right\| = o_p(1) 
	\end{equation}
	where $\tilde B_{0,n} = \{(\theta -\theta^*)^{\intercal} H^*(\theta- \theta^*)  \leq 2 u_n^2\}$. 
	Let $x = u_n^{-1} (\theta - \theta^*)$ be fixed and $x^{\intercal} H^* x < 2$, then for $n$ large enough $x^{\intercal}H^* x \leq 2 - 4M_n^2/(\sqrt{n}u_n^2 )$ and using 
	$$\pi_{u_n^{-1}, \epsilon}(x) = \pi_\epsilon( \theta^* + u_nx | \y) u_n^{k_\theta} $$
	then $$\pi_{u_n^{-1}, \epsilon}(x)  = 1 + o_p(1) .$$
	If $x^{\intercal} H^* x > 2$, then if $\epsilon >0$ is small enough and $n$ is large enough $x^{\intercal}H^* x \geq 2(1-C\epsilon)^{-1} + 4M_n^2/(\sqrt{n}u_n^2 )$ and 
	$$\pi_{u_n^{-1}, \epsilon}(x)=   o_p(1) .$$
	This implies that the ABC posterior distribution of $u_n^{-1}( \theta - \theta^*) $ converges to the Uniform distribution over the ellipsoid $\{ x^{\intercal} H^* x \leq 2\}$ in total variation. 
\end{proof}

\subsection*{Proof of Proposition \ref{prop:counter}}
\begin{proof}
	To prove  Proposition \ref{prop:counter}, we prove that the approximate likelihood 
	$$P_\theta \left( \left\| \tilde Z_n - Z_y +\sqrt{n} (b(\theta) - b_0)\right\|^2\leq \epsilon_n^2 \right) $$
	is highly peaked around $\theta \neq 0$, and, as such, concentration around $\theta^*=0$ can not result. 
	
	As in the proof of Theorem \ref{normal_thm}, writing $Z_n =  Z= (Z_1, Z_2)^{\intercal} $,  we can define $W  = Z/\|Z\|$ and $R = \|Z\|/ v_\theta$ and  we have that  $W$ and $R$ are independent  and that their distribution does not depend on $\theta$. In particular $R^2 \sim \chi^2(2)$. 
	Now, set $h = b(\theta) - b_0 - Z_y/\sqrt{n} $, so that 
	\begin{flalign}
	\left\| Z - Z_y +\sqrt{n} (b(\theta) - b_0)\right\|^2  -n\epsilon_n^2 =& v_\theta^2 R^2  +2 \sqrt{n} Rv_\theta \langle W, h\rangle + n(\|h\|^2 - \epsilon_n^2)\leq 0\label{R_eq}
	\end{flalign}
	if and only if  
	\begin{equation*}
	\Delta (W) = v_\theta^2 n (\langle W, h\rangle^2 - \|h\|^2 + \epsilon_n^2)=  v_\theta^2 n \tilde \Delta (W) \geq 0,\quad R \in (r_1(W), r_2(W)) \cap \mathbb R_+, 
	\end{equation*}
	where 
	$$r_1(W) =  \frac{\sqrt{n}}{ v_\theta} \left[- \langle W, h\rangle  - \sqrt{\tilde \Delta}(W)\right], \quad  r_2(W) = \frac{\sqrt{n}}{ v_\theta} \left[- \langle W, h\rangle  + \sqrt{\tilde \Delta}(W)\right].$$
	
	Note that $\tilde \Delta(W) \leq \epsilon_n^2 $ so that if  $\tilde \Delta (W) \geq 0$ then $|\langle W,h\rangle | = \|h\|(1+ O(\epsilon_n ))$, given that $\|h\|  \asymp  1$ on the event $\| Z_y\| \leq M$ for some arbitrarily large $M$. Therefore if $ \langle- W, h\rangle  \leq 0 $, then $\langle -W, h \rangle \asymp - \|h\|$ and there is no solution for $R$ in \eqref{R_eq}. Hence \eqref{R_eq} holds if and only if  $\langle-W,h\rangle \geq 0$, $\tilde \Delta (W) \geq 0 $ and $R  \in (r_1(W), r_2(W))$.
	
	By symmetry we can set $W = -W$ and,  on the set $\langle W,h \rangle \geq 0$, using the fact that $R^2 \sim \chi^2(2)$, 
	\begin{flalign*}
	P_\theta \left(  \left\| Z - Z_y +\sqrt{n} (b(\theta) - b_0)\right\|^2  \leq n\epsilon_n^2 | W \right) &= e^{- r_1(W)^2/2 } ( 1  -  e^{-\frac{ n \tilde \Delta(W)}{v_\theta^2}} )   \\ r_1(W) &\in   \left( \sqrt{n} \|h\| /v_\theta\left\{1 - 2\epsilon_n\right\}, \sqrt{n} \|h\| /v_\theta\right)
	\end{flalign*}
	To derive an approximation of $P_\theta \left(  \left\| Z - Z_y +\sqrt{n} (b(\theta) - b_0)\right\|^2  \leq n\epsilon_n^2 \right)$ we study more precisely $r_1(W)$. For the sake of simplicity we assume that $\sqrt{n} \epsilon_n = o(1)$, since the case where $\sqrt{n} \epsilon_n = O(1)$ can be treated similarly. 
	Then 
	\begin{equation*}
	\begin{split}
	P_\theta \left(  \left\| Z - Z_y +\sqrt{n} (b(\theta) - b_0)\right\|^2  \leq n\epsilon_n^2  \right) &\leq  e^{- \frac{n\|h\|^2 }{2v_\theta^2} ( 1 - 2\epsilon_n)^2} P_\theta\left( \tilde \Delta (W) \geq 0 \right)\\
	& \geq e^{- \frac{n\|h\|^2 }{2v_\theta^2}} P_\theta\left( \tilde \Delta (W) \geq 0 \right)
	\end{split}
	\end{equation*}
	
	Consider $\theta = r\bar b_0$ with $r \in [-1,1]$ so that  $r= 0$ corresponds to $\theta = \theta*\equiv 0$, then $h = \bar b_0 ( r-1, r+1)^\intercal + O_{P_0}(1/\sqrt{n})$ and 
	\begin{equation*}
	\begin{split}
	P_\theta\left( \tilde \Delta (W) \geq 0 \right) & \leq 2 P_\theta \left(0\leq  \langle W-h/\|h\|, h \rangle \leq \epsilon_n^2/(2\|h\|^2 ) ( 1 + \epsilon_n^2/\|h\|^2)\right)\\
	& \geq 2 P_\theta \left(0\leq  \langle W-h/\|h\|, h \rangle \leq \epsilon_n^2/(2\|h\|^2 ) \right)\\
	& =   \frac{ \epsilon_n^2 }{  (1+r^2)\bar b_0^2 }  g\left\{\frac{(1-r)}{ \sqrt{2 (1+r^2)}} \right\}\left\{ 1 + O(\epsilon_n^2\vee 1/\sqrt{n} ) \right\}
	\end{split}
	\end{equation*}
	where $g(\cdot)$ is the density of $W_1$, with $W = (W_1, W_2)$. 
	We thus obtain that for $n$ large enough
	\begin{equation*}
	\begin{split}
	P_\theta \left(  \left\| Z - Z_y +\sqrt{n} (b(\theta) - b_0)\right\|^2  \leq n\epsilon_n^2  \right) & \leq    \frac{ \epsilon_n^2 }{  (1+r^2)\bar b_0^2 } e^{- \frac{n\bar b_0^2 (1+r^2) }{2v_\theta^2} ( 1 - 3\epsilon_n)^2 }  \times  \\&\;\; g\left\{\frac{\bar b_0(1-r)}{ \sqrt{2 (1+r^2)}} \right\}\left\{ 1 + O(\epsilon_n^2\vee 1/\sqrt{n} ) \right\}\\
	& \geq   \frac{ \epsilon_n^2 }{  (1+r^2)\bar b_0^2 } e^{- \frac{n\bar b_0^2 (1+r^2) }{2v_\theta^2} ( 1 + \epsilon_n)^2 } \times   \\& \;\; g\left\{\frac{\bar b_0(1-r)}{ \sqrt{2 (1+r^2)}} \right\}\left\{ 1 + O(\epsilon_n^2\vee 1/\sqrt{n} ) \right\}\\
	\end{split}
	\end{equation*}
	Take $v_{r \bar b_0}  = \bar b_0 v(r) $ such that $(1+1/4)/v(1/2)^2 \leq 1/(2v(0)^2) $, then for $\delta>0$ small enough,
	\begin{equation*}
	\Pi_\epsilon\left(| \theta - \theta^* | \leq  \delta | \eta(\mathbf y) \right)   = o\left( \Pi_\epsilon\left\{ | \theta - \bar b_0/2| \leq  \delta | \eta(\mathbf y) \right\}  \right)
	\end{equation*}
	since there exists $c>0$ such that 
	\begin{equation*}
	\int_{| \theta | \leq \delta }e^{- \frac{n(\bar b_0^2 +\theta^2) }{2v_\theta^2} }\pi(\theta) d\theta \leq e^{- n  c } \int_{| \theta -\bar b_0/2| \leq \delta }e^{- \frac{n(\bar b_0^2 +\theta^2) }{2v_\theta^2}} \pi(\theta) d\theta.
	\end{equation*}
	
\end{proof}

\subsection*{Proof of Corollary \ref{th:asympAdj}}
\begin{proof}
	The proof is a consequence of Theorem \ref{th:asymp0} and the structure of $\tilde{\theta}=\theta-\hat{\beta}^{\intercal}\{\eta(\mathbf{z})-\eta(\mathbf{y})\}$, and $\tilde{\theta}^{*}=\theta^{*}-\beta_0^{\intercal}\{b(\theta^*)-b_0\}$. Therefore, we only sketch the idea here. 
	
	Take $\delta_{n}\geq M_{n}(\epsilon_{n}-\epsilon^*)\geq M_{n}v_{n}^{-1}$. By assumption $\epsilon^*>0$ and $\|\beta_0\|>0$. Define $\Omega_{d}=\{\mathbf{y}:\|\eta(\mathbf{y})-b_0\|\leq \delta_{n}/u_0\}$ for some $u_0\geq2 (1+\|\beta_0\|)$. By the result of Theorem \ref{th:asymp0} we have that 
	\begin{flalign*}
	\widetilde{\Pi}_{\epsilon_{}}\left[|\tilde{\theta}-\tilde{\theta}^*|>\delta_{n}|\y\right] &={\Pi}_{\epsilon_{}}\left[\{\theta:|\tilde{\theta}-\tilde{\theta}^*|>\delta_{n}\}\cap\{\theta:|{\theta}-{\theta}^*|\leq\delta_{n}/u_{0}\}|\y\right]+o_{P_{0}}(1)
	\\&=\frac{\int_{|\theta-\theta^*|\leq \delta_{n}/u_{0}}\1\left[|\tilde{\theta}-\tilde{\theta}^*|\geq \delta_{n}\right]P_{\theta}\left[\|\eta(\mathbf{z})-\eta(\mathbf{y})\|\leq\epsilon_{n}\right]d\Pi(\theta)   }{\int_{|\theta-\theta^*|\leq \delta_{n}/u_{0}}P_{\theta}\left[\|\eta(\mathbf{z})-\eta(\mathbf{y})\|\leq\epsilon_{n}\right]d\Pi(\theta)}+o_{P_{0}}(1),\end{flalign*} where both equalities follow by posterior concentration of $|\theta-\theta^*|$ at rate $\delta_{n}\gg v_{0,n}^{-1}$. Similar steps to that of Theorem \ref{th:asymp0} yield $$D_{n}=\int_{|\theta-\theta^*|\leq \delta_{n}/u_{0}}P_{\theta}\left[\|\eta(\mathbf{z})-\eta(\mathbf{y})\|\leq\epsilon_{n}\right]d\Pi(\theta)\gtrsim\delta_{n}^{D},$$ under case (i) or case (ii) of \textbf{[A1]}. Define the event $$S(\delta_{n})=\left\{(\mathbf{z},\theta):\{\theta:|\tilde{\theta}-\tilde{\theta}^*|>\delta_{n}\}\cap\{\theta:|{\theta}-{\theta}^*|\leq\delta_{n}/u_{0}\}\cap\{\mathbf{z}:\|\eta(\mathbf{z})-\eta(\mathbf{y})\|\leq\epsilon_{n}\}\right\}$$ Note that
	\begin{flalign*}
	\tilde{\theta}-\tilde{\theta}^*=&\theta-\theta^*+[\hat{\beta}-\beta_0]^{\intercal}\{b(\theta)-b_{0}\}+[\hat{\beta}-\beta_0]^{\intercal}\{b_{0}-\eta(\mathbf{y})\}+[\hat{\beta}-\beta_0]^{\intercal}\{\eta(\mathbf{z})-b(\theta)\}\\&+[\hat{\beta}-\beta_0]^{\intercal}\{b(\theta)-b(\theta^*)\}+\beta_{0}^{\intercal}\{b_{0}-\eta(\mathbf{y})\}+\beta_{0}^{\intercal}\{\eta(\mathbf{z})-b(\theta)\}+\beta_{0}^{\intercal}\{b(\theta)-b(\theta^*)\}
	\end{flalign*}
	For $\mathbf{y}\in\Omega_{d}$, we have
	\begin{flalign*}
	\delta_{n}<|\tilde{\theta}-\tilde{\theta}^*|\leq& |\theta-\theta^*|+\|\hat{\beta}-\beta_0\|\|b(\theta)-b_{0}\|+\|\hat{\beta}-\beta_0\|\|b_{0}-\eta(\mathbf{y})\|+\|\hat{\beta}-\beta_0\|\|\eta(\mathbf{z})-b(\theta)\|\\&+\|\hat{\beta}-\beta_0\|\|b(\theta)-b(\theta^*)\|+\|\beta_{0}\|\|b_{0}-\eta(\mathbf{y})\|+\|\beta_{0}\|\|b(\theta)-b(\theta^*)\|+\|\beta_{0}\|\|\eta(\mathbf{z})-b(\theta)\|\\\leq&\delta_{n}/u_{0}+\|\beta_0\|\delta_{n}/u_{0}+o(\delta_{n})+\left(O(\delta_{n})+\|\beta_{0}\|\right)\|\eta(\mathbf{z})-b(\theta)\|+o_{P_{\theta}}(1)
	\end{flalign*}where the last inequality follows from $\|\hat{\beta}-\beta_0\|=o_{P_{\theta}}(1)$ and concentration of $|\theta-\theta^*|$ at rate $\delta_{n}\gg v_{n}^{-1}$. Therefore, take $u_0\geq2(1+\|\beta_0\|)$ and rearrange the above to obtain $$0<\frac{\delta_{n}}{2(O(\delta_{n})+\|\beta_0\|)}<\|\eta(\mathbf{z})-b(\theta)\|+o(\delta_{n}).$$This then implies that 
	\begin{flalign*}
	\text{Pr}\left[S(\delta_{n})\right]&=\int_{\{\theta:|{\theta}-{\theta}^*|\leq\delta_{n}/u_{0}\}}\1\left[|\tilde{\theta}-\tilde{\theta}^*|>\delta_{n}\right]P_{\theta}\left[\|\eta(\mathbf{z})-\eta(\mathbf{y})\|\leq\epsilon_{n}\right]d\Pi(\theta)\\&\leq\int_{\Theta}P_{\theta}\left[\|\eta(\mathbf{z})-b(\theta)\|>c\cdot\delta_{n}\right]d\Pi(\theta)\\&\lesssim(v_{n}\delta_{n})^{-\kappa}\text{ under case (i) of \textbf{[A1]}}\\ &\lesssim\exp(-cv^{\tau}_{n}\delta^{\tau}_{n})\text{ under case (ii) of \textbf{[A1]}}
	\end{flalign*}
	Recalling that $D_n\gtrsim\delta^{D}$ and using the above, the result follows similarly to Theorem \ref{th:asymp0}.
\end{proof}

\section{Additional Computations for Example 1}
In this section, we consider several additional aspects of the simulation exercise for Example 1 in the main text. For completeness, we now recall the general features of the example. The assumed DGP is $z_1,\dots,z_n$ iid as $\mathcal{N}(\theta,1)$ but the actual DGP is $y_1,\dots,y_n$ iid as $\mathcal{N}(\theta,{\sigma}^2)$. We conduct ABC using the following summary statistics:
\begin{itemize}
	\item the sample mean $\eta_{1}(\mathbf{y})=\frac{1}{n}\sum_{i=1}^{n}{y}_{i}$,
	\item the sample variance $\eta_{2}(\mathbf{y})=\frac{1}{n-1}\sum_{i=1}^{n}({y}_{i}-\eta_{1}(\mathbf{y}))^{2}$.
\end{itemize}
Our prior beliefs are given by $\theta\sim \mathcal{N}(0,25)$. For accept/reject ABC (ABC-AR), we use $N=25, 000$ simulated pseudo datasets generated iid according to $z^{j}_{i}\sim \mathcal{N}(\theta^j,1)$. For both ABC-AR and the local linear regression adjustment (ABC-Reg), we set the tolerance $\epsilon$ to be the 1\% quantile of the simulated distances $\|\eta(\mathbf{y})-\eta(\mathbf{z}^{j})\|$. 

\subsection{Example 1: ABC-AR and ABC-Reg Posterior Comparison}
In this section, we use the exact same data generated for the experiments conducted in Example 1 in Section 1 of the main text, however, we now analyze the posteriors for ABC-AR and ABC-Reg. For clarity, we briefly recall the details of the Monte Carlo design used in this example. A sequence of ``observed'' datasets
for $\y$ are generated, each corresponding to a different value of
${\sigma}^{2}$. To isolate the impact of model misspecification, which occurs if $\sigma^2\neq1$, across the experiments, the same set of random numbers are used to generate
the observed data. The sample size across the experiments is taken to be $n=100$. To further isolate the impact of misspecification, we implement ABC using the same pseudo-data across all experiments.

The example in Section 1 clearly demonstrated that the ABC-AR and ABC-Reg posterior means behave very differently depending on the level of model misspecification. Unsurprisingly, the ABC-AR and ABC-Reg posteriors themselves displaying very different patterns of behavior under model misspecification.

In Figure \ref{fig1} we plot the posterior densities obtained from ABC-AR (Panel A) and ABC-Reg (Panel B) across a subset of the different datasets used in the experiment that correspond to $\sigma^2\in\{1,1.05,1.10,\dots,5\}$.
\begin{figure}[H]
	\centering
	\includegraphics[width=19cm,height=8cm]{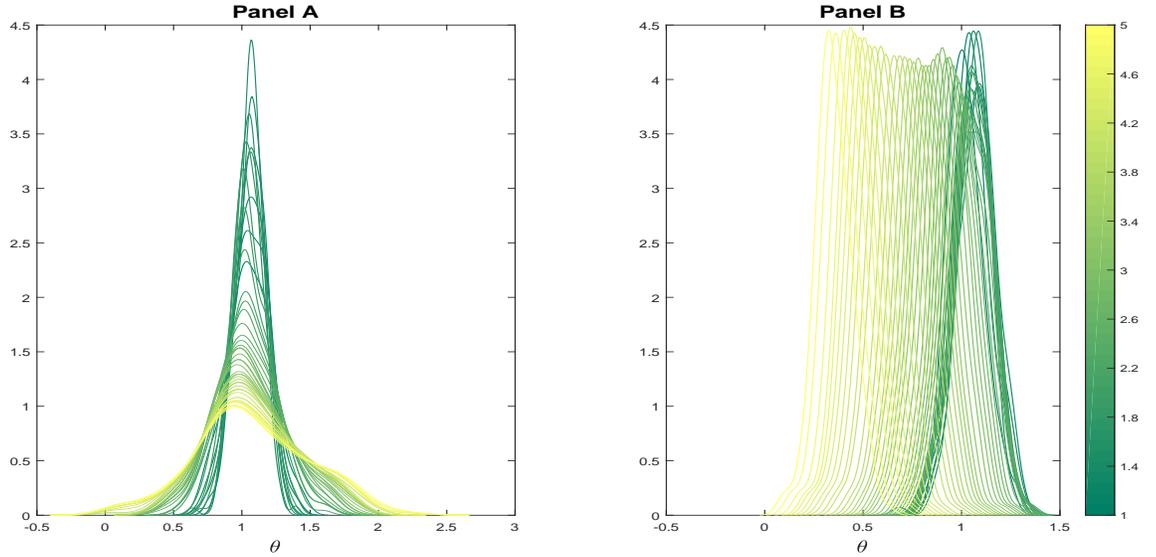} 
	\caption{Comparison of posterior densities for ABC-AR, and ABC-Reg across varying levels of model misspecification. The true value generating the data is $\theta=1$. The posteriors are colored so that darker colors represent less model misspecification ($\sigma^2$ closer to unity), and lighter colors representing larger degrees of model misspecification (larger values of $\sigma^2$).} 
	\label{fig1} 
\end{figure}

Analyzing Panel A in Figure \ref{fig1}, we see that the ABC-AR posteriors remain roughly centered near $\theta=1$, regardless of the value of $\sigma^2$. However, as $\sigma^2$ increases the ABC-AR posteriors noticeably flatten to accommodate the additional variability in the observed data. In contrast to ABC-AR, the mean of the ABC-Reg posterior (Panel B of Figure \ref{fig1}) shifts significantly as $\sigma^2$ increases and the variability of the ABC-Reg posterior remains roughly constant across all levels of misspecification. Therefore, ABC-Reg completely neglects the fact that the variance of the observed data is increasing as $\sigma^2$ increases.

This finding is further visually confirmed by the results in Figure \ref{fig2}, which plots the corresponding 95\% credible intervals (HPD intervals) for ABC-AR and ABC-Reg across the different values of $\sigma^2$ used in the experiments. The results in Figure \ref{fig2} demonstrate that as $\sigma^2$ increases the credible regions for ABC-AR widen to accommodate the increased variance in the observed data, whereas the credible intervals for ABC-Reg maintains roughly the same length as $\sigma^2$ increases. 
\begin{figure}[H]
	\centering
	\includegraphics[width=15cm, height=6cm]{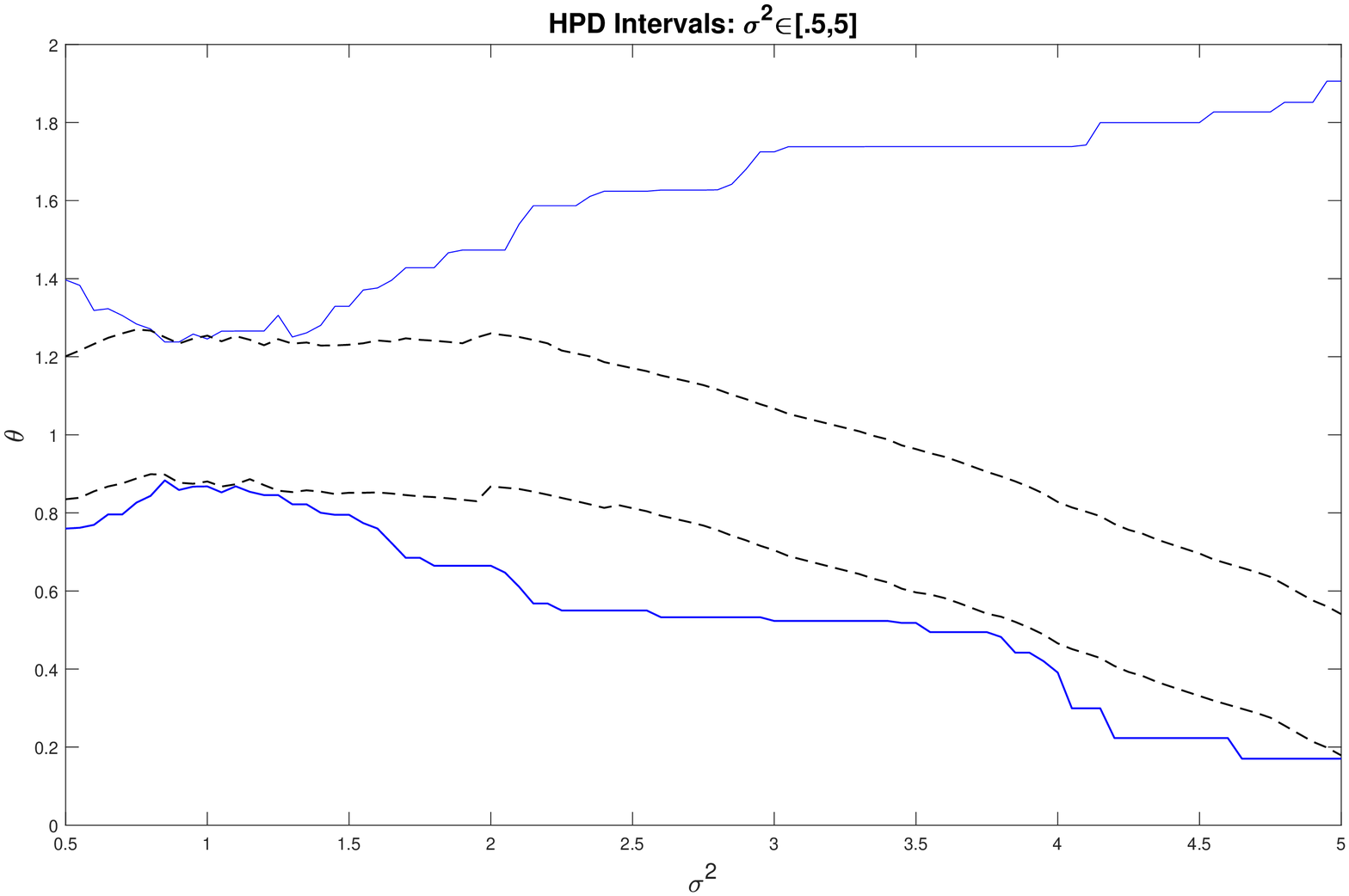} 
	\caption{Comparison of credible intervals (HPD intervals) for ABC-AR, and ABC-Reg across varying levels of model misspecification. The true value generating the data is $\theta=1$. The solid lines represent the HPD intervals for ABC-AR, while the dashed lines are the HPD intervals for ABC-Reg. } 
	\label{fig2} 
\end{figure}

\subsection{Example 1: Alternative Regression Adjustments}
In this section, we further explore the behavior of local regression adjustment approaches under misspecification. In particular, in the confines of Example 1, we will compare the repeated sampling behavior of ABC-AR and three different regression adjustment approaches: the local linear regression adjustment (ABC-Reg), the proposed local linear adjustment approach described in Section 3.2 (ABC-RegN), and the local nonlinear regression adjustment of \cite{blumF2010non} (ABC-NN). For each regression adjustment approach, we consider two versions: one version with the heteroskedasticity adjustment of \cite{blumF2010non} and one version without. All told, we will compare six different regression adjustment approaches across various levels of model misspecification. All adjustment methods are carried out using the \texttt{R} package \texttt{abc}  (\citealp{csillery2012abc}). 

We follow the Monte Carlo design considered in Section 3.2 of the main paper: for each Monte Carlo replication we simulate observed data $y_i\sim\mathcal{N}(1,\sigma^2)$, iid, and consider three different values of $\sigma^2$ corresponding to ${\sigma}^{2}\in \{1,2,3\}$.  For each value of $\sigma^2$ we generate 1,000 artificial observed data sets each of length $n=100$. Every ABC procedure relies on $N=25, 000$ pseudo-data sets generated according to $z_i\sim \mathcal{N}(\theta,1)$, iid, and for each method the tolerance is chosen to the 1\% quantile of the simulated distances $\|\eta(\y)-\eta(\z)\|$. As in the main text, the local regression adjustment procedures use the Epanechnikov kernel.

In Figure \ref{fig3}, we plot the posterior means of the different procedures, across the values of $\sigma^2$ and across the Monte Carlo replications, without the heteroskedasticity correction, and the results in Figure \ref{fig4} consider the case with the heteroskedasticity correction. While the results of Figure \ref{fig3} have already been presented in the main text, we reproduce them here to simplify the comparison. 

The results of this experiment demonstrate that the point estimators obtained from the ABC regression adjustment procedures can display larger variability, in this repeated sampling context, than those obtained from ABC-AR. Moreover, there seems to be no real benefit from the local nonlinear regression adjustment. In addition, the heteroskedasticity correction does not significantly alter the behavior of any local regression adjustment approach. Indeed, the two sets of results are visually very similar.

\begin{figure}[H]
	\centering
	\resizebox{\textwidth}{!}{\input{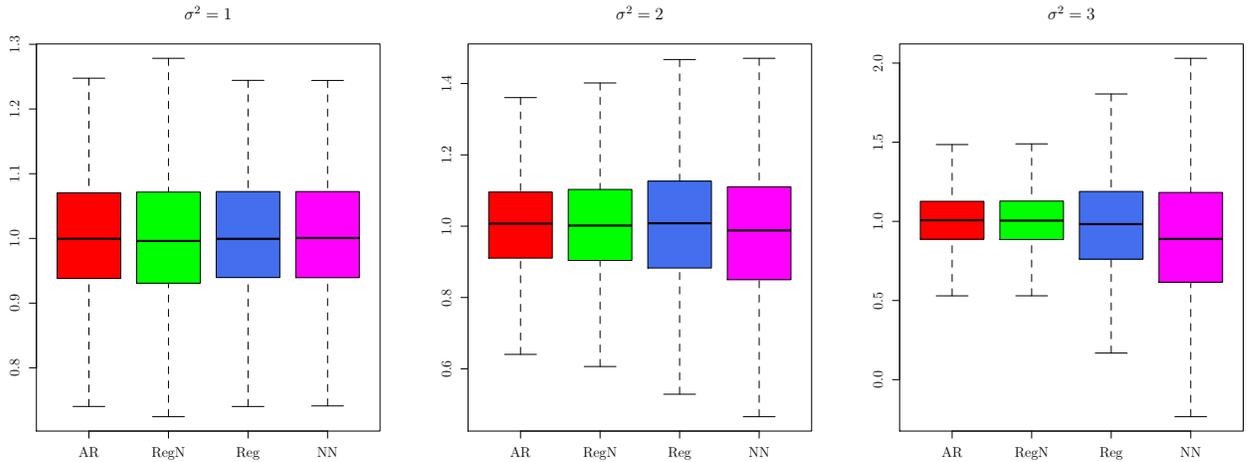}  }  
	\caption{{Posterior mean comparison of ABC-AR (AR), local linear regression adjustment (Reg), the proposed local linear regression adjustment (RegN), and the local nonlinear regression adjustment (NN) across three values of $\sigma^2$: $\sigma^2\in\{1,2,3\}$. Recall that $\sigma^2=1$ corresponds to correct model specification. All plots are presented without outliers.}} 
	\label{fig3}
\end{figure}
\begin{figure}[H]
	\centering
	\resizebox{\textwidth}{!}{\input{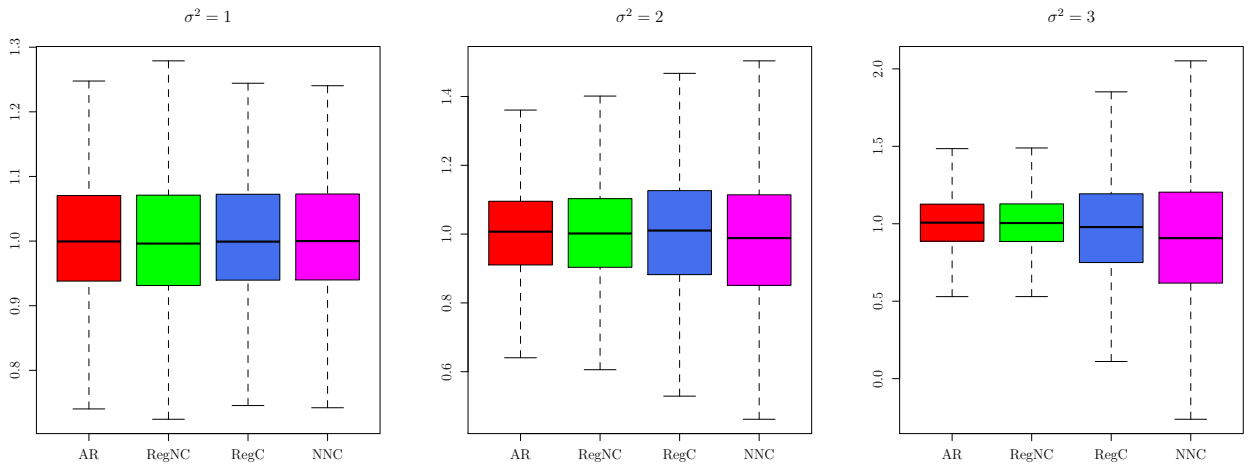}  }  
	\caption{{Posterior mean comparison of ABC-AR (AR), local linear regression adjustment with heteroskedasticity correction (RegC), the proposed local linear regression adjustment with heteroskedasticity correction (RegNC), and the local nonlinear regression adjustment with heteroskedasticity correction (NNC) across three separate values of $\sigma^2$: $\sigma^2\in\{1,2,3\}$. Recall that $\sigma^2=1$ corresponds to correct model specification. All results are presented without outliers.}} 
	\label{fig4}
\end{figure}

Across the Monte Carlo replications and for each ABC procedure, in Table \ref{tab1} we present the median posterior standard deviation, the Monte Carlo coverage, the median length of the 95\% posterior credible interval, as well as the 2.5\% and 97.5\% posterior quantiles (represented as the median value over the replications). {In this experiment, we compare medians instead of means due to the fact that the heteroskedasticity corrected regression adjustment approaches returned several large outliers that could otherwise skew the comparison.} 

The results in Table \ref{tab1} demonstrate that, in terms of posterior variability, as measured by posterior standard deviation, ABC-AR displays the largest variability across all designs. The local linear regression adjustments without heteroskedasticity correction (ABC-Reg, ABC-RegN) have posterior standard deviations that are virtually unchanged across the designs, while the standard deviation of the heteroskedasticity corrected local linear approach (ABC-RegC) is actually decreasing as $\sigma^2$ increases, i.e., as the variability in the data increase. In contrast, the posterior standard deviation for the nonlinear regression adjustment without heteroskedasticity correction (ABC-NN) increases as $\sigma^2$ increases, however, the  posterior standard deviation of the corrected version (ABC-NNC) is stable across the designs. 

As discussed in the main text, the fact that the local regression adjustment has small posterior variability, relative to ABC-AR, results in small credible sets and a false sense of precision. As a direct consequence, all the adjustment procedures, with and without the heteroskedasticity correction, have poor Monte Carlo coverage when $\sigma^2=2,3$.

Overall, these results suggest that, at least in this context, there is no meaningful difference between the results obtained by local linear or local nonlinear regression adjustments. In addition, the heteroskedasticity correction does not improve the behavior of the regression adjustment, and can potentially exacerbate the coverage issues observed in these procedure (see, e.g., the results of ABC-RegC in Table \ref{tab1}). Lastly, these results suggest that our proposed local linear regression adjustment (ABC-RegN) performs well relative to the other local regression adjustments (with and without heteroskedasticity correction).


\begin{table}[H]
	\centering
	\caption{{Monte Carlo coverage (Cov), credible set length (Len), and posterior standard deviation (Std) for the normal example under various levels of model misspecification. Cov is the percentage of times that the 95\% credible set contained $\theta=1$. Len is the median length of the 95\% credible set, across the Monte Carlo trials. Std is the median posterior standard deviation across the Monte Carlo trials. Q025 represents the median (over the replications) of the 2.5\% posterior quantiles across the different methods, while Q975 represents the corresponding 97.5\% quantile.}}
	\begin{tabular}{rlrrrrr}
		&       &       &       &       &       &  \\
		\cmidrule{2-7}          &  {$\sigma^2=1$}     & \multicolumn{1}{l}{Cov} & \multicolumn{1}{l}{Len} & \multicolumn{1}{l}{Std} & \multicolumn{1}{l}{Q025} & \multicolumn{1}{l}{Q975} \\
		\cmidrule{2-7}          & ABC-AR & 0.982 & 0.461 & 0.120 & 0.770 & 1.234 \\
		& ABC-RegN & 0.938 & 0.385 & 0.100 & 0.802 & 1.190 \\
		& ABC-Reg & 0.941 & 0.382 & 0.100 & 0.811 & 1.192 \\
		\multicolumn{1}{l}{} & ABC-NN & 0.950 & 0.385 & 0.100 & 0.811 & 1.194 \\
		& ABC-RegNC & 0.938 & 0.389 & 0.100 & 0.802 & 1.193 \\
		& ABC-RegC & 0.945 & 0.386 & 0.100 & 0.809 & 1.197 \\
		& ABC-NNC & 0.941 & 0.384 & 0.101 & 0.808 & 1.200 \\
		\cmidrule{2-7}          &       &       &       &       &       &  \\
		\cmidrule{2-7}          & {$\sigma^2=2$}      & \multicolumn{1}{l}{Cov} & \multicolumn{1}{l}{Len} & \multicolumn{1}{l}{Std} & \multicolumn{1}{l}{Q025} & \multicolumn{1}{l}{Q975} \\
		\cmidrule{2-7}          & ABC-AR & 0.961 & 0.614 & 0.158 & 0.698 & 1.310 \\
		& ABC-RegN & 0.802 & 0.383 & 0.100 & 0.809 & 1.194 \\
		& ABC-Reg & 0.717 & 0.382 & 0.100 & 0.818 & 1.200 \\
		\multicolumn{1}{l}{} & ABC-NN & 0.729 & 0.428 & 0.113 & 0.766 & 1.204 \\
		& ABC-RegNC & 0.809 & 0.387 & 0.100 & 0.805 & 1.194 \\
		& ABC-RegC & 0.645 & 0.361 & 0.094 & 0.791 & 1.223 \\
		& ABC-NNC & 0.636 & 0.383 & 0.100 & 0.773 & 1.201 \\
		\cmidrule{2-7}          &       &       &       &       &       &  \\
		\cmidrule{2-7}          &  {$\sigma^2=3$}     & \multicolumn{1}{l}{Cov} & \multicolumn{1}{l}{Len} & \multicolumn{1}{l}{Std} & \multicolumn{1}{l}{Q025} & \multicolumn{1}{l}{Q975} \\
		\cmidrule{2-7}          & ABC-AR & 0.913 & 0.613 & 0.157 & 0.694 & 1.311 \\
		& ABC-RegN & 0.707 & 0.383 & 0.099 & 0.816 & 1.195 \\
		& ABC-Reg & 0.460 & 0.381 & 0.100 & 0.791 & 1.180 \\
		\multicolumn{1}{l}{} & ABC-NN & 0.419 & 0.479 & 0.128 & 0.632 & 1.137 \\
		& ABC-RegNC & 0.708 & 0.386 & 0.100 & 0.811 & 1.199 \\
		& ABC-RegC & 0.462 & 0.261 & 0.069 & 0.700 & 1.249 \\
		& ABC-NNC & 0.409 & 0.383 & 0.101 & 0.624 & 1.171 \\
		\cmidrule{2-7}    \end{tabular}%
	\label{tab1}%
\end{table}%

\newpage

\section{Additional Example: Misspecified $g$-and-$k$ Model}
To further demonstrate the behavior of different ABC approaches under model misspecification, we consider an additional example based on the $g$-and-$k$ distribution, which is an oft-used example in the ABC literature to compare the behavior of different ABC approaches (see, e.g., \citealp{drovandi2011likelihood}, \citealp{FP2012}, and \citealp{Bernton2017}). The $g$-and-$k$ model is most commonly stated through it's quantile function: 
$$q \in(0,1) \mapsto a+b\left(1+0.8 \frac{1-\exp (-g z(q)}{1+\exp (-g z(q)}\right)\left(1+z(q)^{2}\right)^{k} z(q),$$where $z(q)$ refers to the $q$-th quantile of the standard normal distribution. The four parameters of the $g$-and-$k$ distributions have specific interpretations. The parameter $a$ represents the location, $b$ the scale, while $g$ and $k$ control the skewness and kurtosis, respectively. 
Following \citet{drovandi2011likelihood}, we consider the following priors on the parameters
\begin{flalign*}
a\sim \mathcal{U}[0,10],\;b\sim \mathcal{U}[0,10],\;\;g\sim \mathcal{U}[0,10]\;b\sim \mathcal{U}[0,10],
\end{flalign*}where $\mathcal{U}[0,1]$ denotes the uniform distribution on $[0,1]$.

ABC-based inference in the $g$-and-$k$ model is usually conducted using the quantiles of the simulated and observed data. Therefore, in what follows we take as our summary statistics for ABC the octiles of the data:
\begin{itemize}
	\item $\eta_{j}(\mathbf{y})=O_{j}(\y)$, for $1,\dots,7$, where $(O_{1},\dots,O_{7})$ partitions the data into eight equal parts.
\end{itemize}

The $g$-and-$k$ distribution is a highly-flexible class of distributions that is capable of modeling data with a complex unconditional distribution. While highly flexible, $g$-and-$k$ distributions are unimodal and are not capable of capturing multi-modality that may exist in the data. 

In this section, we compare the behavior of different ABC-based procedures when the underlying assumption is that the data comes from the $g$-and-$k$ distribution, but when the observed data is actually generated from a distribution with minor bi-modality. In particular, we generate observed data iid from the Gaussian mixture
\begin{flalign}\label{eq:mix}
y_i\sim w\cdot\mathcal{N}(\mu_1,\sigma_1^2)+(1-w)\cdot\mathcal{N}(\mu_2,\sigma_2^2).
\end{flalign}
In what follows, we fix the parameters in equation \eqref{eq:mix} to be $$(\mu_1,\sigma_1^2)^\intercal=(1,2)^\intercal,\;\;(\mu_2,\sigma_2^2)^\intercal=(7,2)^\intercal,\text{ and } w=0.9.$$This DGP produces observed data that exhibits positive skewness and excess kurtosis, owing to the fact that the density exhibits a minor ``hump'' in the right tail of the data. This specification was chosen particularly to generate minor bi-modality in the observed data, which is a feature that the $g$-and-$k$ model can not capture. For illustration, the kernel density of a representative data set of size $n=100$ simulated from this mixture model is given in Figure \ref{fig:mixdata}.

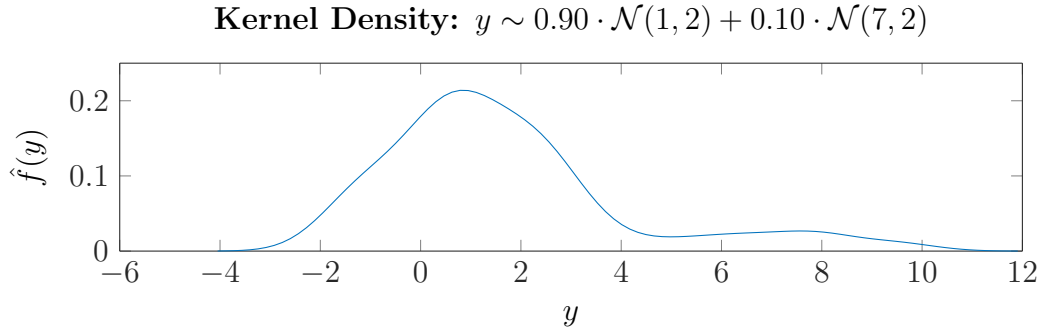
\begin{figure}[H]
	\centering 
	\setlength\figureheight{2.5cm} 
	\setlength\figurewidth{12cm} 
%
%
%
\definecolor{mycolor1}{rgb}{0.00000,0.44700,0.74100}%
\begin{tikzpicture}

\begin{axis}[%
width=\figurewidth,
height=\figureheight,
scale only axis,
separate axis lines,
every outer x axis line/.append style={white!15!black},
every x tick label/.append style={font=\color{white!15!black}},
xmin=-6,
xmax=12,
xlabel={$y$},
every outer y axis line/.append style={white!15!black},
every y tick label/.append style={font=\color{white!15!black}},
ymin=0,
ymax=0.25,
ylabel={$\hat{f}(y)$},
title style={font=\bfseries},
title={Kernel Density: $y\sim  0.90\cdot\mathcal{N}(1,2)+0.10\cdot\mathcal{N}(7,2)$},
legend style={draw=white!15!black,fill=white,legend cell align=left}
]
\addplot [color=mycolor1,solid,forget plot]
  table[row sep=crcr]{%
-4.04805456331554	0.00013115137101497\\
-3.8868806096807	0.000266011595106401\\
-3.72570665604587	0.000516621479648462\\
-3.56453270241103	0.000961105601075039\\
-3.40335874877619	0.00171350760993036\\
-3.24218479514136	0.00292897995357779\\
-3.08101084150652	0.0048025917927551\\
-2.91983688787169	0.00755786328389871\\
-2.75866293423685	0.0114223151917893\\
-2.59748898060201	0.0165903806555469\\
-2.43631502696718	0.0231787103166802\\
-2.27514107333234	0.0311839937943423\\
-2.11396711969751	0.0404568444092821\\
-1.95279316606267	0.0507047922971504\\
-1.79161921242783	0.0615316681450329\\
-1.630445258793	0.0725102785570525\\
-1.46927130515816	0.0832732460659572\\
-1.30809735152333	0.0935977097499023\\
-1.14692339788849	0.103457466051153\\
-0.985749444253655	0.113022975204309\\
-0.824575490618819	0.122603918434513\\
-0.663401536983983	0.132546107698401\\
-0.502227583349147	0.143108809302374\\
-0.341053629714311	0.154355361368578\\
-0.179879676079475	0.166087454568174\\
-0.018705722444639	0.177842654530711\\
0.142468231190198	0.188958816731191\\
0.303642184825033	0.198692036317897\\
0.464816138459869	0.20636072875671\\
0.625990092094705	0.211480785516604\\
0.787164045729542	0.213857757243396\\
0.948337999364377	0.213612088111294\\
1.10951195299921	0.21113053007048\\
1.27068590663405	0.206956620184401\\
1.43185986026889	0.20164982037206\\
1.59303381390372	0.195651416658111\\
1.75420776753856	0.189192753651845\\
1.91538172117339	0.182268539849705\\
2.07655567480823	0.174678864158424\\
2.23772962844307	0.166124252915904\\
2.3989035820779	0.156324468184557\\
2.56007753571274	0.145127542635807\\
2.72125148934757	0.132581215864142\\
2.88242544298241	0.118951791624311\\
3.04359939661725	0.104690882150045\\
3.20477335025208	0.0903637602902253\\
3.36594730388692	0.076560631364431\\
3.52712125752175	0.0638127692872\\
3.68829521115659	0.0525301214749137\\
3.84946916479143	0.0429683604813252\\
4.01064311842626	0.0352247482542354\\
4.1718170720611	0.0292562356919354\\
4.33299102569593	0.0249110079423486\\
4.49416497933077	0.0219655111508837\\
4.65533893296561	0.0201611335007027\\
4.81651288660044	0.0192365823742321\\
4.97768684023528	0.0189530495626364\\
5.13886079387011	0.0191100073789794\\
5.30003474750495	0.0195507785186909\\
5.46120870113979	0.0201591432566283\\
5.62238265477462	0.0208505169588739\\
5.78355660840946	0.0215624562624163\\
5.94473056204429	0.0222485788494335\\
6.10590451567913	0.0228775359808325\\
6.26707846931397	0.0234355572073765\\
6.4282524229488	0.0239287643751834\\
6.58942637658364	0.0243810214984264\\
6.75060033021847	0.0248247984605367\\
6.91177428385331	0.0252857490447531\\
7.07294823748814	0.0257651957468938\\
7.23412219112298	0.0262270247991924\\
7.39529614475782	0.0265954282359897\\
7.55647009839265	0.0267670357448872\\
7.71764405202749	0.0266359339846167\\
7.87881800566233	0.0261246705262371\\
8.03999195929716	0.0252109066521903\\
8.201165912932	0.0239397286033368\\
8.36233986656683	0.0224160622796746\\
8.52351382020167	0.0207786326905566\\
8.6846877738365	0.0191636782953065\\
8.84586172747134	0.0176703805932464\\
9.00703568110618	0.0163391912248971\\
9.16820963474101	0.015149358370965\\
9.32938358837585	0.0140350597495051\\
9.49055754201069	0.0129133883743405\\
9.65173149564552	0.0117142179784691\\
9.81290544928036	0.0104025361095994\\
9.97407940291519	0.00898748913030566\\
10.13525335655	0.0075174094219597\\
10.2964273101849	0.00606454135009241\\
10.4576012638197	0.00470563819699861\\
10.6187752174545	0.00350461770026606\\
10.7799491710894	0.00250153289656141\\
10.9411231247242	0.0017093061797328\\
11.102297078359	0.00111711943874857\\
11.2634710319939	0.000697819041980908\\
11.4246449856287	0.000416392834310377\\
11.5858189392636	0.000237233408109442\\
11.7469928928984	0.000128999099368363\\
11.9081668465332	6.69246624598708e-05\\
};
\end{axis}
\end{tikzpicture}%
	\caption{Kernel density of data simulated from the DGP in equation \eqref{eq:mix}.} 
	\label{fig:mixdata} 
\end{figure}

\subsection{Monte Carlo Experiments}
Similar to Example 1 in the main paper, we now compare the behavior of ABC-AR and various local ABC regression adjustments. In particular, we consider the following local regression adjustments: ABC-Reg (the standard weighted local-linear adjustment, \citealp{Beaumont2025}), ABC-RegN (our weighted local-linear adjustment) and ABC-NN (a local nonlinear regression adjustment using neural nets, \citealp{blumF2010non}). From the DGP in equation \eqref{eq:mix}, we simulate $n=100$ observations for the observed data $\y$. Just as in Example 1 in the main text, each ABC procedure is based on $N=25,000$ simulated draws, where the tolerance $\epsilon_n$ is chosen to be to the 1\% quantile of the overall simulated distances. For ABC-AR we use the Euclidean norm. Figure \ref{fig:gk1} plots the posteriors from the different ABC methods across the four parameters of the $g$-and-$k$ distribution for a representative experiment. 
\begin{figure}[H]
	\centering
	\setlength\figureheight{2cm} 
	\setlength\figurewidth{3cm} 
	\resizebox{\textwidth}{!}{\input{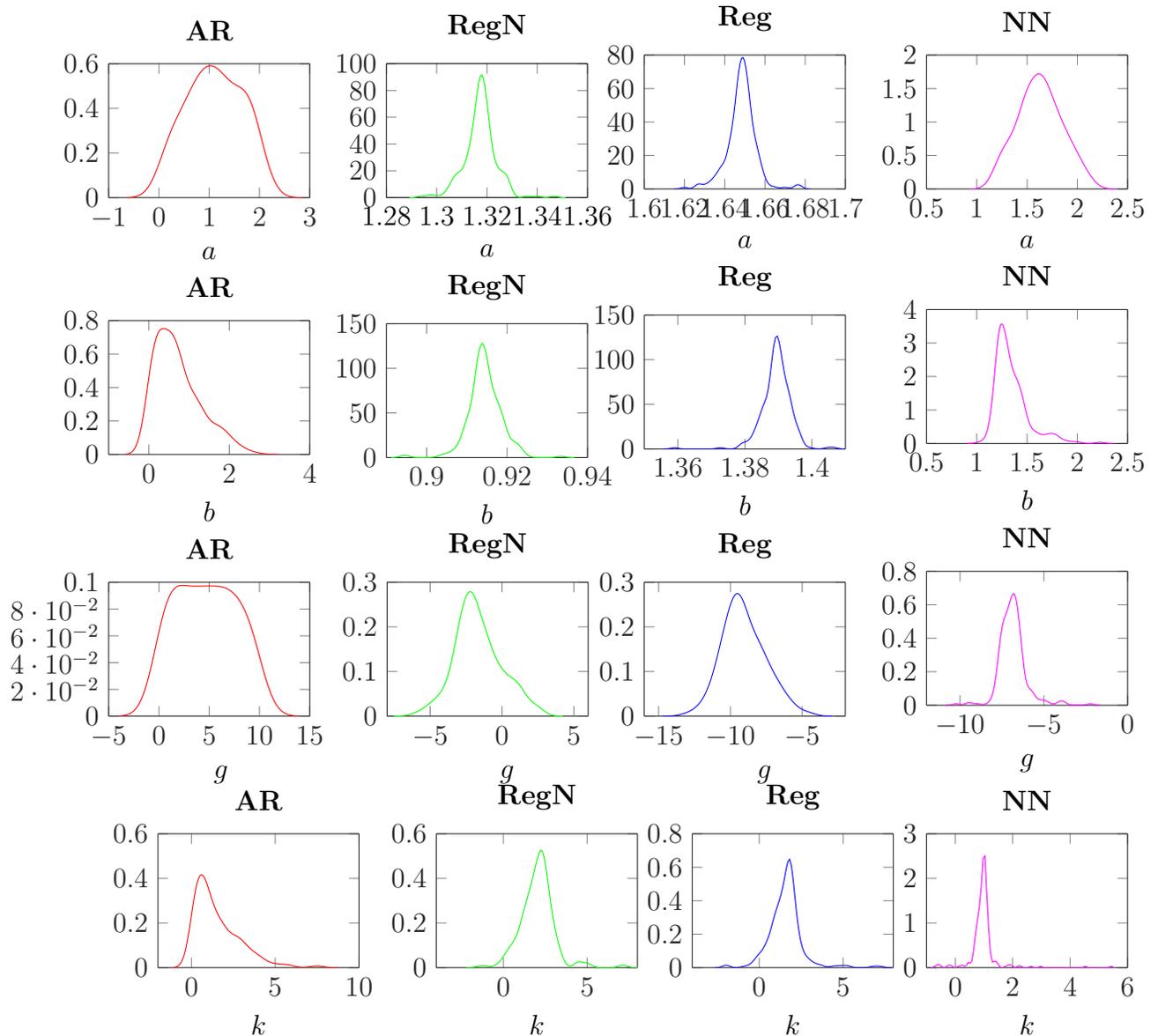}} 
	\caption{Posterior comparison of various ABC procedures in the $g$-and-$k$ model when the true DGP is given by equation \eqref{eq:mix}. } 
	\label{fig:gk1}
\end{figure}
Analyzing Figure \ref{fig:gk1}, three features are immediately in evidence. Firstly, for most of the ABC procedures, and across virtually all the parameters, the posteriors are decidedly non-Gaussian (as suggested by the results of Theorem 2 in the main text). Secondly, across the parameters, the different ABC procedures produce posteriors with very different behavior, which mirrors the results observed in Example 1. As mentioned earlier, this behavior is a direct consequence of the local regression adjustments transformation of the accepted ABC-AR draws. Thirdly, for the $g$ parameter the local regression posteriors place significant mass outside of the original prior space, which was $\mathcal{U}[0,10]$. Similarly, for the parameter $k$ these procedures yield some negative values for $k$. 

This later feature is troubling given the specific nature of the mismatch between the actual DGP and the assumed $g$-and-$k$ distribution. In particular, the theoretical moments in the actual DGP imply that population skewness is greater than 1.5 and that  excess kurtosis is nearly 3, respectively. Recalling that skewness and kurtosis in the $g$-and-$k$ distribution are controlled by $g$ and $k$, respectively, we know that $g>0$ is associated with positive skewness, and $k<0$ is associated with kurtosis that is less than that of the normal distribution. 
Given this, the results of Figure \ref{fig:gk1} demonstrate that the local regression adjustment procedures can place significant posterior mass on $g<0$, even though the observed data is always positively skewed. Similarly, these procedures can place some positive posterior mass on $k<0$, even though the observed data exhibits positive excess kurtosis. 

These behaviors observed in the local regression adjustments are entirely due to the fact that these approaches disregard the nature of the parameter space when adjusting the accepted draws, and, hence, can transport posterior mass outside the original parameter space. As a consequence, in this example the local regression adjustment approaches can place significant posterior mass on values of $g$ (and to a lesser extent $k$) that are incompatible with the observed data.

Similar to the conclusions obtained in Example 1, this behavior of the local regression adjustment is not a feature of any particular dataset but is persistent across different datasets. To demonstrate this fact, we simulate 1000 replications from the true DGP in equation \eqref{eq:mix}, and rerun each ABC procedure on these observed datasets. Across the replications, for each of the different ABC procedures, we record the posterior mean, standard deviation, as well as the length of the 95\% credible regions for each of the parameters in the $g$-and-$k$ distribution, which we calculate via the corresponding 2.5 and 97.5 quantiles of the individual parameter posteriors. In Table \ref{tab:gk1}, we present the averages of these results across the replications.

To help interpret the findings in Table \ref{tab:gk1}, before discussing the results, we first calculate the pseudo-true value $$\theta^*=\arg\min_{\theta\in\Theta}\|b_0-b(\theta)\|,$$ where we reminded the reader that $b_0$ is the probability limit of $\eta(\y)$, given by the octiles of the observed data, and $b(\theta)$ is the probability limit of $\eta(\mathbf{z})$, which corresponds to the octiles of the simulated data. To obtain this pseudo-true value, we must first calculate $b_0$, the population octiles of the Gaussian mixture distribution. The quantile function of the Gaussian mixture has no closed-form but the quantiles, and hence the value of $b_0$, can be obtained by numerically inverting the corresponding CDF of the Gaussian mixture. Give the value of $b_0$, and the fact that the quantiles of the $g$-and-$k$ distribution have an analytical form, in terms of the standard normal quantile function, we can numerically solve for the value of $\theta=(a,b,g,k)^{\intercal}$ that minimizes $\|b_0-b(\theta)\|$. Using this approach, the pseudo-true value under this particular Monte Carlo design is given by
\begin{equation}\label{eq:ptv}
\theta^*=(a^*,b^*,g^*,k^*)^\intercal=(1.17,1.50,0.41,0.23)^\intercal.
\end{equation}
This value of $\theta^*$ corresponds to a $g$-and-$k$ distribution with positive skewness and kurtosis that is larger than that of the normal distribution. Therefore, as one would hope, the pseudo-true value reflects the actual features of the true DGP, namely, positive skewness and excess kurtosis. Using this pseudo-true value, we also calculate the Monte Carlo coverage of the different ABC procedures, which we also display in Table \ref{tab:gk1}. 

\vspace{.3cm}

\begin{table}[H]
	\centering
	\caption{Posterior summaries for accept/reject ABC (AR), local linear regression adjusted ABC (Reg), see \cite{Beaumont2025}, our proposed local linear regression adjusted ABC (RegN), and a local nonlinear regression adjustment (NN) based on neural nets, see \cite{blumF2010non}. Mean is the average posterior mean, Std is the average posterior standard deviation, Len is the average length of the 95\% credible set, Cov is the Monte Carlo coverage calculated using the pseudo-true value defined in \eqref{eq:ptv}, and Q025 and Q975 refer to the average (across the Monte Carlo replications) 2.5\% and 97.5\% quantiles of the posterior.}
	\begin{tabular}{lrrrrrlrrrr}		\hline\hline
		
		$a $    & \multicolumn{1}{l}{AR} & \multicolumn{1}{l}{RegN} & \multicolumn{1}{l}{Reg} & \multicolumn{1}{l}{NN} &       & $b$     & \multicolumn{1}{l}{AR} & \multicolumn{1}{l}{RegN} & \multicolumn{1}{l}{Reg} & \multicolumn{1}{l}{NN} \\\hline
		Mean  & 1.1190 & 1.1568 & 1.2060 & 1.1862 &       & Mean  & 0.8961 & 1.0143 & 1.6342 & 1.6685 \\
		Std   & 0.6564 & 0.0027 & 0.0031 & 0.1333 &       & Std   & 0.7114 & 0.0104 & 0.0113 & 0.0994 \\
		Len   & 2.2517 & 0.0112 & 0.0125 & 0.4893 &       & Len   & 2.5380 & 0.0414 & 0.0439 & 0.3914 \\
		Cov   & 1.0000 & 0.0180 & 0.0140 & 0.5110  &       & Cov   & 0.9990 & 0.1730  & 0.5220  & 0.5600 \\
		Q025  & 0.0831 & 1.1509 & 1.1993 & 0.9419 &       & Q025  & 0.0478 & 0.9940 & 1.6119 & 1.5326 \\
		Q975  & 2.3349 & 1.1621 & 1.2118 & 1.4312 &       & Q975  & 2.5858 & 1.0354 & 1.6558 & 1.9240 \\
		&       &       &       &       &       &       &       &       &       &  \\		\hline\hline
		
		$g $    & \multicolumn{1}{l}{AR} & \multicolumn{1}{l}{RegN} & \multicolumn{1}{l}{Reg} & \multicolumn{1}{l}{NN} &       & $k$     & \multicolumn{1}{l}{AR} & \multicolumn{1}{l}{RegN} & \multicolumn{1}{l}{Reg} & \multicolumn{1}{l}{NN} \\\hline
		Mean  & 4.7758 & 5.3310 & 0.8639 & 0.3707 &       & Mean  & 1.5321 & 2.1790 & 0.7973 & 0.3889 \\
		Std   & 2.9404 & 1.5265 & 1.4632 & 0.7212 &       & Std   & 1.4145 & 0.9426 & 0.7810 & 0.3866 \\
		Len   & 9.4511 & 6.1010 & 5.8798 & 2.9902 &       & Len   & 5.1661 & 3.8949 & 3.2841 & 1.5546 \\
		Cov   &0.8900 &  0.3380  & 0.6320 & 0.6870 &       & Cov   & 1.0000&  0.5760 & 0.8420 & 0.8060 \\
		Q025  & 0.2484 & 2.4728 & -1.7755 & -0.9119 &       & Q025  & 0.0510 & 0.3648 & -0.6637 & -0.3487 \\
		Q975  & 9.6995 & 8.5738 & 4.1043 & 2.0783 &       & Q975  & 5.2171 & 4.2597 & 2.6204 & 1.2059 \\\hline
	\end{tabular}%
	\label{tab:gk1}%
\end{table}%

Analyzing the results in Table \ref{tab:gk1}, we see that all procedures give relatively accurate point estimates of the location parameter $a$. However, the Monte Carlo coverage for this parameter varies drastically across the different ABC procedures, which reflects the extremely small posterior standard deviations of the local regression adjustment approaches. For the scale parameter,  $b$, a similar story is in evidence. Namely, all procedures give point estimators that are not too far from the pseudo-true value, $b^*=1.5$, however, the small posterior standard deviations for the local adjustment procedures leads to significant under-coverage. 

For the parameter $g$, the different ABC procedures have very different posterior behavior. The most striking feature is that, while the average posterior means for the local linear adjustment (Reg) and the nonlinear adjustment (NN) are not too far from the pseudo-true value, $g^*=0.41$, both procedures place a significant amount of posterior mass on $g<0$.\footnote{This point can be seen by analyzing the $2.5\%$ and $97.5\%$ quantiles from these procedures.} We recall that the true DGP is such that the observed data always has positive skewness.

The results for the parameter $k$, which governs the kurtosis of the $g$-and-$k$ distribution, are similar to those obtained for the parameter $g$. Namely, while the local linear and nonlinear regression adjustments generally give point estimators that are closer to the pseudo-true value of $k^*\approx0.23$ than the other ABC procedures, both approaches place a non-negligible amount of posterior mass on $k<0$. This behavior is clearly at odds with the observed data: values of $k<0$ imply that the observed data exhibits kurtosis that is less than that of the normal. 

As demonstrated by Corollary 2 in the main paper, this behavior is a direct consequence of model misspecification and the nature of the local regression adjustment. In essence, under model misspecification, these adjustment procedures can be viewed as taking accepted draws, which asymptotically yield the smallest distance between observed and simulated summaries, and perturbing them according to a criterion that does not respect the original optimality of these draws. 

Given the behavior of the local regression adjustment in this example, as well as the results of Example 1 in the main text, we suggest researchers treat the output of local regression adjustment ABC procedures with a healthy level of skepticism in settings where model misspecification is a possibility.

\end{document}